\documentclass[fleqn,11pt]{article}
\usepackage{geometry}                % See geometry.pdf to learn the layout options. There are lots.
\usepackage{graphicx}
\usepackage{amssymb}
\usepackage{amsmath}
\usepackage{xfrac}
\usepackage{authblk}
\usepackage{setspace}
\usepackage[usenames,dvipsnames]{xcolor}
\usepackage{hyperref}
\usepackage{amsfonts}
\usepackage{pdfpages}
\usepackage{xcolor}
\usepackage{amsthm}
\newtheorem{theorem}{Theorem}[section]

\newtheorem{corollary}{Corollary}[theorem]
\newtheorem{assumption}{Assumption}
\newtheorem{remark}{Remark}

\usepackage{subcaption}
\usepackage{algorithm}
\usepackage{algpseudocode}

\usepackage{bigints}

\usepackage{array}
\newcolumntype{P}[1]{>{\centering\arraybackslash}p{#1}}
\newcolumntype{M}[1]{>{\centering\arraybackslash}m{#1}}

\DeclareMathOperator*{\argmin}{arg\,min}
\title{Optimal control for sampling the transition path process\\ and estimating rates
}

\author[1]{Jiaxin Yuan\thanks{jyuan98@umd.edu}}
\author[2]{Amar Shah\thanks{amarshah1000@berkeley.edu}}
\author[3]{Channing Bentz\thanks{cbentz2845@gmail.com}}
\author[1]{Maria Cameron\thanks{mariakc@umd.edu}}
\affil[1]{\small{Department of Mathematics, University of Maryland, College Park, MD 20742, USA}}
\affil[2]{\small{College of Letters and Science, University of California, Berkeley, CA 94709, USA}}                               
\affil[3]{\small{College of Arts and Sciences, University of Nebraska, Lincoln, NE 68508, USA}}
\begin{document}
\maketitle

\begin{abstract}
 Many processes in nature such as conformal changes in biomolecules and clusters of interacting particles, genetic switches, mechanical or electromechanical oscillators with added noise, and many others are modeled using stochastic differential equations with small white noise. The study of rare transitions between metastable states in such systems is of great interest and importance. The direct simulation of rare transitions is difficult due to long waiting times. Transition path theory is a mathematical framework for the quantitative description of rare events.  Its crucial component is the committor function, the solution to a boundary value problem for the backward Kolmogorov equation. The key fact exploited in this work is that the optimal controller constructed from the committor leads to the generation of transition trajectories exclusively. We prove this fact for a broad class of stochastic differential equations. Moreover, we demonstrate that the committor computed for a dimensionally reduced system and then lifted to the original phase space still allows us to construct an effective controller and estimate the transition rate with reasonable accuracy.  Furthermore, we propose an all-the-way-through scheme for computing the committor via neural networks, sampling the transition trajectories, and estimating the transition rate without meshing the space. We apply the proposed methodology to four test problems: the overdamped Langevin dynamics with Mueller's potential and the rugged Mueller potential in 10D, the noisy bistable Duffing oscillator,  and Lennard-Jones-7 in 2D.

\end{abstract}





%93E20\sep % optimal stochastic control
%60J20\sep % APPLICATIONS OF STOCHASTIC ANALYSIS
%65N99 % numerical methods for PDEs, BVPs

%\end{frontmatter}

\section{Introduction}
\label{sec:intro}
% \item what problem are we solving and molecular simulation, rare event simulation: difficult to characterize (look at Marzouk's paper)\\

The study of rare events in stochastic systems is crucial for understanding natural phenomena such %{ CB: replace "such natural phenomena" with "natural phenomena, such" } 
as conformal changes in biomolecules and clusters of interacting particles, protein folding, noise-driven transitions in nonlinear oscillator systems, genetic switches, and many others. Often rare events in such systems are associated with transitions between metastable states separated by high energetic barriers. Direct simulations of rare events are difficult due to long waiting times. 
Deterministic approaches based on solving partial differential equations are hampered by the high dimensionality of phase space or other numerical issues. In this work, we 
propose an approach for sampling transition trajectories between metastable states based on optimal control and use them to calculate transition rates. 
% {\color{blue} AS: Can we say why we are trying to calculate transition trajectories? I.e. describe what specifically are we trying to calculate ($\rho_{AB}, E[\tau_{AB}, \nu_{AB}$) and why this is useful} 
This approach is inspired by a remarkable fact demonstrated in a recent work by Zhang, Sahai, and Marzouk \cite{Zhang_Marzouk_2022} that a highly effective controller for a broad class of stochastic systems can be obtained using a rough approximation to the solution of an appropriate partial differential equation.
The theoretical foundation of this approach is motivated by work \cite{Lu2015} of Lu and Nolen and 
recent work \cite{Gao2021_OptimalControl} of Gao, Li, Li, and Liu detailing the case of overdamped Langevin dynamics.

\subsection{An overview}
Numerical methods for the study of rare events can be divided into two large classes: deterministic and stochastic. 

Stochastic methods include direct simulation suitable only for the case where the noise is relatively large \cite{BB2021} as well as various enhanced sampling algorithms. These include methods aimed at sampling rare transitions e.g. transition path sampling \cite{TPS_1997,TPS_1999,TPS_2002}, milestoning \cite{milestoning2020}, weighted ensemble \cite{Aristoff2022}, and { adaptive splitting methods \cite{AMS_cerou_guyader,AMS_review},} and methods for the exploration of the configurational space e.g. metadynamics \cite{metadynamics_2002,barducci2008well} and machine learning-assisted techniques -- see \cite{Wang2020,Belkacemi2021} and references therein.
%{  AS: maybe get rid of ``references therein'' since if it is important enough, we should directly cite it}

The class of deterministic methods can be subdivided into several categories.  First, there are methods aiming at finding maximum likelihood transition paths via a gradient descent in the path space \cite{NEB1998,string2002,MAM2004,stringCV2006,string2007,gMAM_Vanden_Eijnden} and via using control theory and a shooting approach \cite{beri2005,CCB2022}. While they are computationally cheap and suitable for high dimensions, these methods produce a single transition path along which the transition flux is focused as the noise amplitude tends to zero, and 
{ are only able to give an asymptotic scaling for the exponential factor of the transition rate in the limit of noise coefficient approaching zero.}
%are unable to give an estimate of the transition rate. 
Second, methods for computing the \emph{quasipotential} on a mesh are useful for visualization of the effective potential for two- or three-dimensional systems with nongradient drifts -- see \cite{yang2019computing,cameron2019computing,Paskal2022,Kikuchi2020RitzMF}.
% see e.g. { AS: replace `` -- see e.g.'' with ``: see''} \cite{yang2019computing,cameron2019computing,Paskal2022,Kikuchi2020RitzMF} and references therein. { AS: delete ``and references therein''}
Third, there are methods using the framework of transition path theory \cite{VE2006,VE2010} where the key component is the numerical solution of the \emph{committor problem}, a boundary-value problem for the stationary backward Kolmogorov equation. Novel techniques developed for accomplishing this task in dimensions higher, and even in some cases, much higher than three, include those based on training neural networks \cite{khooluying, LiLinRen2019,Rotskoff2020}, crafting diffusion maps \cite{Banisch2020,evans2021computing,Evans2022JCP}, or representing the solution by tensor trains \cite{TTLindsey}.

%Freidlin, M. Some Remarks on the Smoluchowski?Kramers Approximation. Journal of Statistical Physics 117, 617?634 (2004). https://doi.org/10.1007/s10955-004-2273-9

\subsection{Optimal control}
\label{sec:optimal_control}
 The first two categories of deterministic methods, i.e. methods for finding maximum likelihood transition paths in the vanishing noise limit 
 and methods for computing the quasipotential, are connected via the \emph{deterministic optimal control}. 
 For example, consider a system evolving according to an SDE of the form
\begin{equation}
\label{sde1}
dX_t = b(X_t)dt + \sigma(X_t)\sqrt{\epsilon}dW_t,
\end{equation}
where the drift field $b$ is smooth and has a finite number of attractors lying within a ball of a 
finite radius around the origin, the matrix function $\sigma$ is smooth, and $\epsilon$ is a small parameter.
Let $\mathcal{A}$ be an attractor of the corresponding ODE $\dot{x} = b(x)$.
The escape problem from the attractor $\mathcal{A}$ can be viewed as an optimal control problem where \emph{an optimal realization of the Brownian motion $W_t$ driving the process out of the basin of $\mathcal{A}$ is sought}. 
% { This is worded a little awkwardly. Maybe replace ``to drive'' with ``that drives''. Also why is this italicize?} 
Therefore, the controlled ODE with a controller $u$ is
\begin{equation}
\label{ode1}
\dot{x}= b(x) + \sigma(x)u,\quad x(0)\in\mathcal{A}.
\end{equation}
The cost functional is derived in the large deviations theory \cite{freidlin2004,FW2012}. If $\sigma$ is nonsingular everywhere, the cost functional defined for all absolutely continuous paths $\phi(\cdot)$ is given by 
\begin{equation}
\label{cost1}
S_T(\phi) = \frac{1}{2} \int_0^{T} \|u(\phi(t))\|^2dt \equiv  \frac{1}{2} \int_0^{T} \|\sigma^{-1}(\phi) [\dot{\phi}- b(\phi)]\|^2dt.
\end{equation}
The last expression is exactly the Freidlin-Wentzell action functional for SDE \eqref{sde1} \cite{FW2012}. 
The optimal controller $u$ is given by
\begin{equation}
\label{qpot}
u =\sigma^\top\nabla U,\quad{\rm where}\quad U(x) = \inf_{\phi,T}\left\{S_t(\phi)~|~\phi(0)\in\mathcal{A},\phi(T) = x\right\}
\end{equation}
is the \emph{quasipotential}. 
The infimum in (\ref{qpot}) taken over all paths and all final times is always achieved at $T=\infty$ since $\phi(0)\in\mathcal{A}$ \cite{FW2012}.
% { AS: I am a little confused about what this sentence is saying. What is achieved at $T = \infty$?}\jx{resolved}.
Plugging this optimal controller to ODE \eqref{ode1} results in the ODE governing the optimal escape path from the attractor $\mathcal{A}$ \cite{cameron2012,dahiya2018ordered1}
\begin{equation}
\label{ode2}
\dot{x}= b(x) + \sigma(x)\sigma(x)^\top \nabla U(x).
\end{equation}
In practice, if the quasipotential is found, one can find the optimal escape path from the basin of $\mathcal{A}$ by integrating ODE \eqref{ode2} backward in time starting at the point at the boundary of the basin of $\mathcal{A}$ where the quasipotential $U$ is minimal. 
Therefore, \emph{the quasipotential determines the optimal controller for finding the most probable escape path in the zero noise limit}, and the escape path is governed by ODE \eqref{ode2}.

A similar connection via \emph{stochastic optimal control} exists between transition path sampling and methods relying on solving the committor problem. 
% { AS: the ``,and transition path sampling'' reads kinda awkwardly. It might be better to say ``exists between transition path sampling and the third category ...'' } 
These two approaches address the case where the noise is small but finite. The study of stochastic optimal control problems in the context of molecular dynamics applications was started by Hartmann and collaborators in 2012 \cite{Hartmann2012}.
% { AS: Is the year super important here?}
In contrast to seeking an optimal realization of the Brownian motion as in the deterministic optimal control problem outlined above, the setting of the stochastic optimal control problem leaves the noise term unchanged. Instead, it aims at \emph{finding a minimal modification to the drift term of the governing SDE that would make all trajectories accomplish the desired transition}.  
% \jx{Maybe mention Lu and Nolen here: For a system evolving according to the process \eqref{sde1}, \cite{Lu2015} derived the SDE governing dynamics of transition trajectories. Gao et al. in \cite{Gao2021_OptimalControl} showed a similar result for overdamped Langevin equations from the perspective of control theory.}
For example, suppose a system is governed by the overdamped Langevin dynamics
\begin{equation}
\label{sde2}
dX_t = - \nabla V(X_t)dt + \sqrt{2\beta^{-1}}dW_t,
\end{equation} 
where $V(x)$ is a smooth and coercive potential with a finite number of isolated minima and { $\beta^{-1}$} 
is a small parameter often interpreted in chemical physics applications as the temperature.  
Let $A$ and $B$ be open disjoint sets surrounding local minima $x_A$ and $x_B$ of $V(x)$, and let $q(x)$ be the committor function,
% { AS: It might be worth defining the committor function more formally/generally since it is pretty important for the paper.} 
i.e., the probability that the process governed by \eqref{sde2} and starting at $x$ will first reach $B$ rather than $A$.
The \emph{committor function determines the optimal controller for stochastic dynamics}. 
% { AS: maybe its worth it to define (optimal) controller as well.} \jx{resolved}
Precisely, the dynamics of transition paths from $A$ to $B$ are governed by~\cite{Lu2015} 
\begin{equation}
\label{sde3}
% dX_t = \left[- \nabla V(X_t) + \frac{\beta}{2}\nabla\log q(x)\right]dt + \sqrt{2\beta^{-1}}dW_t, \quad X_0\in\partial A.
dX_t = \left[- \nabla V(X_t) + {2\beta^{-1}}\nabla\log q(x)\right]dt + \sqrt{2\beta^{-1}}dW_t, \quad X_0\in\partial A.
\end{equation} 
The function {$2\beta^{-1}\nabla\log q(x)$} is the optimal controller with respect to the cost functional~\cite{Hartmann2013,Gao2021_OptimalControl} 
\begin{equation}
\label{cost2}
C(u) = \mathbb{E}_P\left[\frac{1}{2}\int_0^{\tau_{AB}} \|u(X_t)\|^2dt   + g(X_{\tau_{AB}})\right],\quad {\rm where}\quad g(x) = \begin{cases}+\infty,& x\in \partial A\\0,&x\in\partial B\end{cases}
\end{equation}
is the exit cost, $\tau_{AB} = \inf\{t>0~|~X_t\in \bar{A}\cup\bar{B}\}$ is the stopping time, and $P$ is probability measure on the path space of SDE \eqref{sde2}.

% \jx{The following sentence is moved.}The SDE governing the dynamics of transition trajectories of the system evolving according to a more general SDE \eqref{sde1} was derived in \cite{Lu2015} { AS: I'm a little confused about what this sentence is saying. What was derived?} \jx{resolved}
% Somewhat different optimal control problems for SDE \eqref{sde1} were studied in \cite{hartmann2017,hartmann2019,Zhang_Marzouk_2022}.
% { AS: how were these control problems different?} 

\subsection{Applications}
In this work, we are especially interested in applications coming from chemical physics and mechanical engineering. 

In chemical physics, molecular motion is often modeled by Langevin dynamics in $\mathbb{R}^{6N}$ where $N$ is the number of atoms. 
To alleviate the problem of high dimensionality and make results more interpretable, 
physically motivated collective variables are often introduced. Collective variables (CVs) are functions of atomic coordinates effectively capturing the main dynamical modes.
% { AS: maybe say a little more about what collected variables are, somethign like ``parameters that capture some of the essential behavior of a complex system''}.  
Common choices of CVs are dihedral angles along the backbone of a studied biomolecule, interatomic distances between particular key atoms, etc.
The dynamics in collective variables  { $Z_t\in\mathbb{R}^d$, $d\ll 6N$, }are often modeled by the \emph{overdamped Langevin dynamics in collective variables}
\begin{equation}
\label{sde4}
{dZ_t =\left[-M(Z_t)\nabla F(Z_t) +\beta^{-1}\nabla\cdot M(Z_t)\right]dt + \sqrt{2\beta^{-1}}M^{1/2}(Z_t)dW_t}
% dX_t =\left[-M(X_t)\nabla F(X_t) +\beta^{-1}\nabla\cdot M(X_t)\right]dt + \sqrt{2\beta^{-1}}M^{1/2}(X_t)dW_t,
\end{equation}
where $M(Z_t)$ and $F(Z_t)$ are the diffusion tensor and the free energy found by standard techniques using molecular dynamics (MD) data \cite{stringCV2006}. { The computation of $M$ and $F$ is detailed in Appendix A of \cite{evans2021computing}.}
We note that the dynamics of SDE \eqref{sde4} do not necessarily accurately represent the dynamics of the CVs evaluated along the trajectory governed by the original SDE even if the original SDE is merely the overdamped Langevin dynamics \eqref{sde2} and with only one collective variable \cite{LL2010}. 
It is shown in \cite{LL2010} that \emph{the level sets of the collective variable should be normal to the manifold} along which the dynamics of the original system are focused for an accurate estimate of residence time near metastable states using the reduced system, i.e., SDE \eqref{sde4} with $X_t$ being a scalar function. 

Among mechanical engineering models, we are interested in nonlinear oscillators with small added noise:
\begin{equation}
\label{sde5}
\begin{cases}dX_t &= m^{-1}P_tdt\\
dP_t &= \left[-\gamma P_t - \nabla V(X_t)\right]dt + \sqrt{2\gamma m \epsilon} dW_t,\end{cases}
\end{equation}
where $\epsilon$ is the parameter regulating the noise amplitude, $\gamma$ is the friction coefficient, and $V$ is the potential energy function. 
%Note that SDE \eqref{sde5} can be rewritten in the form of SDE \eqref{sde1} with a \emph{singular diffusion matrix}. { a repetition}
%{ Maybe combine these two sentences: Note that both SDEs \eqref{sde4} and \eqref{sde5} are of the form of SDE \eqref{sde1} with a diffusion matrix function $\sigma(\cdot)$ being ${d\times r}$ where $d$ is the dimension of the phase space and $r\le d$ is the rank of $\sigma$, i.e., the matrix function $\sigma$ has linearly independent columns. In this case, $W_t$ is the standard $r$-dimensional Brownian motion.  }

Note that SDEs \eqref{sde4} and \eqref{sde5} are of the form of SDE \eqref{sde1} with a matrix function $\sigma(\cdot)$ being ${d\times r}$ where $d$ is the dimension of the phase space and $r\le d$ is the rank of $\sigma$, i.e., the matrix function $\sigma$ has linearly independent columns. In this case, $W_t$ is the standard $r$-dimensional Brownian motion.

%%%%%%%%%%%
\subsection{Goals and summary of main results}
\label{sec:goal}
%{ 
%Maybe we could change the title to 'summary of main results'}
The goal %{ contribution} 
of this work is two-fold.
The first objective is to establish the connection between the committor and the optimal control for a broad class of SDEs. The second objective is to develop a methodology based on optimal control and transition path theory 
for sampling transition paths from a metastable region $A$ to a metastable region 
$B$ and finding the transition rate from $A$ to $B$.
% {\color{blue} AS: This is very good, but I think you need to emphasize this earlier in the paper.}
%We are especially interested in applications coming from chemical physics and mechanical engineering. 

% { AS: you basically restate this last sentence in the next paragraph. Perhaps, find a way to combine them} \jx{I don't see it as repetition, previous case is a remark on the SDEs, next paragraph is an assumption}

%To achieve our goal, \textcolor{red}{AS: maybe restate the goal here} 
% we have accomplished the following steps.
 We develop a workflow that allows one to generate the transition trajectories and compute the transition rate without meshing the ambient space. Our Python codes implementing this workflow are posted on GitHub~\cite{margotyjx,mar1akc}. 
%The proposed methodology for sampling transition paths and computing transition rates consists of the following steps.

%{\bf 1. Set up and solve the optimal control problem.} 
{\bf 1. Theoretical result: the solution to the optimal control problem.} 
We have proven a theorem (Theorem \ref{thm:T1}) that established the relationship between the committor and the optimal control
% generalized the result proven in \cite{Gao2021_OptimalControl} 
for SDE \eqref{sde1} with $\sigma(\cdot)$ being ${d\times r}$, ${\sf rank}(\sigma) = r$. 
    The optimally controlled dynamics are found to be of the form
    \begin{equation}
    \label{sde6}
    dX_t = \left[ b(X_t) + \sigma \sigma^\top v^*(X_t)\right] dt + \sigma(X_t)dW_t,
    \end{equation}
    where the optimal control  $\sigma^\top v^*$ satisfies 
    \begin{equation}
    \label{null1}
    \sigma^\top v^* =\sigma^\top \nabla \log q^+
    \end{equation} 
    where $q^+$  is the forward committor  (\ref{eq: q+q- def}). 
 This is a generalization of { Theorem 3.3 in \cite{Gao2021_OptimalControl} and is related to the results in \cite{Lu2015} and \cite{Zhang_Marzouk_2022}.}
 
    % { AS: we need to define forward committor before this} \jx{resolved}
%    A similar result for a fixed stopping time optimal control problem is proven in \cite{Zhang_Marzouk_2022}.
    
{\bf 2. Compute the committor.} 
In MD applications, it is challenging to compute the committor $q(x)$ accurately for the dynamics of interest due to high dimensionality. Therefore, we compute it for the reduced dynamics in CVs \eqref{sde4} and then lift it to the original phase space assuming that the dynamics in it is overdamped Langevin \eqref{sde2}.
It is also difficult to obtain an accurate solution to the committor problem for SDE \eqref{sde5} due to the degeneracy of the elliptic PDE. Therefore, in both cases, we expect to have an approximate solution to the committor problem. 

%We advocate neural network-based solvers for the committor problem because they have several advantages. 
{ Neural network-based solvers for the committor problem have several advantages.} First, they find a globally defined smooth solution function whose gradient needed for the controller $v^*$ is readily accessible via automatic differentiation. Second, they do not require artificial boundary conditions on the outer boundary of the computational domain unlike finite difference and finite element methods. Finally, they do not require meshing the space which makes them more amenable to promotion to higher dimensions. Our neural network-based solver for the committor problem for SDE \eqref{sde4} is similar to the one by Li, Lin, and Ren \cite{LiLinRen2019} that exploits the variational formulation and sets up a solution model that automatically satisfies the boundary conditions. The committor problem for SDE \eqref{sde5} does not admit a variational formulation. Therefore, we use the PINNs framework by Raissi, Perdikaris, and Karniadakis \cite{Raissi2017PhysicsID}. 
    
   {\bf 3. Sample transition trajectories.}
    The transition path theory framework allows one to compute the transition rate once the committor is available. However, if the committor is inexact e.g., due to suboptimal or insufficient set of CVs, the transition rate determined in this way is likely to be highly inaccurate \cite{LL2010,Evans2022JCP}. On the other hand, % as it is shown in \cite{Zhang_Marzouk_2022}, 
    even a rough approximation to the solution of the backward Kolmogorov equation yields a very good controller. Therefore, we use the found committor to construct a controller according to \eqref{sde6}--\eqref{null1} and sample transition trajectories. 
    
{\bf 4. Estimate the transition rate.}
 The transition rate $\nu_{AB}$ is defined as the average number of transitions from $A$ to $B$ observed per unit time. We propose to estimate it as
   \begin{equation}
   \label{nuAB1}
   \nu_{AB} = \frac{\rho_{AB}}{\mathbb{E}[\tau_{AB}]},
   \end{equation}
   where $\rho_{AB}$ is the probability of a trajectory being reactive i.e. on its way from $A$ to $B$, and $\mathbb{E}[\tau_{AB}]$ is the \emph{expected crossover time} from $A$ to $B$ found by simulating the controlled process \eqref{sde6}--\eqref{null1}. 
   % \jx{I'm not sure if we need to clarify or define 'reactive' here as it is introduced later and seems obvious to our target readers. If necessary, can do 'i.e. from A to B without returning to A'.}
   The probability $\rho_{AB}$ is estimated using the computed committor.
  % {\color{gray} This approach gives estimates for the transition rate consistent with those found by brute force for the test problems studied in this work. On the other hand, the transition rate estimated directly using the computed committor is notably less accurate.}
       
{\bf 5. Validation.} 
We apply the proposed methodology to four test problems: 
the overdamped Langevin equation with Mueller's potential in 2D and { with the rugged Mueller potential in 10D as in \cite{LiLinRen2019}}, a single bistable Duffing oscillator as in \cite{Zhang_Marzouk_2022}, and Lennard-Jones cluster of 7 particles in 2D (LJ7) as in \cite{DicksonMakarov_2009,evans2021computing}. We assess the accuracy and demonstrate the efficacy of the proposed methodology. In all test cases, the estimates of the transition rate { by formula \eqref{nuAB1} } are consistent with those found by brute force { even if the committor is not very accurate as in the case of LJ7.} On the other hand, the transition rate estimated directly using the computed committor is
notably less accurate. An explanation for this phenomenon is offered. 
% \emph{ These test cases have shown that using formula \eqref{nuAB1} is robust to the estimation of the committor function. In particular, even if the committor estimate is imperfect, based on a reasonable estimation of $\rho_{AB}$, transition rates using \eqref{nuAB1} still attain high accuracy which is in the ballpark of ground truth transition rate.}

%This phenomenon is explained by the fact that the estimate for $\rho_{AB}$ integrates the committor while the formula for the transition rate integrates the gradient of the committor in which the numerical errors and the approximation errors are amplified. }

The paper is organized as follows. Section \ref{sec:background} provides the necessary background on transition path theory and the transition path process. Section \ref{sec:optimalcontrol} and \ref{sec:transition_rate} contains our theoretical results. Section \ref{sec:numerics} describes the numerical methods used in this work. Section \ref{sec:tests} presents the application of the proposed methodology to three benchmark test problems. Section \ref{sec:conclusion} summarizes the results and gives perspectives for future work. The proof of the main theorem (Theorem \ref{thm:T1}) as well as a number of technical aspects are elaborated in appendices.

%The background on the transition path theory, control theory, and problem formulation is in Section \ref{sec: Background}. Theoretical results concerning the optimal controller for generalized Ito's diffusion process are presented in Section \ref{sec: optimal control}. In Section \ref{sec: transition rate}, we introduce how the transition rate can be computed utilizing the information sampled from the controlled process. Section \ref{sec: numerics} presents how do we solve for committor function through neural network, and numerical results are shown in Section \ref{sec: experiments} . We conclude the paper with future work in Section \ref{sec: conclusion}.

%%%%%%%%%%%%%%%%%%%%%%%%%%%%%%%
%%%%%%%%%
%%%%%%%%%             B A C K G R O U N D
%%%%%%%%%
%%%%%%%%%%%%%%%%%%%%%%%%%%%%%%%

\section{Background} 
\label{sec:background}
In this section, we will provide the necessary theoretical background on transition path theory (TPT) and the transition path process.
%In this section, we first give a quick review of the dynamical system that we consider and then present transition path theory (TPT) and the committor function, together with statistics to characterize rare transitions with committors, i.e., reactive current and transition rate. Then we briefly introduce the importance sampling for rare events, from perspectives of optimal control \cite{Gao2021_OptimalControl}\cite{Hartmann2013}\cite{Hartmann_2016} and Koopman framework \cite{Zhang_Marzouk_2022}.

%%%%%%%%%%%%%%%%%%%%%%%
\subsection{Transition Path Theory}
% Assuming that function $f$ is twice continuously differentiable, the infinitesimal generator for the overdamped Langevin equation is defined to be 
% \begin{equation}\label{eq: generator}
%     Lf = b \cdot \nabla f + \frac{1}{2} \sigma \sigma^T : \nabla \nabla f = \beta^{-1} \Delta f - \nabla f \cdot \nabla U
% \end{equation}
Transition Path Theory (TPT) is a celebrated mathematical framework for the quantitative description of rare transitions in stochastic systems  \cite{VE2006,VE2010}. Suppose a system is evolving according to SDE
\begin{equation}
\label{sdeB1}
dX_t = b(X_t)dt + \sigma(X_t)dW_t, \quad X_t\in\Omega\subseteq\mathbb{R}^d, \quad x\in\Omega\subset\mathbb{R}^d.
\end{equation}
Throughout this work, we will adopt the following assumptions about SDE \eqref{sdeB1}.
\begin{assumption}
\label{A1}
The domain $\Omega$ is either $\mathbb{R}^d$ or a manifold without boundary with metric being locally Euclidean. 
\end{assumption}
For example, $\Omega$ can be  a $d$-dimensional ``flat" torus, $\mathbb{T}^d$, or a direct product $\mathbb{T}^k\times\mathbb{R}^{d-k}$, $1\le k\le d-1$.  
\begin{assumption}
\label{A2}
The drift field $b:\Omega\rightarrow\mathbb{R}^d$ is a smooth vector function. 
The corresponding ODE $\dot{x} = b(x)$ has a finite number of attractors and all its trajectories approach one of the attractors as $t\rightarrow\infty$.  
\end{assumption}
\begin{assumption}
\label{A3}
The matrix function $\sigma:\Omega\rightarrow \mathbb{R}^{d\times r}$, $r:={\sf rank}(\sigma)$, is smooth. 
The entries of $\sigma(x)$ are bounded, and ${\sf rank}(\sigma)$ is the same for all $x\in\Omega$. 
The singular values of $\sigma$ are bounded from above and from below in $\Omega$.
\end{assumption}
\begin{assumption}
\label{A4}
There exists a unique invariant density $\mu(x)$, and the system is ergodic. 
\end{assumption}
The infinitesimal generator $\mathcal{L}$ of the process governed by SDE \eqref{sdeB1} is defined as
\begin{equation}\label{eq: generator}
        \mathcal{L}f(x) = b \cdot \nabla f + \frac{1}{2}{\sf tr}\left( \sigma \sigma^T \nabla \nabla f\right).
\end{equation}
The invariant density is the solution to $\mathcal{L}^{\ast}\mu = 0$,  $\int_{\Omega}\mu dx = 1$, where $\mathcal{L}^{\ast}$ is the adjoint to the generator \eqref{eq: generator}.
% Recall that the invariant density is the solution to $\mathcal{L}^{\ast}\mu = 0$,  $\int_{\Omega}\mu dx = 1$, where $\mathcal{L}^{\ast}$ is the adjoint to the generator of the process \eqref{sdeB1} { AS: do we want to define generator or can it be assumed that people know it?}.

Suppose we want to study transitions between two disjoint open sets $A$ and $B$ in $\Omega$. 
For example, if we are interested in transitions between neighborhoods of two distinct attractors of the ODE $\dot{x} = b(x)$, we choose $A$ and $B$ to be these neighborhoods. Let $\{X_t~|~-\infty< t< \infty\}$ be an infinitely long trajectory of \eqref{sdeB1}. Due to ergodicity, it will visit $A$ and $B$ infinitely many times. TPT studies statistics of pieces of such a trajectory that start at $\partial A$ and next hit $\partial B$ without returning to $\bar{A}$\footnotemark[1] in between. \footnotetext[1]{The bar above a set denotes its closure.}
% { AS: I think it might be worth adding that $\delta A$ is the boundary of $A$ in the footnote}
Such pieces are called the \emph{reactive trajectories}, and $A$ and $B$ are called the \emph{reactant} and \emph{product sets} respectively. Key concepts of TPT are \emph{forward and backward committors} $q^{+}(x)$ and $q^{-}(x)$. The forward committor $q^{+}(x)$ is defined as the probability that the process starting at $x$ will first hit $\bar{B}$ rather than $\bar{A}$.
The backward committor is the probability that the process arriving at $x$ has hit $\bar{A}$ last rather than $\bar{B}$. Specifically, let $\tau^+_D$  ($\tau^-_D$)  be the first (last) hitting time of region $\bar{D}$ for trajectory starting (arriving) at $x$, and $\tau^-_D$  be the last hitting time of region $\bar{D}$ for trajectory arriving at $x$. Equivalently, $\tau^-_D$ is the first hitting time of $\bar{D}$ for the time-reversed process of $\hat{X}_t = X_{-t}$:
\begin{equation}
    \begin{cases}
        \tau^+_D(x) = \inf\{t \geq 0: X_t \in \bar{D}, x(0) = x\}\\
        \tau^-_D(x) = \inf\{t \geq 0: \hat{X}_t \in \bar{D}, \hat{x}(0) = x\}.
    \end{cases}
\end{equation}
Given two disjoint regions $A$ and $B$, the forward committor function $q^+(x)$ and backward committor function $q^-(x)$ are defined as follows
\begin{equation}\label{eq: q+q- def}
\begin{cases}
    q^+: \Omega \rightarrow [0,1], q^+(x) = \mathbb{P}\{\tau^+_B(x) < \tau^+_A(x)\}\\
    q^-: \Omega \rightarrow [0,1], q^-(x) = \mathbb{P} \{\tau_A^-(x) < \tau_B^-(x)\}.
\end{cases}
\end{equation}The region $\Omega$ with removed sets $\bar{A}$ and $\bar{B}$ will be denoted by $\Omega_{AB}$:
\begin{equation}
\label{OAB}
\Omega_{AB}:= \Omega\backslash(\bar{A}\cup\bar{B}).
\end{equation}
The forward and backward committors $q^+$ and $q^-$ are the solutions to the following boundary value problems (BVPs) \cite{VE2006}
\begin{equation}
\label{eg:commproblem}
\begin{cases}
\mathcal{L}q^+ = 0 & x \in \Omega_{AB}\\
q^+(x) = 0 & x \in \partial A\\
q^+(x) = 1 & x \in \partial B,\\
\frac{\partial q^+}{\partial \hat{n}} = 0,& x\in\partial\Omega,
\end{cases}\hfill \quad
\begin{cases}
    \mathcal{L}^{\dagger}q^- = 0 & x \in \Omega_{AB}\\
q^-(x) = 1 & x \in \partial A\\
q^-(x) = 0 & x \in \partial B,\\
\frac{\partial q^-}{\partial \hat{n}} = 0,& x\in\partial\Omega.
\end{cases}
\end{equation}
We will refer to the BVP for the forward committor as the \emph{committor problem}.
Here, $\mathcal{L}$ is the infinitesimal generator \eqref{eq: generator}.
% of the process governed by SDE \eqref{sdeB1}{ AS: We should have this definition earlier since we use the generator at the start of this section}
$\mathcal{L}^{\dagger}$ is the infinitesimal generator of the corresponding time-reversed process~\cite{VE2010} 
\begin{equation}
    \mathcal{L}^\dagger f(x) = -b \cdot \nabla f + \frac{1}{\mu}{\sf div}( \sigma \sigma^T \mu) +\frac{1}{2} {\sf tr}\left( \sigma \sigma^T \nabla \nabla f\right),
\end{equation} 
where ${\sf div}( \sigma \sigma^T \mu)$ is the divergence of the matrix $\sigma \sigma^T \mu$, i.e., a vector with components $\left[{\sf div}( \sigma \sigma^T \mu)\right]_i =\sum_ j\partial_{x_j}( \sigma \sigma^T \mu)_{ji}$, $1\le i\le d$.
The homogeneous Neumann boundary conditions on $\partial\Omega$ in \eqref{eg:commproblem} are relevant if $\Omega$ 
has a reflecting boundary; $\hat{n}$ is the external unit normal to $\partial \Omega$.
Note that for both SDEs of our interest, \eqref{sde4} and \eqref{sde5}, the invariant density  $\mu$ is known.
For overdamped Langevin dynamics in CVs \eqref{sde4},
\begin{equation}
 \mu(x) = Z_Fe^{-\beta F(x)},~~Z_F = \int_\Omega e^{-\beta F(x)}dx\quad\text{for SDE \eqref{sde4}}\label{mu4}
\end{equation}
and for Langevin dynamics \eqref{sde5},
\begin{align}
 \mu(x,p)& = Z_He^{- H(x,p)/\epsilon},\quad H(x,p) = \frac{p^2}{2m}+V(x), \label{mu5}\\
  Z_H &= \int_\Omega e^{-H(x,p)/\epsilon}dxdp\quad\text{for SDE \eqref{sde5}}.\notag
\end{align}
The generators for SDEs \eqref{sde4} and \eqref{sde5}, respectively, are 
% {\color{blue} AS: We reference SDE (9) and (10) enough that it might be better just to refer to them by name (e.g. ``overdamped Langevin dynamics'' and ``nonlinear oscillators'') instead, so that the reader doesnt have to click to see what they are each time. This is just a stylistic choice though. Feel free to ignore this.}
\begin{equation}
\label{gen4}
\mathcal{L}f = \beta^{-1}e^{\beta F}\nabla\cdot\left(e^{-\beta F}M\nabla f\right)
\end{equation}
and
\begin{equation}
\label{gen5}
\mathcal{L}f = \frac{p}{m}\cdot\nabla_xf -\nabla_x V\cdot\nabla_p f - \gamma p\cdot\nabla_p f +\gamma m\epsilon\Delta_p f,
\end{equation}
where $\nabla_x$ and $\nabla_p$ denote gradients with respect to the coordinates and momenta respectively, and $\Delta_p$ is the Laplacian with respect to momenta.

If the governing SDE is time-reversible, e.g. the overdamped Langevin dynamics \eqref{sde2} and the overdamped Langevin dynamics in collective variables \eqref{sde4},
% as the overdamped Langevin dynamics and SDE \eqref{sde4}, the overdamped Langevin dynamics in collective variables, { CB: the phrasing of this is very awkward. Can it be left just as "If the governing SDE is time-reversible," ? If not maybe replace it with "If the governing SDE is time-reversible in the form of overdamped Langevin dynamics and SDE[9],"}
the backward committor is readily found from the forward committor: $q^{-}(x) = 1-q^{+}(x)$. Langevin dynamics \eqref{sde5}
% SDE \eqref{sde5} 
is not time-reversible, however, there is a nice relationship between the forward and backward committors \cite{VE2010}: $q^{-}(x,p) = 1 - q^{+}(x,-p)$.

The time-reversible SDE \eqref{sde4} admits a variational formulation for the committor problem \cite{stringCV2006} that motivates the construction of the loss function for neural network-based committor solvers \cite{khooluying, LiLinRen2019}:
\begin{equation}
\label{eg:variational}
    q^+(x) = \argmin_{f\in \mathcal{Q}} \int_{\Omega_{AB}} \nabla f(x)^\top M(x)\nabla f(x)\mu(x)dx,
\end{equation}
where $\mu$ is the invariant density given by \eqref{mu4} and $\mathcal{Q}$ is the set of all continuously differentiable functions $f$ satisfying the boundary conditions 
$f(x) = 0$, for $x \in \partial A$, and $f(x) = 1$, for $x \in \partial B$.

The probability density of reactive trajectories is given by
\begin{equation}
\label{muAB}
\mu_{AB} = \mu q^{+}q^{-}.
\end{equation}
The integral of $\mu_{AB}$ over $\Omega_{AB}$,
\begin{equation}
\label{rhoAB}
\rho_{AB} =\int_{\Omega_{AB}} \mu q^{+}q^{-}dx
\end{equation}
 is the probability that a stochastic trajectory at a randomly picked time is reactive, i.e., is on its way from $A$ to $B$. 
 % { AS: can we give a name to this since we are going to use it later (when talking about expected crossover time)?}

The transition rate from $A$ to $B$, $\nu_{AB}$, is defined as
\begin{equation}
\label{nuAB}
\nu_{AB} = \lim_{T\rightarrow\infty}\frac{N_{AB}}{T},
\end{equation} 
where $N_{AB}$ is the number of transitions from $A$ to $B$ observed during the time interval $[0,T]$.
The escape rate from $A$ is defined as
 \begin{equation}
\label{kAB}
k_{AB} = \lim_{T\rightarrow\infty}\frac{N_{AB}}{T_A}\equiv \frac{\nu_{AB}}{\rho_A},
\end{equation} 
where $T_A$ is the total time within $[0,T]$ during which the system last hit $A$ and 
$\rho_A$ is the probability that an infinitely long trajectory at any randomly picked moment of time has last hit $\bar{A}$. 
% rather than $\bar{b}B$.{ CB: typo? should just be $\bar{B}?$} \jx{References/appendix?}
The probability $\rho_A$ is equal to
\begin{equation}
\label{rhoA}
\rho_A = \int_{\Omega}\mu(x)q^-(x)dx.
\end{equation}

The reactive current is a vector field such that its flux through any surface separating $A$ and $B$ is the transition rate $\nu_{AB}$. 
For SDEs \eqref{sde4} and \eqref{sde5}, it is given, respectively, by
% { AS: maybe its worth stating the general formula, so readers know where these came from?}
\begin{equation}
\label{rcur4}
    {J}_{AB} =  \beta^{-1}Z^{-1} e^{-\beta F}M \nabla q^+
\end{equation}
and
\begin{equation}
\label{rcur5}
    {J}_{AB} = Z^{-1}_H e^{- H/\epsilon}q^+q^-\begin{pmatrix}
p\\
-\nabla_x V
    \end{pmatrix} + \epsilon \gamma Z^{-1}_H e^{- H/\epsilon}\begin{pmatrix}
        0\\
        q^- \frac{\partial q^+}{\partial p} - q^+ \frac{\partial q^-}{\partial p}
    \end{pmatrix}.
\end{equation}
A natural choice of a surface separating $A$ and $B$ is an isocommittor surface $\Sigma_\alpha:=\{x\in\Omega~|~q^+(x) = \alpha\}$ for any $\alpha\in[0,1]$. 
The transition rate can be expressed as (Proposition 6 in \cite{VE2006}):
% calculated by averaging the reactive flux through all isocommittor surfaces, resulting in the following result for SDE \eqref{sdeB1}:
\begin{equation}
\label{nu1}
\nu_{AB} = \frac{1}{2}\int_{\Omega_{AB}}\mu[\nabla q^+]^\top\sigma\sigma^\top\nabla q^+ dx = \frac{1}{2}\int_{\Omega_{AB}}\mu[\nabla q^-]^\top\sigma\sigma^\top\nabla q^- dx.
\end{equation}
% { AS: I'm a little confused about this equation. It's not super clear to me why this intgral is ``averaging the reactive flux through all isocommittor surfaces.'' Also, I am not sure about why the second equality holds: we know that $\nabla q^+(x) = \nabla q^-(x)$ for time-reversible equation and for SDE (10) from what you said earlier, but we don't know this for the general case right?} \jx{resolved in \cite{VE2006}}
In particular,
for SDEs \eqref{sde4}  and \eqref{sde5}, the transition rates, respectively, are
\begin{equation}
\label{nu4}
    \nu_{AB} = \beta^{-1} Z_F^{-1} \int_{\Omega_{AB} }(\nabla q^+)^\top M \nabla q^+ e^{-\beta F}dx
\end{equation}
and
\begin{equation}
\label{nu5}
     \nu_{AB} = \epsilon \gamma Z^{-1}_H \int_{\Omega_{AB}} \left[\nabla_p q^+\right]^\top m\left[ \nabla_p q^+ \right]e^{- H/\epsilon}dxdp,
\end{equation}
where $m$ is the diagonal mass matrix.

%%%%%%%%%%%%%%%%%%%%%%%%%%%%%%%
\subsection{Transition path process}
\label{sec:TPP}
TPT received an interesting development in the paper ``Reactive trajectories and the transition path process" by Lu and Nolen (2015) \cite{Lu2015}.
% \jx{Move these sentences}The results in this paper are important for the present work. While the dynamics in assumed to be governed by SDE \eqref{sdeB1} with the matrix $\sigma\sigma^\top$ satisfying the strong ellipticity condition, { most results in this paper can be easily extended to the class of processes considered in this paper} (CHECK THIS), i.e., evolving according to SDE \eqref{sdeB1} with Assumptions \ref{A1}--\ref{A4}. { AS: I don't think this paragraph adds much. It would be better to directly state the results from the Lu/Nolen paper and then say that they can be extended to the processes in this paper.}
It is proven in \cite{Lu2015} that the dynamics of the reactive trajectories for SDE \eqref{sdeB1}, with matrix $\sigma \sigma^\top$ satisfying the strong ellipticity condition,  are governed by the SDE
% { AS: as a clarification, is this sentence just saying that if we take trajectories given by equation (14) and we take only the reactive trajectories, this is the same as taken trajectories given by equation (36)? If so, that makes sense, but then we should say something about equation (14) here}
\begin{equation}
\label{sdeB1c}
        d{X}_t = \left[b({X}_t) + \sigma({X}_t) \sigma^\top({X}_t) \nabla\log q^+({X}_t)\right] dt + \sigma({X}_t) dW_t,
\end{equation}
and this expression was suggested by the Doob $h$-transform. 
% { The reactive trajectories appear at the boundary $\partial A$ and disappear at the boundary $\partial B$ at the same rate. One can imagine that the reactive trajectories are teleported from $\partial B$ to $\partial A$ as illustrated in Fig. \ref{fig:tpp}.}
The process governed by \eqref{sdeB1c} is called \emph{the transition path process}.
Though we do not assume strong ellipticity condition for $\sigma \sigma^\top$, most results in \cite{Lu2015} can be extended to the processes considered in this paper i.e., evolving according to SDE \eqref{sdeB1} with Assumptions \ref{A1}--\ref{A4}.

% Furthermore, 
The equilibrium unnormalized density of reactive trajectories on $\partial A$ is 
\begin{equation}
\label{distrA}
\eta_A = \mu(x)\mathcal{L}q^+(x) = - \mu(x)\hat{n}(x)^\top\sigma(x)\sigma^\top(x)\nabla q^+(x),\quad x\in\partial A,
\end{equation}
where $\hat{n}(x)$ is the unit normal pointing inside $A$. 
Note that while $\mathcal{L}q^+ =0$ in $\Omega_{AB} \equiv \Omega\backslash (\bar{A}\cup\bar{B})$, 
it does not need to be zero on $\partial A$. Equation \eqref{distrA} should be understood as the rate at which transition paths exit $A$.
The equilibrium unnormalized density of reactive trajectories on $\partial B$ is
\begin{equation}
\label{distrB}
\eta_B = -\mu(x)\hat{n}(x)^\top\sigma(x)\sigma^\top(x)\nabla q^-(x),\quad x\in\partial B,
\end{equation}
where $\hat{n}(x)$ is the unit normal pointing inside $B$.

The expected crossover time $\mathbb{E}[\tau_{AB}]$ of the reactive trajectories is defined as 
\begin{equation}
\label{cross}
\mathbb{E}[\tau_{AB}] = 
\lim_{N_{AB}\rightarrow\infty}
\frac{1}{N_{AB}}
\sum_{k=0}^{N_{AB}-1}
\left(\tau_{B,k}^+ - \tau_{A,k}^-\right),
\end{equation}
where $N_{AB}$ is the number of transitions from $A$ to $B$ and $\tau_{B,k}^+ $ and $ \tau_{A,k}^-$ are $k$th entrance time to $B$ and $k$th exit time from $A$, 
respectively, registered for a long trajectory of SDE \eqref{sdeB1} initiated at $x\in A$.
The \emph{expected crossover time} is the ratio of the probability $\rho_{AB}$ \eqref{rhoAB}
% that an infinitely long trajectory is reactive at a randomly picked moment of time { AS: we should have a shorter name for this} 
and the transition rate \cite{Lu2015} (also see \ref{sec:app0})
\begin{equation}
\label{cross1}
\mathbb{E}[\tau_{AB}] = \frac{\rho_{AB}}{\nu_{AB}}.
\end{equation}
This relationship is very important for this work as it will be used for restoring the transition rate using sampled transition trajectories (see Section \ref{sec:transition_rate}). 

The transition path process governed by SDE \eqref{sdeB1c} can be thought of as a process in which the trajectories are killed as they reach $\partial B$ 
and reintroduced at $\partial A$ after a waiting time at the rate $\nu_{AB}$. Moreover, the trajectories entering $\Omega_{AB}$ 
through the boundary of $A$ are distributed according to \eqref{distrA}. 

The probability density of reactive trajectories $\mu_{AB} $ is an invariant measure of the transition path process \cite{VE2006,Lu2015}.
Indeed, the backward and forward Kolmogorov operators for \eqref{sdeB1c} are, respectively,
\begin{align}
\mathcal{L}_{AB}f &= \left[b + \sigma\sigma^\top \nabla \log q^+\right]\cdot\nabla f + \frac{1}{2}{\sf tr}\left(\sigma\sigma^\top \nabla\nabla f\right)\quad{\rm and}\label{genAB}\\
\mathcal{L}^*_{AB}f &= -\nabla\cdot\left[(b + \sigma\sigma^\top \nabla \log q^+)f - \frac{1}{2}{\sf div}\left(\sigma\sigma^\top f\right)\right].\label{gen*AB}
\end{align}
One can check that 
\begin{align}
\mathcal{L}^*_{AB}\mu_{AB} &= -\nabla\cdot \left[(b + \sigma\sigma^\top \nabla \log q^+)\mu q^+q^- - \frac{1}{2}{\sf div}\left(\sigma\sigma^\top \mu q^+q^-\right)\right] \notag \\
& = -\nabla\cdot\left[ \left(b\mu -\frac{1}{2}{\sf div}\left(\sigma\sigma^\top \mu \right)\right)q^+q^- 
+\frac{1}{2}\mu\sigma\sigma^\top\left(q^-\nabla q^+  - q^+\nabla q^-\right)\right] \label{JAB} \\
& = -\nabla\cdot J_{AB} = 0.
\end{align}

% { AS: should this be $\mathcal{L}_{AB}^*\mu_{AB}$ on the first line?}

\begin{figure}[htbp]
\begin{center}
\includegraphics[width=0.5\textwidth]{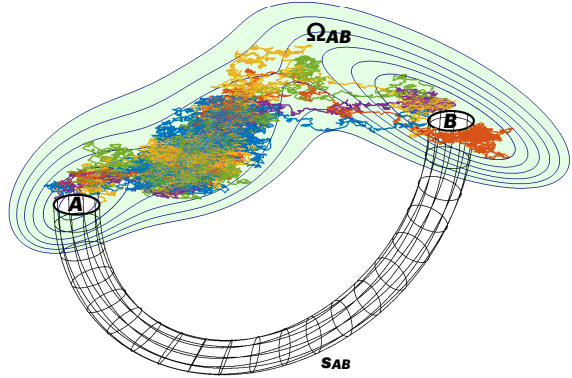}
\caption{An illustration for the transition path process. { Transition trajectories appear at the boundary of the region $A$, travel to the region $B$ without returning to $A$, and disappear at the boundary of $B$. To complete the dynamics of transition trajectories to an equilibrium process, the trajectories absorbed at $\partial B$ must be transported back to $\partial A$. The pipe in the figure symbolizes this reverse transport.}}
\label{fig:tpp}
\end{center}
\end{figure}

The symbol $J_{AB}$ is the reactive current \cite{VE2010}.
A calculation showing that $-\nabla\cdot J_{AB} = 0$ in $\Omega_{AB}$ is detailed in  \ref{sec:appA}.  
In order to make the transition path process an equilibrium process
we { assume that the transition trajectories that have reached $\partial B$ are transported back to $\partial A$ as shown in Fig. \ref{fig:tpp}.
A similar construction was used for Markov jump processes in \cite{CVE2014}.}
%This construction is illustrated in Fig. \ref{fig:tpp}. 

The probability to find the system at $s_{AB}$ at a randomly picked moment of time is $1-\rho_{AB}$. 
As a result, the invariant measure of the transition path process becomes
\begin{equation}
\label{muAB1}
\mu_{AB} = \begin{cases} \mu q^+q^-,& x\in\Omega_{AB}\\
1-\rho_{AB},&x= s_{AB}\end{cases}.
\end{equation}
It follows from \eqref{JAB} and \eqref{eg:commproblem} that the reactive current $J_{AB}$ on the boundaries $\partial A$ and $\partial B$ of $s_{AB}$ is \cite{Lu2015} 
% {  AS: Can we have this citation earlier in the sentence. This just reads a little awkwardly}
\begin{equation}
\label{JABbdry}
J_{AB}=\begin{cases} \frac{1}{2}\mu\sigma\sigma^\top\nabla q^+,& x\in\partial A\\
- \frac{1}{2}\mu\sigma\sigma^\top\nabla q^-,& x\in\partial B
\end{cases}.
\end{equation}
Since $\nabla\cdot J_{AB} = 0$ in $\Omega_{AB}$ we have 
$$
\int_{\partial A} J_{AB}\cdot\hat{n} ds  + \int_{\partial B} J_{AB}\cdot\hat{n} ds = 0
$$
where $\hat{n}$ is the outer unit normal and $ds$ is the surface element. Equation \eqref{JABbdry} is consistent with \eqref{distrA} and \eqref{distrB}.

% In this work, we later show that the transition rate $\nu_{AB}$ can be estimated by an alternative method, with the aid of the optimal controller.
%%%%%%%%%%%%%%%%%%%%%%%%%%%%%%%%%%%%%%%%%%%%%%%%%%%%%%%%%
 %%%%%%%%%%%%%%%%%%%%%%%%%%%%%%%%%%%%%%%%%%%%%%%%%%%%%%%%%
 %%%%%%%%%%%%%%%%%%%%%%%%%%%%%%%%%%%%%%%%%%%%%%%%%%%%%%%%%

\section{Optimally controlled dynamics}
\label{sec:optimalcontrol}
In this section, we show that SDE \eqref{sdeB1c} governing the transition path process can be obtained as the solution to a stochastic optimal control problem.
Precisely, we generalize Theorem 3.3 {  proven in~\cite{Gao2021_OptimalControl} for the case of the overdamped Langevin dynamics.
}
%addressing the case of the overdamped Langevin dynamics.
% in~\cite{Gao2021_OptimalControl} proven for the case of the overdamped Langevin dynamics.

\subsection{The general Ito SDE}
We consider the dynamics governed by the general Ito SDE \eqref{sdeB1} where the drift field $b$ and 
the diffusion matrix $\sigma$ are as described at the beginning of Section \ref{sec:background}. 
The controlled dynamics are set to be of the form
\begin{equation}
\label{eq:sde1c}
        dY_t = \left[b(Y_t) + \sigma(Y_t) \sigma^\top(Y_t) v(Y_t)\right] dt + \sigma(Y_t) dW_t
\end{equation}
where $v:[0,\infty)\rightarrow\Omega\subset\mathbb{R}^d$ is to be chosen in an optimal manner. 
{ The form of the modification to the drift, $\sigma\sigma^{\top} v$, is borrowed from \cite{Zhang_Marzouk_2022} and is motivated as follows. We are seeking a drift term that will make all trajectories of the resulting process obey the same statistics as the transition trajectories of the original process \eqref{sdeB1}. This means that all trajectories of the controlled process can be observed in the original process at particular noise realizations. Hence, the span of the modification to the drift at every point should match the span of the noise term which is the column space of the matrix-valued function $\sigma(\cdot)$. That's why there is the factor $\sigma$ on the left of $\sigma\sigma^{\top} v$. The factor $\sigma^\top$ is used for convenience so that the control $v$ is of the same dimension as the process $X_t$.}

We have chosen the cost functional to be 
% {  AS: maybe we should talk about why this makes sense to be the cost functional}
\begin{equation}
\label{eq:costfunction}
    C_x[v(\cdot)] = \mathbb{E}_P\left[\frac{1}{2} \int_0^{\tau_{AB} }\|\sigma^\top(Y_s) v(Y_s)\|^2 ds + g(Y_\tau)~|~X_0=x\right],
\end{equation}
where $P$ is the probability measure on the path space of SDE \eqref{sdeB1}, 
\begin{equation}
\label{tauAB}
\tau _{AB}= \inf\{t>0~|~Y_t\in(\bar{A}\cup\bar{B})\}
\end{equation}
is the stopping time,
and $g$ is the exit cost defined by 
\begin{equation}
\label{eq:exitcost}
        g(x) = 
        \begin{cases}
            +\infty, & x \in \overline{A} \\
            0, & x \in \overline{B}
        \end{cases}.
\end{equation}
Cost functional \eqref{eq:costfunction} gives finite cost only if the trajectory leaves $\Omega_{AB}$ via the boundary of $B$. The function $\sigma^\top v$ is the standard form in which the optimal control is sought~\cite{Fleming2006}. 
The optimal control problem is to find the function $v(\cdot)$ that minimizes the cost functional \eqref{eq:costfunction}. 
Its solution is given in the following theorem.
\begin{theorem} 
\label{thm:T1}
Let $X_t$ be a process governed by SDE \eqref{sdeB1} satisfying Assumptions \ref{A1}--\ref{A4} and let $Y_t$ be the corresponding controlled process governed by SDE \eqref{eq:sde1c}.
In addition, we assume that $\Omega$ is compact with a reflecting boundary and $\partial A$ and $\partial B$ are smooth.
 The infimum of the cost functional \eqref{eq:costfunction} is given by
 \begin{equation}
 \label{mincost}
 c^*(x):=\inf_{v\in\mathcal{V}}C_x[v(\cdot)]  = -\log q^+(x),
 \end{equation}
 where $\mathcal{V}$ is the set of admissible controls 
 \begin{equation}
 \label{admiss}
 \mathcal{V}:=\overline{\left\{v\in C^1(\Omega_{AB})~|~\mathbb{E}_P\left[\exp\left(\int_0^{\tau{AB}}\frac{1}{2}\|\sigma^T(Y_s) v(Y_s)\|^2ds\right)\right]<\infty\right\}},
 \end{equation} 
$P$ is the measure on the path space of \eqref{sdeB1}, and $q^+$ is the forward committor for SDE \eqref{sdeB1}.
 The corresponding optimal control $v^*$ satisfies
     \begin{equation}
    \label{vstar}
    \sigma^\top v^*  = \sigma^\top \nabla \log q^+.
    \end{equation} 
\end{theorem}
The proof of Theorem \ref{thm:T1} combines ideas from Gao et al. (\cite{Gao2021_OptimalControl}, the proof of Theorem 3.3) 
and L.~C. Evans's notes on the control theory \cite{LCEvans_notes}. 
% {  AS: I don't think we need to mention Evan's Notes on control theory by name since we only use basic control theory facts/ideas from there. We can still cite it in the proof. I think we could maybe just delete this sentence}.
It is found in  \ref{sec:AppB}. 

\begin{remark}
\label{remark1}
Theorem \ref{thm:T1} shows that the optimally controlled dynamics are governed by SDE \eqref{sdeB1c} for the transition path process.
Therefore, all facts about the transition path process stated in Section \eqref{sec:TPP} are valid for the optimally controlled dynamics.
\end{remark}
% {  AS: Doesn't Section 2.2 use that Equation 36 governs the dynamics of reactive trajectories, not necessarily optimally controlled reactive trajectories. Why are these the same?} \jx{resolved}

\begin{remark}
The requirement that the domain $\bar{\Omega}$ is compact is often implemented in numerical simulations: 
the computational domain is always bounded. Deterministic techniques always use meshes or point clouds of finite size. 
In stochastic simulations, particles are often put in a box. 
Therefore, the assumption that $\Omega$ is compact is not practically restrictive.
\end{remark}
{ \begin{remark}
Theorem \ref{thm:T1} suggests a new approach for finding the committor via minimization of the cost functional \eqref{mincost}. We leave the investigation into the viability of this approach for the future.
\end{remark}
}

%%%%%%%%%%%%%%%%%%%%%%%%%%%%%%%%%%%%%%%%%%%%%%%%%%%%%%%%%%
%%%%%%%%%%%%%%%%%%%%%%%%%%%%%%%%%%%%%%%%%%%%%%%%%%%%%%%%%
\subsection{Overdamped Langevin equation in collective variables}
Theorem \ref{thm:T1} has an immediate application to a practical scenario: overdamped Langevin equation in collective variables  \eqref{sde4}.
The diffusion matrix $M(x)$ in \eqref{sde4}  is symmetric positive definite everywhere in $\Omega$.
Applying Theorem  \ref{thm:T1} to SDE \eqref{sde4} results in the following corollary.
\begin{corollary}
For the  overdamped Langevin equation in collective variables \eqref{sde4} the optimal controller is
\begin{equation}
\label{vstar4}
    v^*(x) = \nabla \log q^+(x),
\end{equation}
and the controlled process is
\begin{align}
    d{X}_t &= [-M({X}_t) \{\nabla F({X}_t) - 2\beta^{-1}\nabla \log q^+\}+\beta^{-1} \nabla \cdot M({X}_t)] dt \notag\\
    &+ \sqrt{2 \beta^{-1}}M({X}_t)^{\frac{1}{2}}dW_t. \label{sde4c}
\end{align}
\end{corollary}
%%%%%%%%%%%%%%%%%%%%%%%%%%%%%%%%%%%%%%%%%%%%%%%%%%%%%%%%%%
%%%%%%%%%%%%%%%%%%%%%%%%%%%%%%%%%%%%%%%%%%%%%%%%%%%%%%%%%%
%%%%%%%%%%%%%%%%%%%%%%%%%%%%%%%%%%%%%%%%%%%%%%%%%%%%%%%%%%
\subsection{Full Langevin equations}
The Langevin dynamics \eqref{sde5} can be written as follows 
\begin{equation}
\label{sde51}
    d\begin{bmatrix}
X_t\\
P_t
\end{bmatrix} = \begin{bmatrix}
m^{-1}P_t\\
-(\nabla U(X_t) + \gamma P_t)
\end{bmatrix}dt + \sqrt{2\gamma \epsilon m}\begin{pmatrix}
0\\
I
\end{pmatrix}dW_t
\end{equation}
where the diffusion matrix is
\begin{equation}
\label{sigma5}
    \sigma(X_t) = \sqrt{2\gamma \epsilon m}\begin{pmatrix} 0\\ I \end{pmatrix} \in \mathbb{R}^{2d \times d}.
\end{equation}
The application of Theorem \ref{thm:T1} to SDE \eqref{sde51} gives the following controlled process.
\begin{corollary}
For the Langevin dynamics \eqref{sde51}, the optimal controller $v^*$ can be chosen to be
\begin{equation}
\label{vstar5}
v^*(x,p)=\nabla_p\log q^+(x,p)
\end{equation}
and the corresponding controlled process is
\begin{equation}
\label{sde5c}
    \begin{cases}
    d{X}_t = m^{-1} {P}_t dt\\
    d{P}_t = \left(-\nabla U({X}_t) - \gamma {P}_t + 2\gamma \epsilon m \nabla_p \log q\right)dt + \sqrt{2\gamma \epsilon m}dW_t.
    \end{cases}
\end{equation}
\end{corollary}

The form of the optimal controller $v^*$ follows from \eqref{vstar} and \eqref{sigma5}. 
Note that \eqref{vstar} and \eqref{sigma5} do not define $v^*$ uniquely as $v^*$ can have arbitrary first $d$ components corresponding to the coordinate subspace. However, the control in SDE \eqref{sde5c} is of the form $\sigma\sigma^\top v^*$ where
$$
\sigma\sigma^\top = 2\gamma \beta^{-1}m \left[\begin{array}{cc}0&0\\0&1\end{array}\right].
$$
Hence the components of $v^*$ in the coordinate subspace are eliminated in SDE \eqref{sde5c}.

%%%%%%%%%%%%%%%%%%%%%%%%%%%%%%%%%

\section{Estimation of the transition rate}
\label{sec:transition_rate}
The problem of estimating the transition rate between metastable states has been one of the central problems addressed by chemists, physicists, and mathematicians working on quantifying rare events and remains a subject of active research \cite{Lelievre2023Hill,Lelievre2023Est}. Practical methods for finding transition rates can be roughly divided into two categories, splitting and reweighting. 

Splitting methods, e.g. transition interface sampling  \cite{TIS_Bolhuis_2003} and forward flux sampling \cite{FFS_rate_Allen2005},  stratify $\Omega_{AB}$ using level sets of a reaction coordinate. The level sets are denoted by $\lambda_i$, $i=0,\ldots,n$, where $\lambda_0 = \partial A$ and $\lambda_n = \partial B$. Next, the transition probabilities $\mathbb{P}(\lambda_{i+1}|\lambda_i)$ to reach $\lambda_{i+1}$ starting from $\lambda_i$ before returning to $A$ are estimated.  Then the escape rate from $A$ to $B$ is calculated according to the formula
\begin{equation}
\label{rate_split}
k_{AB} = \frac{\nu_{A,\lambda_1}}{\rho_A}\prod_{i=1}^{n-1}\mathbb{P}(\lambda_{i+1}|\lambda_i),
\end{equation}
where $\nu_{A,\lambda_1}$ is the average number of transitions from $A$ to $\lambda_1$ per unit time and $\rho_A$ is the probability that the system has last visited $A$ rather than $B$. It was shown in  Ref. \cite{DicksonMakarov_2009} that these methods can suffer from an unfortunate choice of the reaction coordinate. 

Reweighting methods use enhanced sampling techniques and restore the statistics of unbiased transition paths using an appropriate reweighting scheme. These include weighted ensemble \cite{Aristoff2022} and the Girsanov reweighting \cite{girsanov_keller_2018}. The Girsanov reweighting was used in \cite{Zhang_Marzouk_2022} in combination with optimal control with a fixed stopping time.

The settings in which we need to determine the transition rate are different than those in the works mentioned above. 
We plan to compute the committor for the system under consideration or for its reduced model. 
This means that we can compute the transition rate using  \eqref{nu1} provided that the invariant density $\mu$ is known.

However, the rate computed in this manner may be inaccurate due to 
\begin{itemize}
\item 
an inaccurate estimate of the normalization constant for the invariant density (see Section \ref{sec:mueller}) and/or
\item
 a suboptimal choice set of collective variables when model reduction is used. 
\end{itemize}
The issue with the normalization constant for the invariant density can be eliminated if the escape rate from the set $A$ is computed instead -- see equation \eqref{EtauA} below. 
The problem with the choice of collective variables can be hard to overcome if the system is complicated. If the underlying dynamics are time-reversible, i.e. given by \eqref{sde4},
and $\xi(x)$ is the set of collective variables then the transition rate $\tilde{\nu}_{AB}$ estimated in the space of collective variables $\xi$ 
is always exaggerated and relates to the original transition rate $\nu_{AB}$ via (Proposition 6 in Zhang, Hartmann, Schuette (2016)\cite{ZhangHartmannSchutte_2016})
\begin{equation} 
\label{ZHS2016}
\nu_{AB} \le\tilde{\nu}_{AB}  = \nu_{AB} +\frac{1}{\beta}\int_{\Omega_{AB}}\nabla[q(x)-\tilde{q}(\xi(x))]^\top M(x)\nabla[q(x)-\tilde{q}(\xi(x))] \mu(x) dx,
\end{equation}
where $q$ and $\tilde{q}$ are the committors computed for the original and reduced systems respectively. 
Equation \eqref{ZHS2016} shows that if the set of the collective variables $\xi$ were perfect, i.e., if $q(x) = \tilde{q}(\xi(x))$, then  $\nu_{AB} = \tilde{\nu}_{AB}$. In particular, this means that the lowering dimensionality per se does not lead to an error in the transition rate. 
Otherwise, there will be a model reduction error in the transition rate. 

We propose the following scheme { utilizing the fact that the controlled SDE \eqref{eq:sde1c} with the optimal controller \eqref{mincost} exactly matches SDE \eqref{sdeB1c} that governs the transition trajectories of SDE \eqref{sdeB1}. In particular, this means that the expected crossover time $\mathbb{E}[\tau_{AB}]$ for trajectories of SDE \eqref{eq:sde1c} with \eqref{mincost} is the same as that for the transition trajectories of \eqref{sdeB1}.}  Therefore, first one needs to generate a set of transition trajectories using  { the optimally controlled dynamics \eqref{eq:sde1c} with \eqref{mincost}. This allows us to compute the expected crossover time $\mathbb{E}[\tau_{AB}]$ and find the transition rate $\nu_{AB}$:}
\begin{equation}
\label{cross2}
\nu_{AB} = \frac{\rho_{AB}}{\mathbb{E}[\tau_{AB}]}.
\end{equation}
The expected escape time from $A$ can be readily found as well  using \eqref{kAB} and \eqref{cross2}:
\begin{equation}
\label{EtauA}
    \mathbb{E}[\tau_{A}] = \mathbb{E}[\tau_{AB}]\frac{\rho_{A}}{\rho_{AB}}.
\end{equation}
The probabilities $\rho_A$ and $\rho_{AB}$ will be estimated using the computed committors and formulas \eqref{rhoAB} and \eqref{rhoA}. 
The error in their estimates due to model reduction via imperfect collective variables is less impaired than the error in the rate.
The reason is that the formula for the rate uses the gradient of the committor while the formulas for $\rho_A$ and $\rho_{AB}$ involve only the committors. 
The error in the gradient of the committor is amplified due to the differentiation. 
This effect can be illustrated using an example from \cite{LL2010} found in \ref{app:LLexample}. { Moreover, the estimate of the expected crossover time $\mathbb{E}[\tau_{AB}]$ from the controlled dynamics \eqref{eq:sde1c} remains reasonably accurate even if the estimate of the committor is rough. This issue is explored in \ref{app:E}. }

%%%%%%%%%%%%%%%%%%%%%%%%%%%%%%%%%%%%%%%%%%%
%%%%%%%%%%%
%%%%%%%%%%%.  N U M E R I C S
\section{Numerical solution to the committor problem}
\label{sec:numerics}
This section offers descriptions of numerical methods that we have used for finding the committors for the test problems reported in Section \ref{sec:tests}.
As mentioned in Section \ref{sec:goal}, neural network-based (NN-based) solvers have several advantages. 
First, they yield a globally defined smooth solution whose derivative is readily available due to automatic differentiation.
Second, they do not require artificial boundary conditions on the outer boundary of $\Omega_{AB}$.
Finally, they do not require meshing the space which makes them more amenable for promotion to higher dimensions.
The finite element method (FEM) is used for validation of the committors computed using NN-based methods. Its implementation for the time-reversible dynamics \eqref{sde4} and for the Langevin dynamics is detailed in \ref{sec:appFEM}.

\subsection{NN-solver based on the variational form of the committor problem}
\label{sec:NN1}
Two NN-based solvers for the committor problem \eqref{eg:commproblem} for the overdamped 
Langevin dynamics \eqref{sde2} based on the variational formulation \eqref{eg:variational} were proposed by Khoo, Lu, and Ying (2018) \cite{khooluying} and Li, Lin, and Ren (2019) \cite{LiLinRen2019}. 
These solvers both use the loss function motivated by the minimizing property of the committor \eqref{eg:variational} and require only the first derivatives of the committor. 
They use different solution models for the committor and enforce the boundary conditions at $\partial A$ and $\partial B$ in different ways. 
In \cite{khooluying}, the solution model  involves 
Green's function for Laplace's equation and the boundary conditions are implemented via penalty terms in the loss function.
In \cite{LiLinRen2019}, the solution model is designed similarly to the first NN-based PDE solver (Lagaris et al. (1998) \cite{Lagaris1998}) 
so that it automatically satisfies the boundary conditions. In this work, we chose to use the NN-based solver by Li, Lin, and Ren \cite{LiLinRen2019} as it is simpler and its extension to the committor problem for the overdamped Langevin dynamics in collective variables is straightforward. 

Since the overdamped Langevin dynamics \eqref{sde2} and the overdamped Langevin dynamics in collective variables \eqref{sde4} are time-reversible, the forward and backward committor are related via $q^- = 1-q^+$. Therefore, for brevity, we use the notation $q$ for the forward committor whenever the underlying dynamics is time-reversible.

Following \cite{LiLinRen2019}, we use the following solution model to the committor problem \eqref{eg:commproblem}
\begin{equation}
\label{solmodel}
    q(x;\theta) = (1 - \chi_A(x))[(1-\chi_B(x))\mathcal{N}(x;\theta) + \chi_B(x)], \quad x \in \Omega_{AB},
\end{equation}
where $\mathcal{N}(x;\theta)$ is the output of a fully connected neural network (NN) and $\chi_A(x)$ and $\chi_B(x)$ are smooth approximations to the indicator functions of $\partial A$ and $\partial B$. In this work, we use fully connected neural networks with $L$ hidden layers with ${\sf tanh}$ 
activation functions and the outer layer with the sigmoid function ${\sf sigmoid}(x) = (1+e^{-x})^{-1}$. For example, for $L = 1$ and $L = 2$, the neural networks are
\begin{equation}
\label{eq:NN1}
\begin{cases}
 \mathcal{N}(x;\theta) ={\sf sigmoid}\left[W_2{\sf tanh}(W_1 x+b_1)+b_2\right],& L = 1,\\
 \mathcal{N}(x;\theta) = {\sf sigmoid}\left(W_3 {\sf tanh}\left[W_2{\sf tanh}(W_1 x+b_1)+b_2\right] +b_3\right), & L = 2.\\
\end{cases}
\end{equation}
The argument $\theta$ comprises all entries of the matrices $W_j$ and the shift vectors $b_j$.

The loss function is derived from the minimizing property of the committor \eqref{eg:variational}:
\begin{equation}
\label{var1}
q(x) = \arg\min\int_{\Omega_{AB}}\nabla f(x)^\top M(x) \nabla f(x) e^{-\beta F(x)} dx,
\end{equation}
where the minimum is taken among all functions $f\in \mathcal{H}^1(\Omega_{AB})$ (the Sobolev space) such that $f(\partial A) = 0$ and $f(\partial B) = 1$.
The integral in \eqref{var1} is the expectation of $\nabla f(x)^\top M(x)\nabla f(x)$ where $x$ is a random variable  
with invariant density  proportional to $e^{-\beta F(x)}$ supported in $\bar{\Omega}_{AB}$:
\begin{equation}
\label{var2}
\int_{\Omega_{AB}}\nabla f(x)^\top M(x) \nabla f(x) e^{-\beta F(x)} dx = \mathbb{E}_{\substack{x\in\Omega_{AB}\\x\sim e^{-\beta F}}} \left[\nabla f^\top M\nabla f\right].
\end{equation}
If we want $x$ to be distributed according to a density $\rho$, we rewrite this expectation as
\begin{equation}
\label{var3}
\mathbb{E}_{\substack{x\in\Omega_{AB}\\x\sim e^{-\beta F}}} \left[\nabla f^\top M\nabla f\right] =
\mathbb{E}_{\substack{x\in\Omega_{AB}\\x\sim \rho}} \left[\nabla f^\top M\nabla f\frac{e^{-\beta F}}{\rho}\right].
\end{equation}
The last expectation can be approximated as a sample mean. 
Hence, if the training points $x_k\in\Omega_{AB}$, $1\le k\le K$,  are sampled from a density $\rho$, then the loss function is defined by
\begin{equation}
\label{loss1}
{\sf Loss}(\theta) = \frac{1}{K} \sum_{k=1}^K\left[\nabla q(x_k;\theta)^\top M(x_k) \nabla q(x_k;\theta) \frac{e^{-\beta F(x_k)} }{\rho(x_k)} \right].
\end{equation}

We found that it is advantageous to create the set of training points in two stages. 
First, a large point cloud is generated using the metadynamics algorithm \cite{metadynamics_2002} (also see \cite{LiLinRen2019}).
Then this point cloud is rarefied into a {\it delta-net} \cite{crosskey2017,evans2022}, i.e., a spatially quasiuniform set of points obtained as follows. 
Let $\{x_i\}_{i=1}^N$ be the generated point cloud and $\delta$ be the desired distance between the nearest neighbors in the training set. 
Initially, we assign labels {\sf 0} to all points. 
We take $x_1$, compute the distance from $x_1$ to all other points in the point cloud, assign label {\sf keep} to $x_1$ and labels {\sf discard} to all points at distance less than $\delta$. 
Then we find the point with the smallest index that has label {\sf 0},  
compute the distance from it to all points with label {\sf 0},  and change its label to {\sf keep} and labels of all points at distance less than $\delta$ from it to {\sf discard}. 
And so we continue until there are no more points with labels {\sf 0}. All points labeled as {\sf keep} will form the training set.
The resulting set of training points is spatially quasiuniform. Therefore we set $\rho(x_k) = 1$ in the loss function \eqref{loss1}.

The loss function is minimized using the stochastic Adam optimizer \cite{Kingma2014AdamAM}.
Our codes are written in Python and use the PyTorch library for the implementation of neural networks and automatic differentiation \cite{margotyjx}.

%%%%%%%%%%%%%%%%%%%%%%%%%%%%%%%

\subsection{PINNs for Langevin dynamics}
\label{sec:PINN}
The Langevin SDE \eqref{sde5}  is not time-reversible
and its generator \eqref{gen5} is not self-adjoint. As a result, the variational formulation of the committor problem is not available.
In this case, we opt to use the NN-based solver called the \emph{physics-informed neural networks} (PINNs)  proposed by Raissi, Perdikaris, and Karniadakis (2019) \cite{Raissi2017PhysicsID} (also see  \cite{Karniadakis2021}). In this solver, the loss function is the sum of the mean squared discrepancy between the left- and right-hand side of the PDE and the mean squared error at the Dirichlet boundary  of the domain:
\begin{align}
{\sf Loss}(\theta) & = \frac{1}{K}\sum_{(x_k,p_k)\in\Omega_{AB}}\left| \mathcal{L}\mathcal{N}(x_k,p_k;\theta) \right|^2  \label{pinns} \\
&+
 \frac{1}{K_{\partial A}}\sum_{(x_k,p_k)\in{\partial A}} |\mathcal{N}(x_k,p_k)|^2 +
 \frac{1}{K_{\partial B}}\sum_{(x_k,p_k)\in{\partial B}} |\mathcal{N}(x_k,p_k) - 1|^2.\notag 
 \end{align}
 Here, $K$, $K_{\partial A}$, and $K_{\partial B}$ are the numbers of training points in $\Omega_{AB}$, $\partial A$, and $\partial B$ respectively, $\mathcal{L}$ is the generator given by 
 \eqref{gen5}, and $\mathcal{N}(x_k,p_k;\theta)$ is a neural network defined similar to \eqref{eq:NN1}.
 As before, the loss function is minimized using the stochastic Adam optimizer.

%%%%%%%%%%%%%%%%%%%%%%%%%%%%%%%%%%%%%%%%%%%%

\section{Test problems}
\label{sec:tests}
In this section, we demonstrate the effectiveness of the proposed methodology on { the following four} test problems:
\begin{enumerate}
\item the overdamped Langevin dynamics with Mueller's potential, a common 2D test problem in chemical physics,
\item { the overdamped Langevin dynamics with the rugged Mueller potential in 10D with settings as in Ref. \cite{LiLinRen2019},}
\item the bistable Duffing oscillators with added white noise, and
\item the overdamped Langevin dynamics in collective variables for the Lennard-Jones-7 system in 2D.
\end{enumerate}

The codes for these examples are available on GitHub:
\begin{itemize}
    \item neural network-based committor solvers and sampling algorithms for the transition path process \cite{margotyjx};
    \item Finite element committor solvers \cite{mar1akc}.
\end{itemize}

{ The dynamics in test problems 1, 2, and 4 are time-reversible. Therefore, as mentioned in Section \ref{sec:NN1}, the forward and backward committors are related via $q^{-} = 1-q^+$. Hence, it suffices to compute only the forward committor in these problems. For brevity, we will denote the forward committor simply by $q$ in Sections 6.1, 6.2, and 6.4 containing test problems 1, 2, and 4 respectively.}

\subsection{Mueller's potential}
\label{sec:mueller}
We first consider an overdamped Langevin equation with Mueller's potential (see Fig. \ref{fig:Mueller1})
\begin{equation}
\label{mueller}
    V(x_1, x_2) = \sum\limits_{i=1}^4 D_i \exp\left\{a_i(x_1 - X_i)^2 + b_i(x_1 - X_i)(x_2 - Y_i) + c_i(x_2 - Y_i)^2 \right\}
\end{equation}
where
\begin{align*}
    [a_1, a_2, a_3, a_4] &= [-1, -1, -6.5, 0.7] \\
    [b_1, b_2, b_3, b_4] &= [0, 0, 11, 0.6] \\
    [c_1, c_2, c_3, c_4] &= [-10, -10, -6.5, 0.7] \\
    [D_1, D_2, D_3, D_4] &= [-200, -100, -170, 15] \\
    [X_1, X_2, X_3, X_4] &= [1, 0, -0.5, -1] \\
    [Y_1, Y_2, Y_3, Y_4] &= [0, 0.5, 1.5, 1]
\end{align*}

% {  AS: I know you have this later as well, but it might be nice have a graph of Muller's potential side by side with these parameters so that the reader can get an idea of what this potential actually is} \jx{Resolved}

\subsubsection{Computing the committor}\label{sec: mueller committor}
The two deepest minima of Mueller's potential are located near $a = (-0.558, 1.441)$ and $b = (0.623, 0.028)$.  
Following \cite{LiLinRen2019}, the sets $A$ and $B$ are chosen to be the balls centered at $a$ and $b$ respectively with radius $r = 0.1$, and the smoothed indicator function functions of $\partial A$ and $\partial B$ are defined as 
\begin{align*}
    \chi_A(x) & = \frac{1}{2} - \frac{1}{2}\tanh \left[1000(\|x - a\|^2 - (r+0.02)^2) \right],\\
    \chi_B(x) & = \frac{1}{2} - \frac{1}{2}\tanh \left[1000(\|x - b\|^2 - (r+0.02)^2) \right].
\end{align*}
The temperature is set to be  $\beta^{-1} = 10$ as in \cite{LiLinRen2019}. At this temperature, the transitions between $A$ and $B$ are rare.

\begin{figure}[hbt!]
  \begin{subfigure}{.5\linewidth}
  \includegraphics[width=\linewidth]{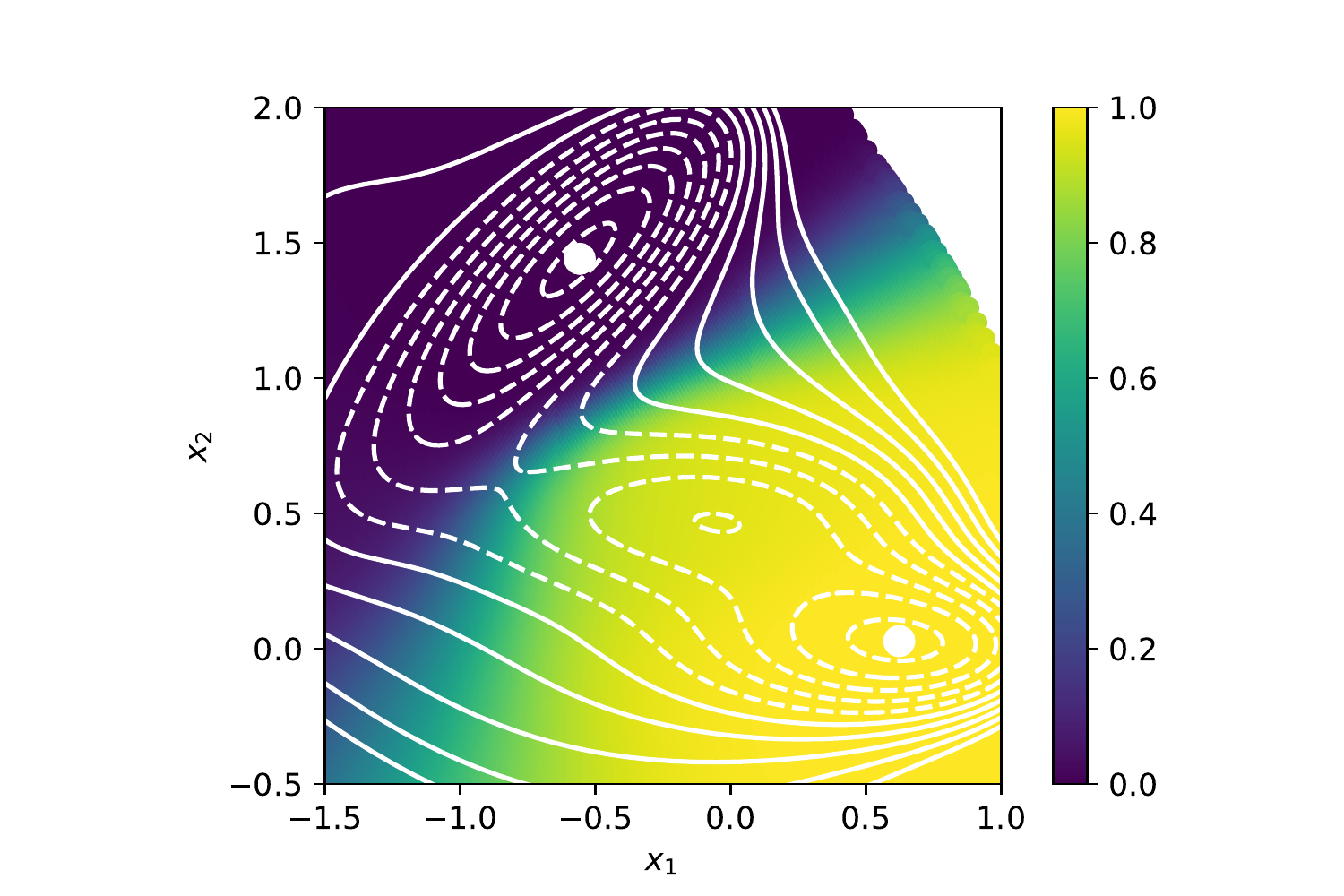}
  \caption{FEM}
\end{subfigure}\hfill % <-- "\hfill"
\begin{subfigure}{.5\linewidth}
  \includegraphics[width=\linewidth]{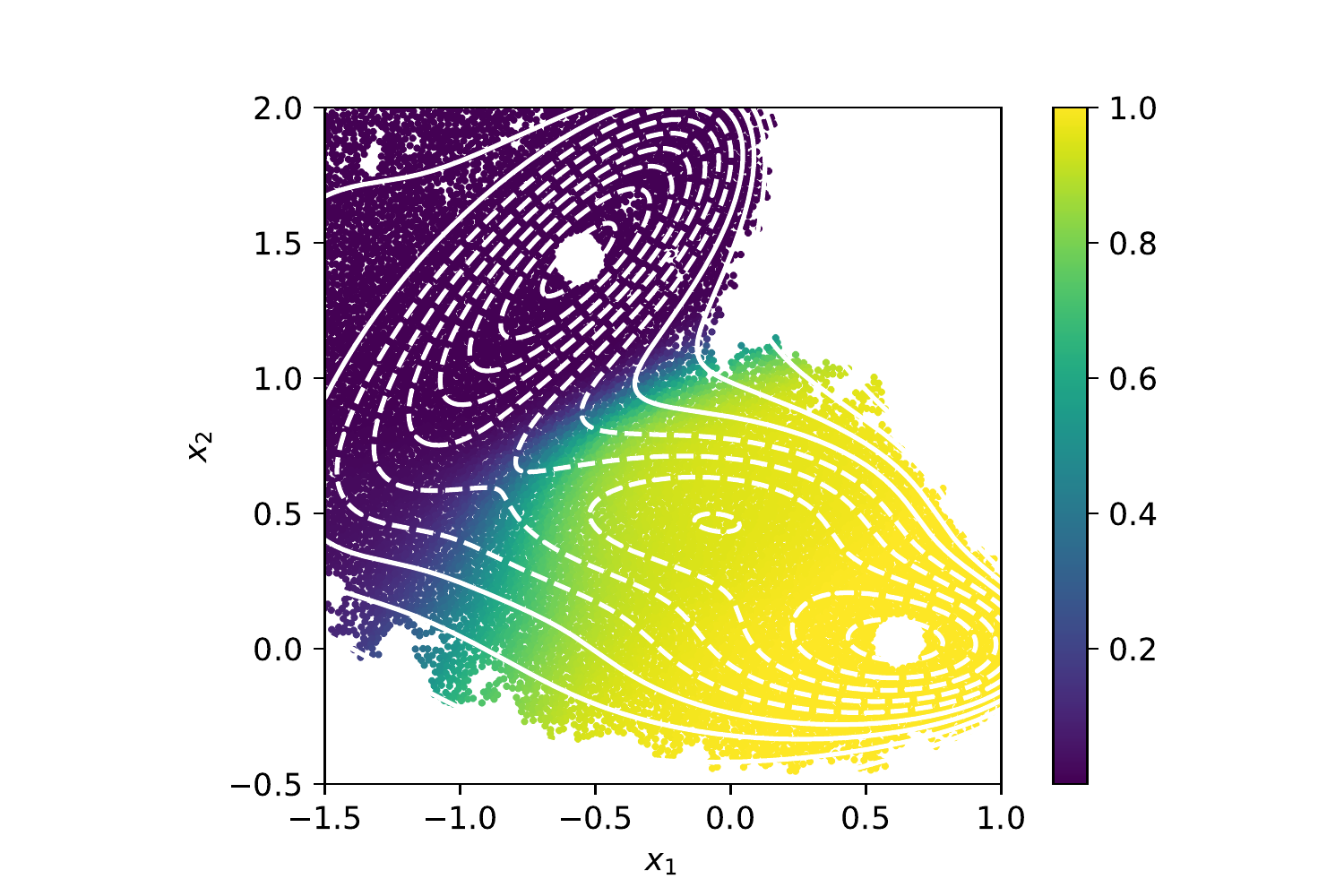}
  \caption{Variational NN}
\end{subfigure}
\caption{The committor for the overdamped Langevin dynamics with Mueller's potential is computed by (a) FEM and (b) NN at $\beta^{-1} = 10$.  The level sets of Mueller's potential are superimposed. The computational domain is larger than the shown box in both cases.}\label{fig:Mueller1}
\end{figure}

We compute the committor using FEM (see \ref{sec:appFEM1}) and the NN-based approach employing the variational formulation of the committor problem (\emph{variational NN} -- see Section \ref{sec:NN1}).
For FEM, the domain $\Omega$ is defined as
\begin{equation}
\label{OmegaMueller}
\Omega = \{x\in\mathbb{R}^2~|~V(x) \le 250\}
\end{equation}
and triangulated using the {\tt DistMesh} algorithm  \cite{Distmesh}.
 We also discretized $\Omega$ using  {\tt mesh2d} \cite{engwirda2009mesh2d} and found that the difference between the committor 
 computed on these two meshes was about $10^{-4}$ in the max norm. The committor computed using FEM is displayed in Fig. \ref{fig:Mueller1}(a).

{ At $\beta = 0.1$, sampling from the invariant Gibbs density leaves the transition region severely underresolved. Therefore, a set of training points  for the variational NN-based solver was generated using a standard enhanced sampling algorithm called \emph{metadynamics} \cite{metadynamics_2002} with settings used in \cite{LiLinRen2019}. Metadynamics is implemented as follows.}
The overdamped Langevin dynamics \eqref{sde1} is simulated with the time step $\Delta t = 10^{-5}$, and Gaussian functions of the form
\begin{equation*}
g_j(x) = w\exp\left(-\frac{(x_1-x_1(t_j))^2}{2\sigma_1^2} - \frac{(x_2-x_2(t_j))^2}{2\sigma_2^2}\right)
\end{equation*}
with height $w = 5$ and $\sigma_1 = \sigma_2 = 0.05$ are added to the potential at times $t_j = 500j\Delta t$, $1\le j\le N_g = 2000$. 
Then, the overdamped Langevin dynamics in the modified potential $V + \sum_j g_j$ is simulated with the same time step $\Delta t = 10^{-5}$, and the initial set of points is recorded.
Finally, the obtained set of points is converted into a spatially quasi-uniform set, a \emph{delta-net} with $\delta = 0.015$, as described in Section \ref{sec:NN1}. 
The resulting training set contains a total of $N_{\sf train} = 15466$ points. 
   
The solution model is given by \eqref{solmodel} with a neural network \eqref{eq:NN1} with $L = 2$ hidden layers   and  $N = 40$ neurons in each layer. The neural network was trained for 1000 epochs at learning rate $\eta = 10^{-4}$. The resulting solution is shown in Fig. \ref{fig:Mueller1}(b). 

To assess numerical errors in the computed committor, 
we use error measures weighted by the probability density of transition trajectories: the weighted mean absolute error (wMAE) and the weighted root mean squared error (wRMSE):
\begin{align}
\label{eq:errormetric}
    \text{wMAE} &= \sum_{i = 1}^{N_{\sf test}} w(x^i)\left|q_{\sf nn}(x^i) - q_{{\sf fem}}(x^i) \right|\\
    \label{eq: weighted RMSE}
    \text{wRMSE} &= \sqrt{\sum_{i = 1}^{N_{\sf test}} w(x^i)\left|q_{\sf nn}(x^i) - q_{\text{fem}}(x^i) \right|^2}
\end{align}
where $q_{\sf fem}$ and $q_{\sf nn}$ are the forward committors computed by FEM and the NN-based solver respectively,  $x^i$, $1\le i\le N_{\sf test}$ are the test points, and the weights $w(x^i)$ are defined so that they are proportional to $\mu_{AB}$ and their sum is one:
\begin{equation}
    w(x^i)= \frac{q_{\sf fem}(x^i)(1-q_{\sf fem}(x^i))\mu(x^i)}{\sum_{j=1}^{N_{\sf test}}q_{\sf fem}(x^j)(1-q_{\sf fem}(x^j))\mu(x^j)}.
\end{equation}
The subset of the nodes of the FEM mesh lying within the box $[-1.5,1]\times[-0.5,2]$ was used as the test point set.

Table \ref{table:Mueller1} shows the wMAE and wRMSE for the variational NN solver with $L=2$ hidden layers and $W=40$ neurons per layer. 
%{ Since the overdamped Langevin dynamics is a time-reversible process, the backward committors can be obtained via $q^- = 1 - q^+$. }
\begin{table}[h]
    \centering
    \begin{tabular}{|c|c|c|c|}
    \hline
    Temperature & NN structure & wMAE & wRMSE  \\
    \hline
    $\beta^{-1} = 10$ & $L = 2$, $W=40$ & 2.6e-3 & 4.1e-3\\
    \hline
    \end{tabular}
    \caption{Errors wMAE and wRMSE in the forward committor for the overdamped Langevin dynamics with Mueller's potential at temperature $\beta^{-1}=10$ 
    computed using the variational NN solver with a spatially quasiuniform set of approximately $16000$ training points. The parameter $L$ indicates the number of hidden layers (see \eqref{eq:NN1}), and $W$ indicates the number of neurons per hidden layer.}
    \label{table:Mueller1}
\end{table}

%%%%%%%%%%%%%%%%%%%%

\subsubsection{Estimation of the transition rate using the controlled process}
\label{sec:Mueller_rate}
The transition rate $\nu_{AB}$ is computed by equation \eqref{cross2}.  
% the ratio of $\rho_{AB}$, the probability that an infinitely long trajectory at a randomly picked moment of time is reactive, and $\mathbb{E}[\tau_{AB}]$, the expected crossover time from $A$ to $B$, -- see equation \eqref{cross2}. {  AS: Can we write the equation instead of saying this as a sentence?}
The expected crossover time $\mathbb{E}[\tau_{AB}]$ is calculated by averaging crossover times of 250 sampled transition trajectories governed by the controlled process
\begin{equation}
\label{MuellerOC}
    dX_t = \left(-\nabla V(X_t) + 2\beta^{-1} \frac{\nabla q_{\sf nn}(X_t)}{q_{\sf nn}(X_t)} \right) dt +\sqrt{ 2\beta^{-1}} dW_t.
\end{equation}
% {  AS: I don't think we ever actually say that $\nabla \log(q) = \frac{\nabla q}{q}$. It's a little obvious, but I think it is still important to say.}

\begin{figure}[hbt!]
\begin{center}
\includegraphics[width=0.5\linewidth]{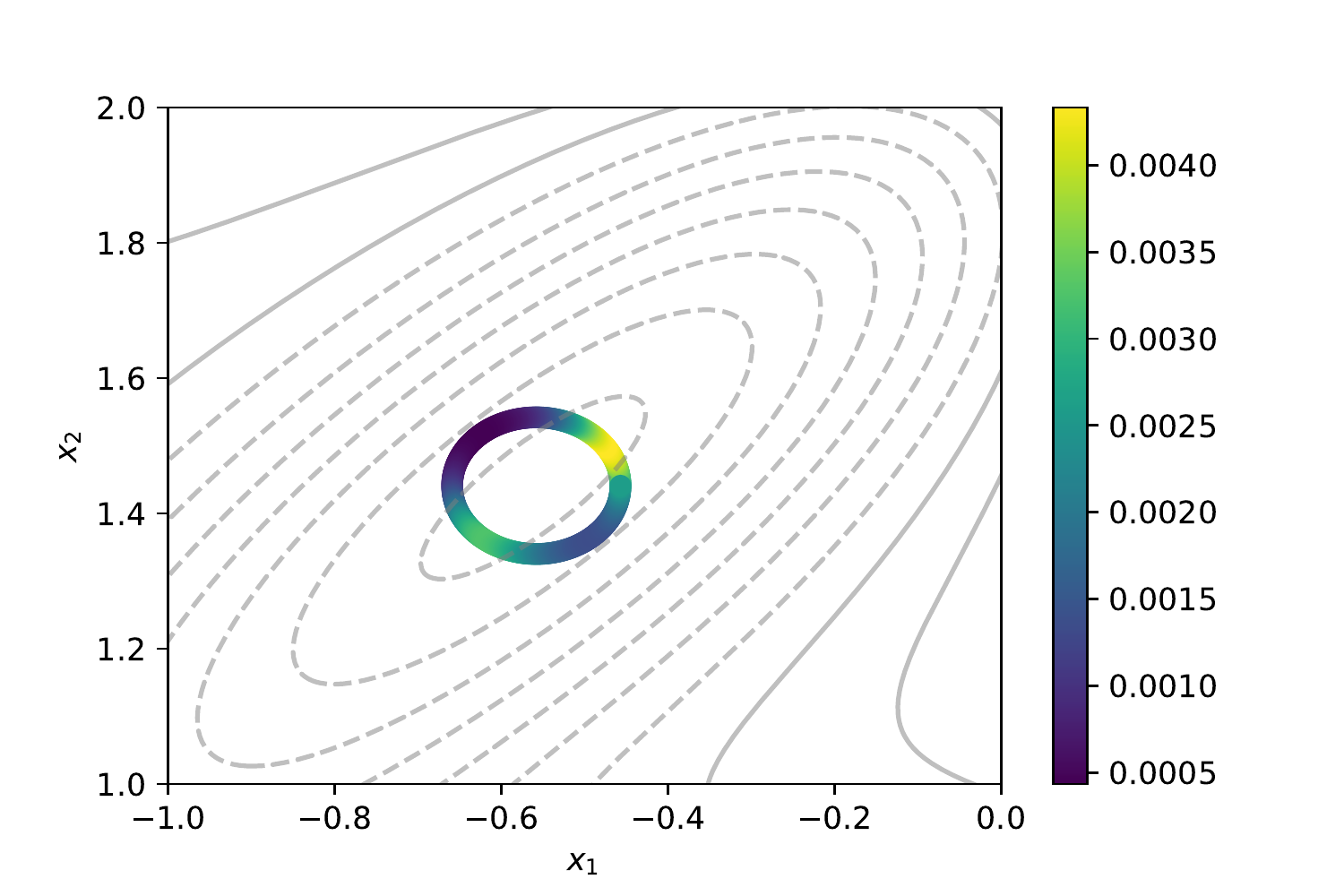}
\caption{The probability distribution of the starting points of transition trajectories on the boundary of $A$ for the overdamped Langevin dynamics in Mueller's potential is shown by color code.  
The gray curves are level sets of Mueller's potential.}
\label{fig:Mueller_Abdry}
\end{center}
\end{figure}

The initial points of these trajectories are sampled according to \eqref{distrA} as follows. First, $N_{\partial A}$ points $x_j$ 
equispaced on a circle of radius $r+\delta r$ where $r$ is the radius of $A$ and $\delta r$ is a small positive number. Then weights are assigned $w_j$ to these points according to
\begin{equation}
\label{w_bdryA}
w_j = \frac{e^{-\beta V(x_j)} |\hat{n}(x_j)\cdot\nabla q_{\sf nn}(x_j)|}{\sum_{k=1}^{N_{\partial A}}e^{-\beta V(x_k)} |\hat{n}(x_k)\cdot\nabla q_{\sf nn}(x_k)|}
\end{equation}
where $\hat{n}(x_j)$ is the outer unit normal to $B_{r+\delta r}(a)$ at $x_j$.
Then the points $x_j$ are sampled according to their probability weights $w_j$ visualized in Fig. \ref{fig:Mueller_Abdry}. We used $r=0.1$, $\delta r = 10^{-3}$, and $N_{\partial A} = 1000$.

\begin{figure}[hbt!]
\begin{subfigure}{.5\linewidth}
\centerline{
  \includegraphics[width=\linewidth]{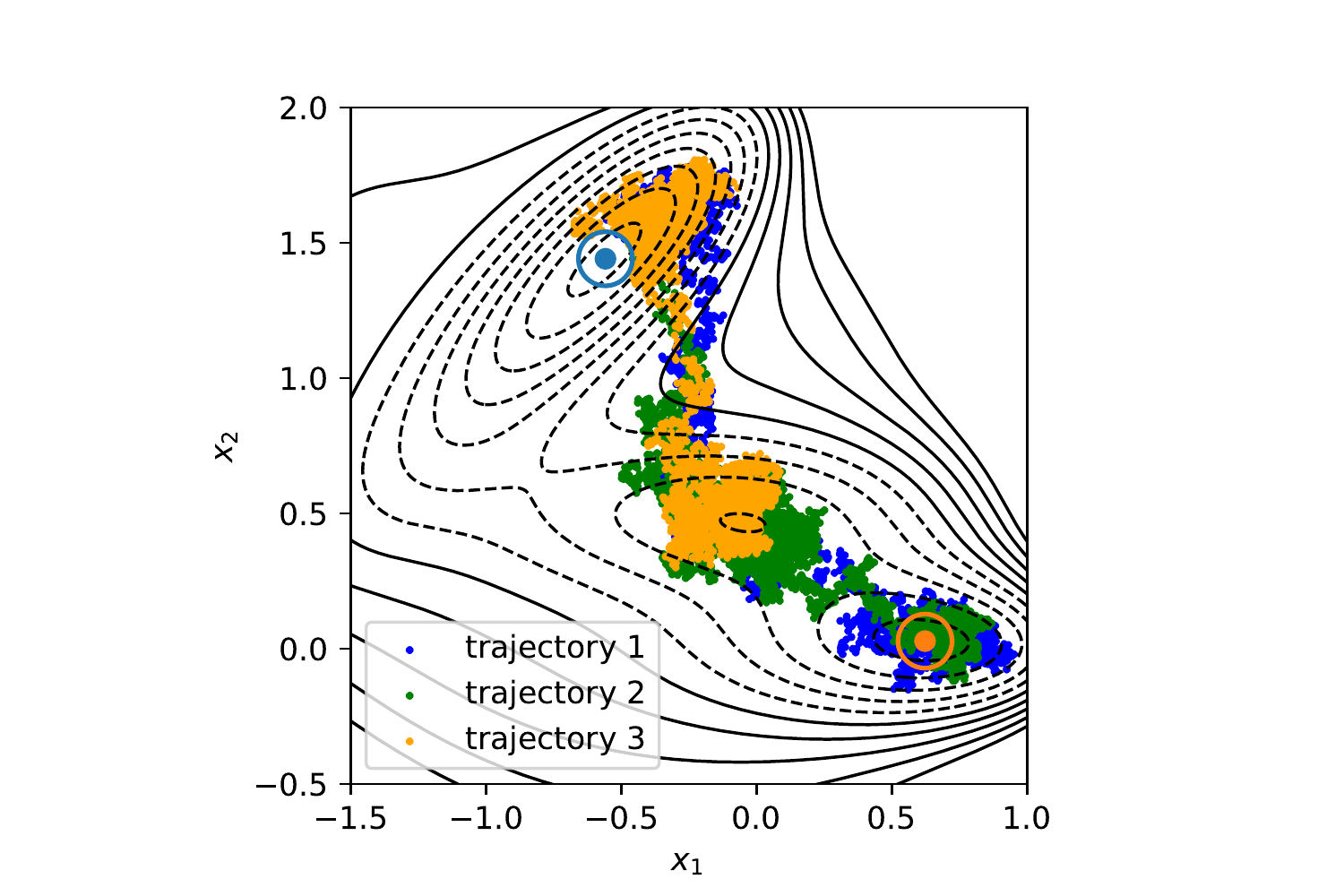}
  }
    \caption{Controlled SDE \eqref{MuellerOC}}
\end{subfigure}
\begin{subfigure}{.5\linewidth}
\centerline{
  \includegraphics[width=\linewidth]{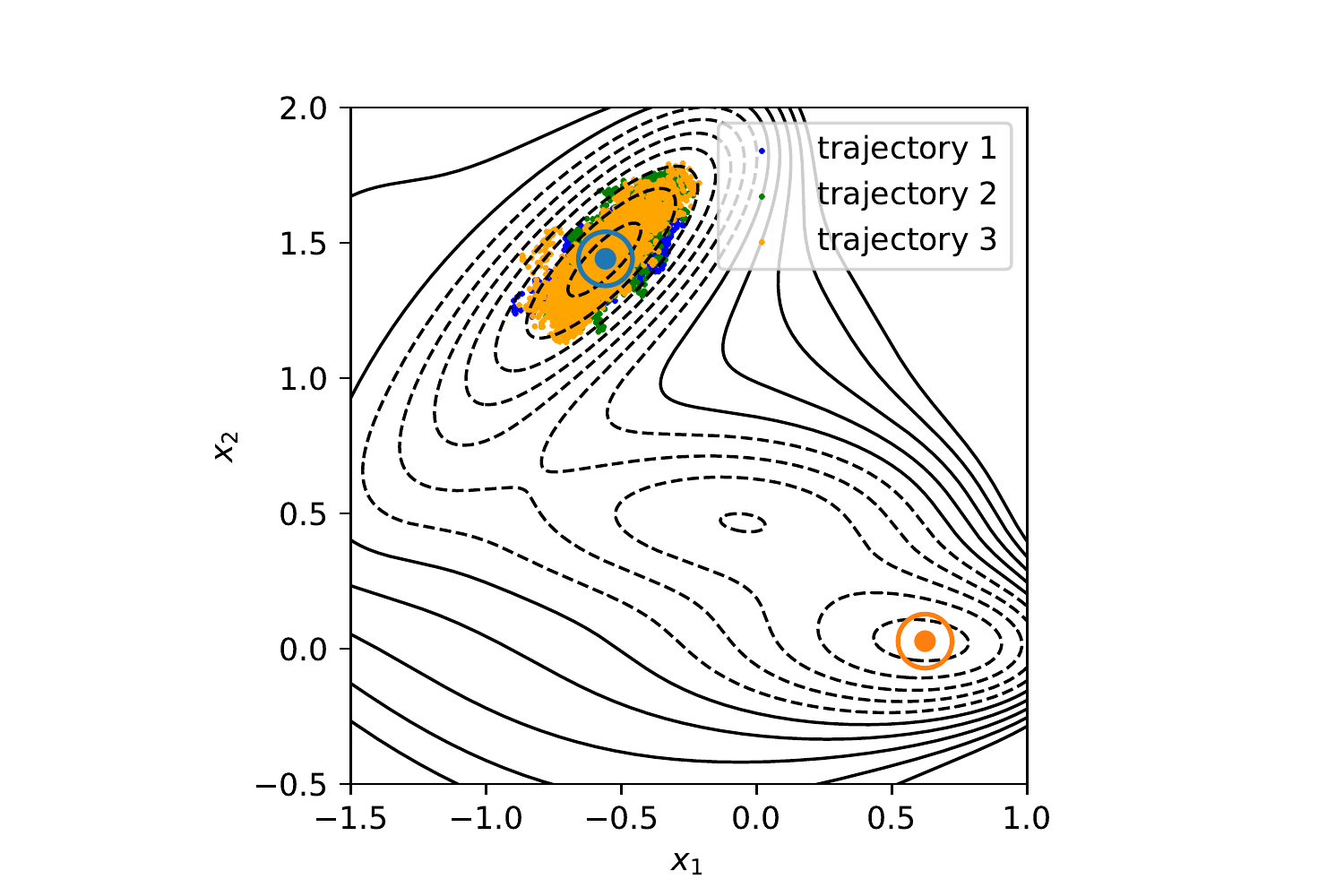}
  }
    \caption{Uncontrolled SDE \eqref{MuellerWOC}}
\end{subfigure}\hfill % <-- "\hfill"
\caption{Comparison of three controlled (a) and uncontrolled (b) processes with the same initial positions and noise realizations at $\beta^{-1} = 10$. Each trajectories consist of 5000 time steps, with time step $\Delta t = 10^{-5}$. The sets $A$ and $B$ are marked by blue and orange circles of radius $r = 0.1$ respectively.}
\label{fig:Mueller2}
\end{figure}

Three samples of trajectories of \eqref{MuellerOC} are displayed in Fig. \ref{fig:Mueller2}(a). 
Three trajectories of 
\begin{equation}
\label{MuellerWOC}
dX_t = -\nabla V(X_t) dt +\sqrt{ 2\beta^{-1}} dW_t
\end{equation}
with the same realizations of the Brownian motion are shown in Fig. \ref{fig:Mueller2}(b) for comparison.

\begin{figure}[htbp]
\begin{center}
  \includegraphics[width=0.5\linewidth]{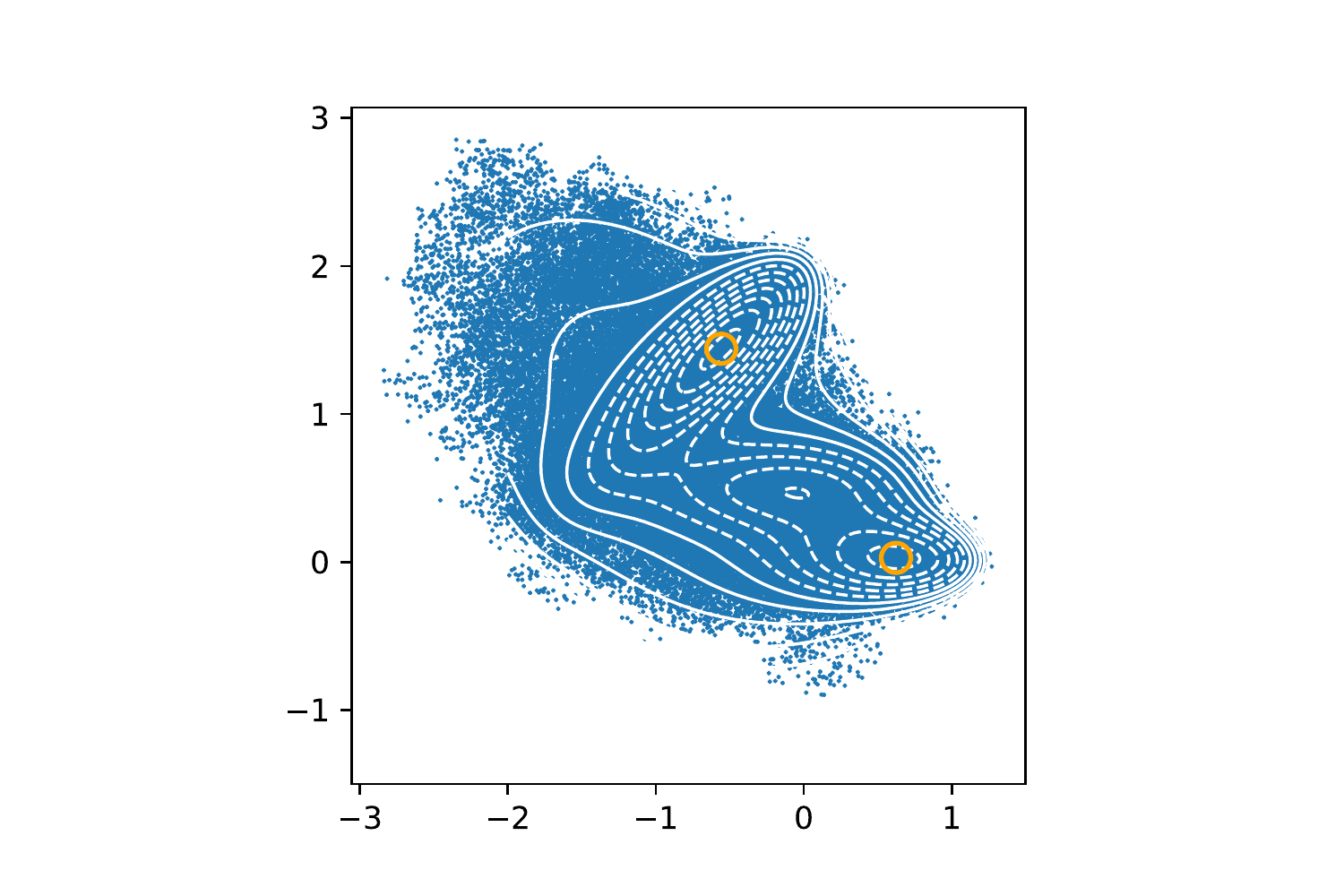}
\caption{The set of points for Mueller's potential used for finding the 
normalization constant for the invariant density and the probability $\rho_{AB}$.}
\label{fig:Mueller_delta_net}
\end{center}
\end{figure}

% The probability $\rho_{AB}$ is estimated using the computed committor according to \eqref{rhoAB}.
Given the FEM committor $q_{\sf fem}$ and the triangulated domain $\Omega _{AB}$, $\rho_{AB}$ is computed directly from \eqref{rhoAB}.
% is computed according to \eqref{rhoAB}.
Given $q_{\sf nn}$ and the set of training points $x_j$ quasiuniformly distributed in $\Omega_{AB}$, $\rho_{AB}$ is obtained by means of Monte Carlo integration: 
% { AS: you should say that you are using Equation (27) (or whatever equation you are using)} \jx{resolved}
\begin{equation}
\label{rhoABfromNN}
\rho_{AB} = % \frac{
\frac{1}{N_{\sf train}}\sum_{j=1}^{N_{\sf train}} Z_V^{-1}e^{-\beta V(x_j)} q_{\sf nn}(x_j)(1-q_{\sf nn}(x_j)).
\end{equation}
The normalization constant $Z_V$ is computed as
\begin{equation}
\label{ZVnn}
Z_V = \frac{1}{N_{\sf train}} \sum_{j=1}^{N_{\sf train}}e^{-\beta V(x_j)}.
\end{equation}
We first generated a  total of $10^5$ points by running metadynamics \cite{metadynamics_2002} and then subsampling it into a spatially quasi-uniform delta-net with $\delta = 0.005$. The resulting set of $56322$ points is shown in Fig. \ref{fig:Mueller_delta_net}.

{ We also calculate the transition rate $\nu_{AB}$ using \eqref{nu1} 
adapted for the overdamped Langevin dynamics \eqref{MuellerWOC}
\begin{equation}
    \label{nu_integral}
    \nu_{AB} = \beta^{-1}Z_V^{-1}\int_{\Omega_{AB}}\|\nabla q(x)\|^2e^{-\beta V(x)} dx
    \end{equation}
with the committors $q_{\sf fem}$ and $q_{\sf nn}$. }
The results are presented in Table \ref{table:Mueller2}. The computation of the 95\% confidence interval is detailed in  \ref{app:confidence_interval}.

The discrepancy between the estimates for $\rho_{AB}$ obtained using the committors $q_{\sf nn}$ and $q_{\sf fem}$ and equations \eqref{rhoABfromNN} and \eqref{rhoAB} is about 3\%.
The difference between the transition rates $\nu_{AB}$ obtained using $q_{\sf nn}$ with Monte Carlo integration and $q_{\sf fem}$ with formula \eqref{nuAB} is less then 2\%.
The expected crossover time $\mathbb{E}[\tau_{AB}]$ computed using the controlled process \eqref{MuellerOC} and $q_{\sf nn}$ is used for computing the transition rate using $\rho_{AB}$ obtained using \eqref{cross1} with $q_{\sf nn}$ and Monte Carlo integration. The result differs from the other rate values by less than 4\%.
Thus, we conclude that
\begin{itemize}
\item the values for the transition rate obtained in three different ways are all consistent; all of them fall into the 95\% confidence interval of the rate value computed using \eqref{cross1};
\item the method for estimating $\rho_{AB}$ and $\nu_{AB}$ using the neural network solver and Monte Carlo integration can be promoted to higher dimensions.
\end{itemize}
Finally, we remark that the estimate for $\nu_{AB}$ is sensitive to the normalization constant $Z_V$ in \eqref{ZVnn}.
It is important to estimate it accurately. For example, Mueller's potential has rather large region in which the potential energy $V$ is relatively low -- see Fig. \ref{fig:Mueller_delta_net}. Sampling points for determining $Z_V$ from a smaller region (a lower sublevel set of $V$) leads to a notable discrepancy in the transition rate estimate. At the same time, the estimate for $\rho_{AB}$ is much less sensitive to the accuracy of $Z_V$. The reason in the difference of sensitivity is that $\nu_{AB}$ has the gradient of the committor in its integral, while $\rho_{AB}$ has the committor itself. 

\begin{table}[h]
    \centering
    \begin{tabular}{|c|c|c|c|}
    \hline
         & Simulations, optimal control & TPT, NN &TPT, FEM\\
    \hline
        $\rho_{AB}$ & NA &2.43e-4& 2.36e-4\\
    \hline
        $\mathbb{E}[\tau_{AB}]$ & 5.05e-2 $\pm$ 0.42e-2 & NA & NA\\
    \hline
        $\nu_{AB}$ & 4.80e-3, [4.43e-3, 5.23e-3] & 4.99e-3 & 4.93e-3\\
    \hline
    \end{tabular}
    % FEM TPT: used to be 2.35e-4 and 4.90e-3
    \caption{Estimates for $\rho_{AB}$, $\mathbb{E}[\tau_{AB}]$ and $\nu_{AB}$ for the overdamped Langevin dynamics in Mueller's potential at $\beta^{-1} = 10$ computed using three schemes. The value for $\nu_{AB} = 4.80e-3$ is obtained with $\rho_{AB}= 2.43e-4$. The 95\% confidence interval is obtained using the error estimate for $\mathbb{E}[\tau_{AB}]$.
    }
    \label{table:Mueller2}
\end{table}
        %%%%%%%%%%%%%%%%%%%%%%%%%%%%%%%%%%%%%%%%%%%%%%%%%%%%%%%%%%%%%%%
        %%%%%%%%%%%%%%%%%%%%%%%%%%%%%%%%%%%%%%%%%%%%%%%%%%%%%%%%%%%%%%%
        %%%%%%%%%%%%%%%%%%%%%%%%%%%%%%%%%%%%%%%%%%%%%%%%%%%%%%%%%%%%%%%

{
\subsection{The rugged Mueller potential in 10D}
The test problem with Mueller's potential can be upgraded by making it 10-dimensional and perturbing its energy landscape with an oscillatory function:
\begin{equation}
\label{eq:Vrugged}
    V(x) = V_0 (x_1,x_2) +
     \gamma \sin(2k\pi x_1) \sin(2k\pi x_2) + \frac{1}{\sigma^2}\sum_{i = 3}^{10} x_i^2.
\end{equation}
Here, $V_0(x_1,x_2)$ is Mueller's potential \eqref{mueller} and $\gamma = 9, k = 5, \sigma = 0.05$ as in \cite{LiLinRen2019}. Following \cite{LiLinRen2019}, the set $A$ and $B$ are chosen to be cylinders centered at $a = (-0.558, 1.441)$ and $b = (0.623, 0.028)$ with radius $r = 0.1$. The exact solution to the committor problem for with $V$ given by \eqref{eq:Vrugged} and such sets $A$ and $B$ is independent of $x_3,\ldots,x_{10}$. This allows us to use the FEM solver in 2D to test the solution computed using the variational NN-based solver in 10D.

\subsubsection{Computing the committor}
We compute the committor using the same procedure as detailed in Section \ref{sec: mueller committor}. For the FEM solver, the computational domain is $\Omega = \{x \in \mathbb{R}^2~|~ V(x) \leq 250\}$. For the variational NN-based solver, a training set of $N_{\sf train} = 64882$ points is generated by sampling $2.1\times 10^5$ points in 10D and rarefying them into a delta-net with $\delta = 0.005$. The neural network  in the solution model \eqref{solmodel} has $L = 3$ hidden layers and $N = 10$ neurons in each layer. The committors computed by the FEM and variational NN-based solvers are shown in Fig. \ref{fig:RuggedMueller1}(a) and (b) respectively.

\begin{figure}[hbt!]
\begin{subfigure}{.45\linewidth}
\centerline{
  \includegraphics[width=\linewidth]{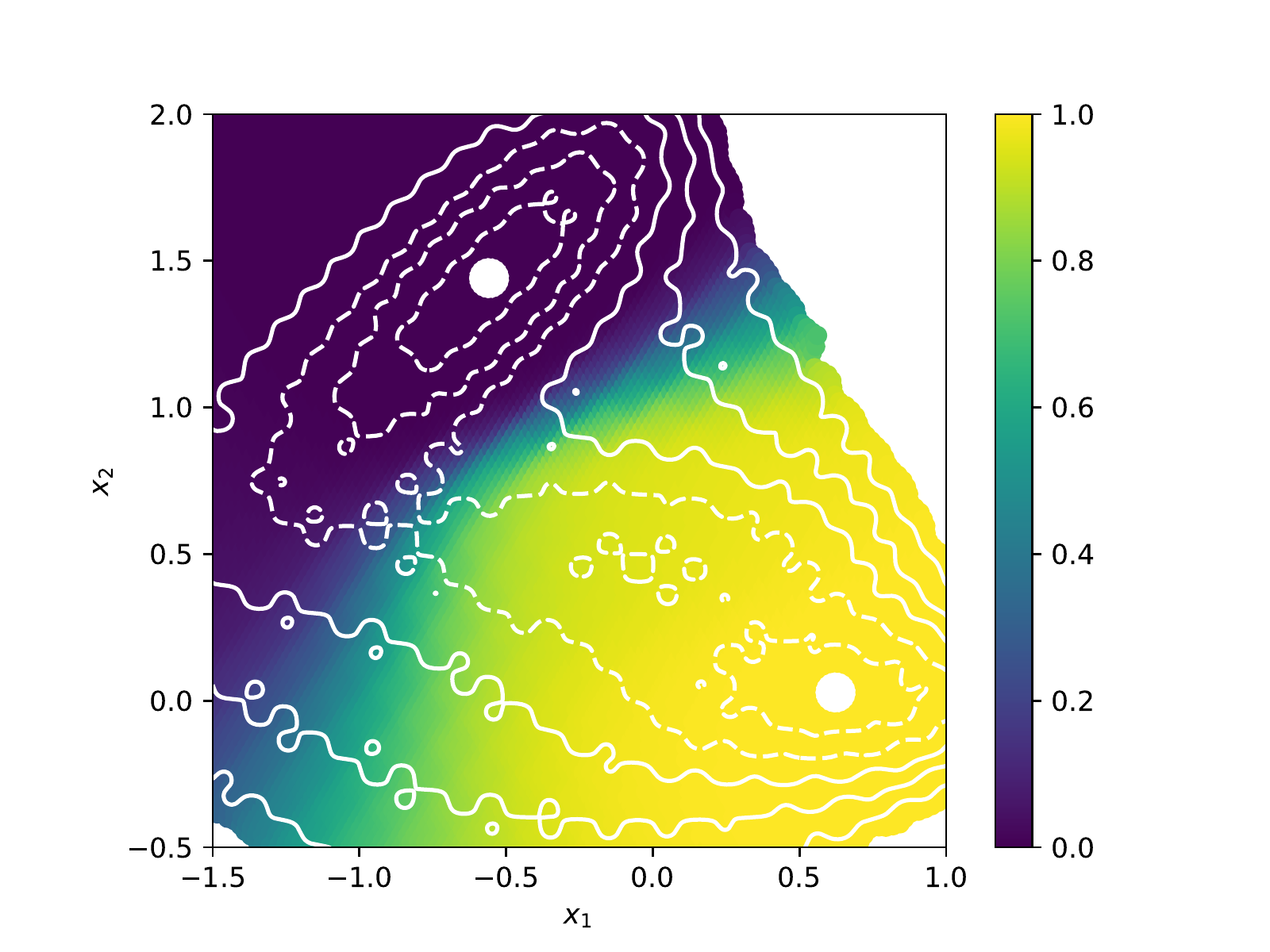}
  }
    \caption{FEM}
\end{subfigure}
\begin{subfigure}{.45\linewidth}
\centerline{
  \includegraphics[width=\linewidth]{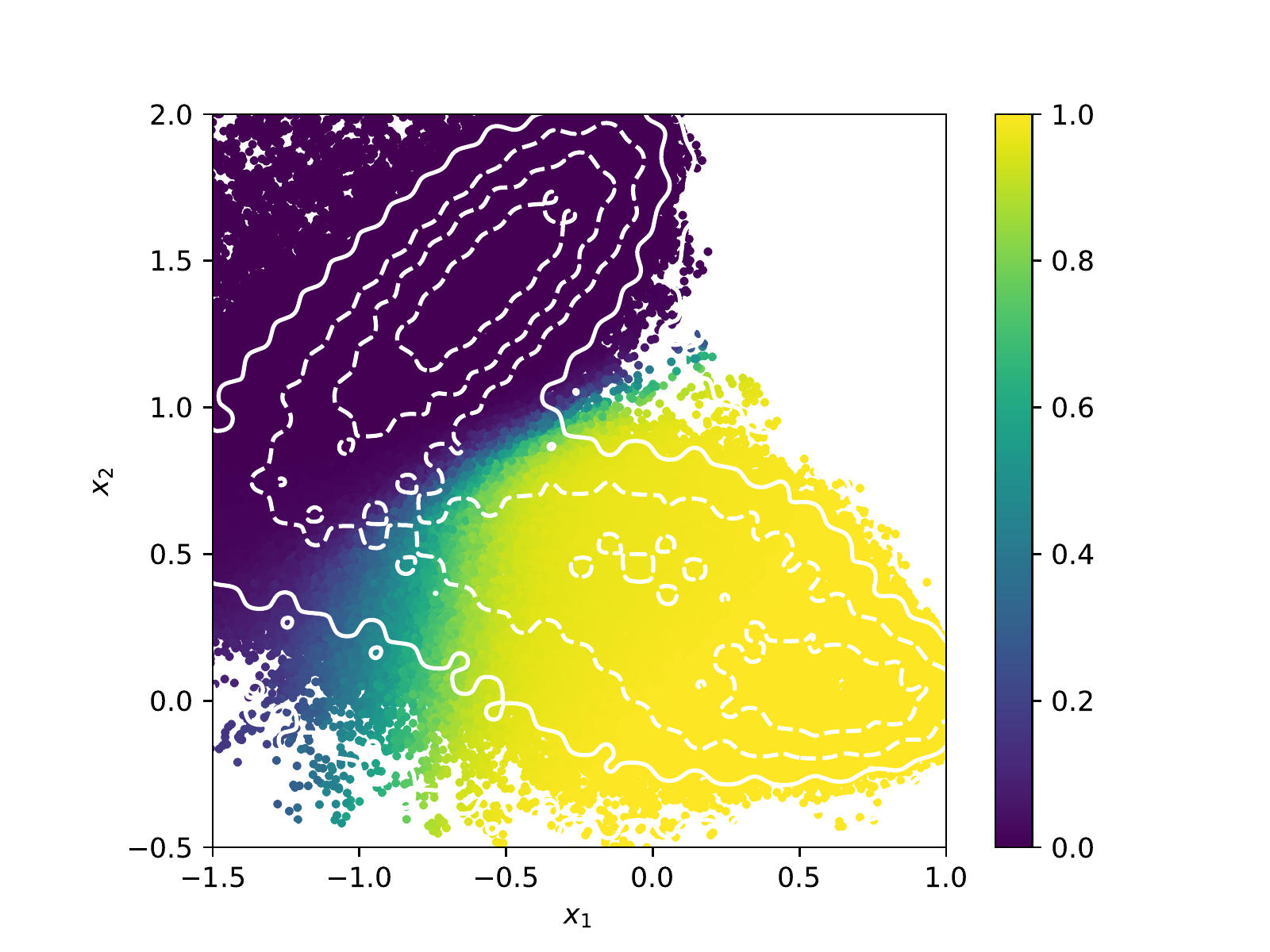}
  }
    \caption{Variational NN}
\end{subfigure}\hfill % <-- "\hfill"
\caption{The committor for the overdamped Langevin dynamics with the rugged Mueller potential \eqref{eq:Vrugged} is computed by (a) FEM  in 2D and (b) by variational NN  in 10D. The temperature is $\beta^{-1} = 10$. The level sets of the potential are superimposed. 
% The computational domain is larger than the shown box in both cases.
}
\label{fig:RuggedMueller1}
\end{figure}

We compute wMAE \eqref{eq:errormetric} and wRMSE \eqref{eq: weighted RMSE} to assess numerical errors. A set of $N_{\sf test} = 44938$ test points in 10D is generated using metadynamics and delta-net postprocessing. The variational NN solution is evaluated at these test points and projected onto the 2D space $(x_1,x_2)$ to compare with the FEM solution. The resulting wMAE and wRMSE are reported in Table \ref{table:ruggedMueller1}.

\begin{table}[h]
    \centering
    \begin{tabular}{|c|c|c|c|}
    \hline
    Temperature & NN structure & wMAE & wRMSE  \\
    \hline
    $\beta^{-1} = 10$ & $L = 3$, $W=10$ & 1.08e-2 & 3.13e-2\\
    \hline
    \end{tabular}
    \caption{Errors wMAE and wRMSE in the forward committor for the overdamped Langevin dynamics with rugged Mueller's potential in 10D at temperature $\beta^{-1}=10$ 
    computed using the variational NN solver with a spatially quasiuniform set of $64882$ training points. The parameter $L$ indicates the number of hidden layers, and $W$ indicates the number of neurons per hidden layer.}
    \label{table:ruggedMueller1}
\end{table}

\subsubsection{Estimation of the transition rate using the controlled process}
The transition rate $\nu_{AB}$ is found using equation \eqref{cross2}.

The expected crossover time $\mathbb{E}[\tau_{AB}]$ is estimated by averaging crossover times of 250 transition trajectories sampled using the controlled process \eqref{MuellerOC} in 10D. The initial points of the trajectories are sampled according to \eqref{distrA} as in the previous test problem. First, $N_{\partial A}$ points are equispaced on the circle of radius $r + \delta r$ centered at $a$ lying in the subspace $(x_1,x_2)$. The weights of these points assigned according to \eqref{w_bdryA} are shown in  Fig. \ref{fig:ruggedMueller_Abdry}. Then these points are sampled based on their probability weight. Three samples of trajectories of the controlled and uncontrolled processes in 10D with the same initial state and the same realizations of the Brownian motion projected onto the $(x_1,x_2)$-subspace are visualized in Fig. \ref{fig:ruggedMueller2} for comparison.

\begin{figure}[hbt!]
\begin{center}
\includegraphics[width=0.5\linewidth]{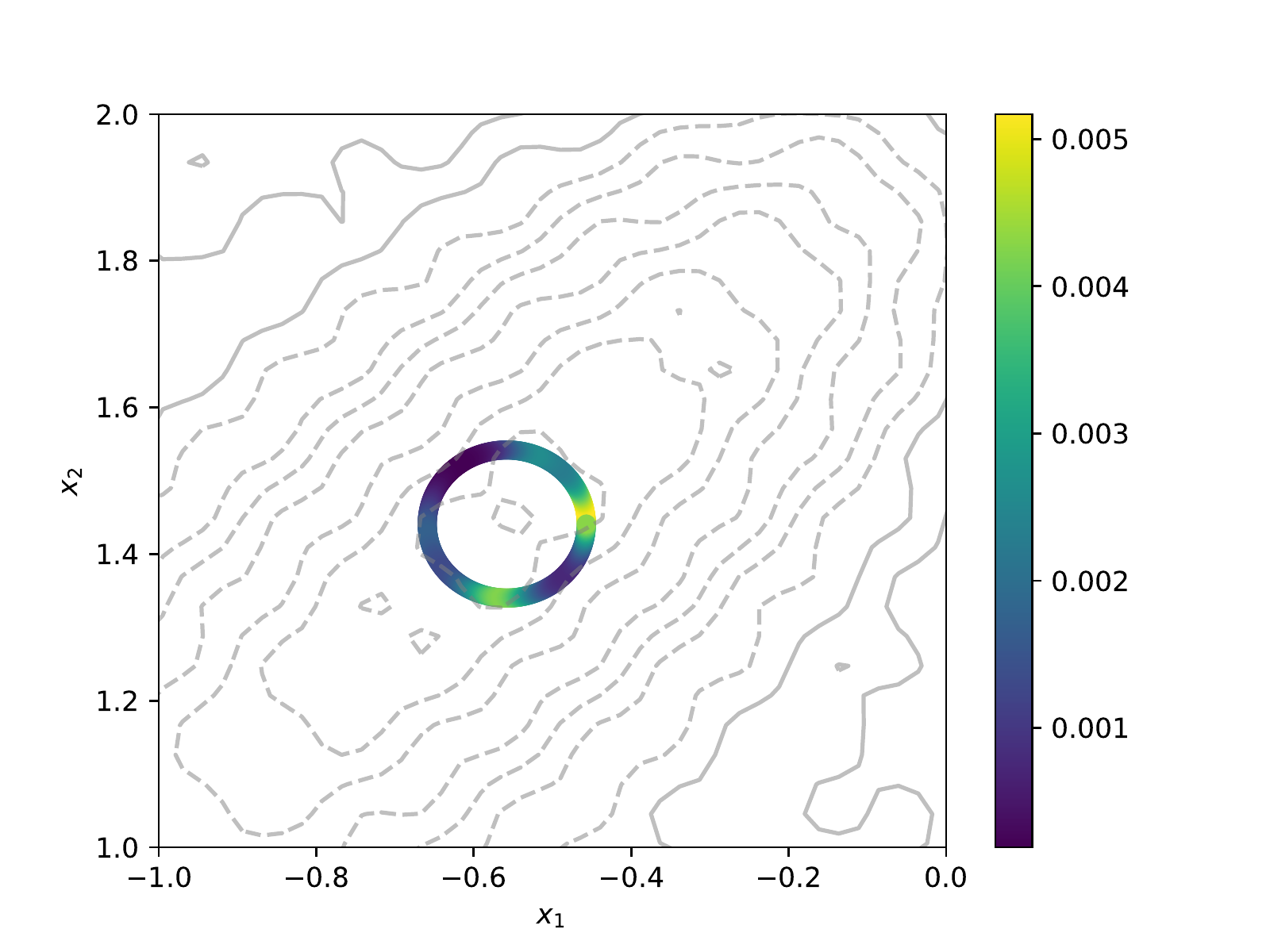}
\caption{The probability distribution of the starting points of transition trajectories on the boundary of $A$ for the overdamped Langevin dynamics in the rugged Mueller potential \eqref{eq:Vrugged} is shown by color code.  The gray curves are level sets of the potential. }
\label{fig:ruggedMueller_Abdry}
\end{center}
\end{figure}

\begin{figure}[hbt!]
\begin{subfigure}{.45\linewidth}
\centerline{
  \includegraphics[width=\linewidth]{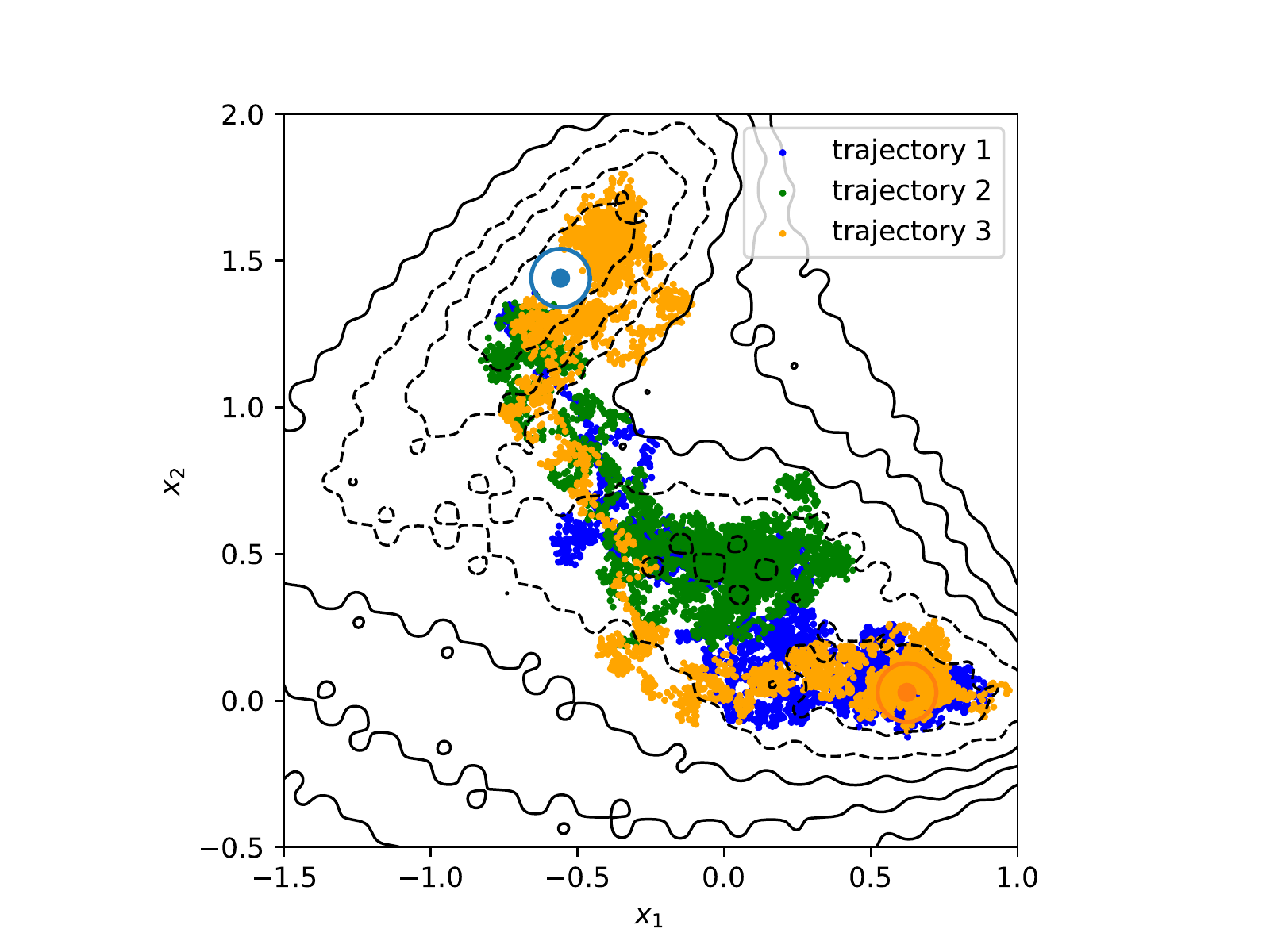}
  }
    \caption{Controlled SDE \eqref{MuellerOC}}
\end{subfigure}
\begin{subfigure}{.45\linewidth}
\centerline{
  \includegraphics[width=\linewidth]{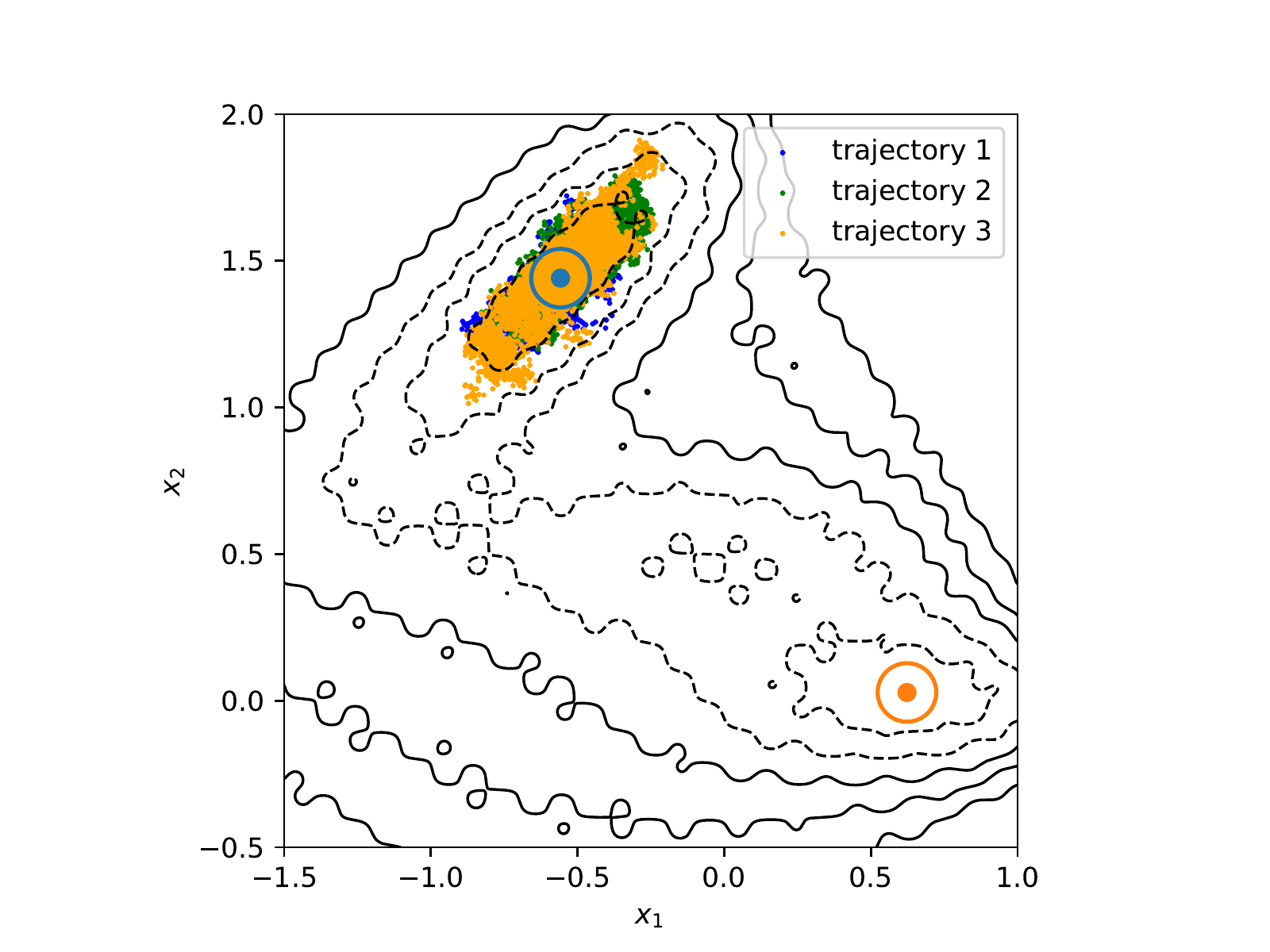}
  }
    \caption{Uncontrolled SDE \eqref{MuellerWOC}}
\end{subfigure}\hfill % <-- "\hfill"
\caption{Comparison of three controlled (a) and uncontrolled (b) trajectories for the rugged Mueller potential in 10D with the same initial positions and the same noise realizations at $\beta^{-1} = 10$ projected onto the $(x_1,x_2)$-subspace. Each trajectory consists of 5000 time steps, with time step $\Delta t = 10^{-5}$. The sets $A$ and $B$ are marked by blue and orange circles of radius $r = 0.1$ respectively.}
\label{fig:ruggedMueller2}
\end{figure}

The probability $\rho_{AB}$ is computed using \eqref{rhoAB} for both committors $q_{\sf fem}$ and $q_{\sf nn}$. Monte Carlo integration is used with $q_{\sf nn}$ over the set of $N_{\sf test} = 44938$ test points as described in Section \ref{sec:Mueller_rate}. 

The transition rate $\nu_{AB}$ is also estimated using \eqref{nu_integral} and the committors $q_{\sf fem}$ and $q_{\sf nn}$ for comparison. 
The results are presented in Table \ref{table:ruggedMueller2}. 

The discrepancy between the estimates for $\rho_{AB}$ obtained using the committors $q_{\sf nn}$ in 10D and $q_{\sf fem}$ in 2D and equations \eqref{rhoABfromNN} and \eqref{rhoAB} is about 4.5\%.
The difference between the transition rates $\nu_{AB}$ obtained using $q_{\sf nn}$ with Monte Carlo integration and $q_{\sf fem}$ with formula \eqref{nuAB} is approximately 13\%. On the other hand, the transition rate computed via \eqref{cross2} uses the expected crossover time and $\rho_{AB}$ acquired from $q_{\sf nn}$, resulting in a transition rate that differs from the FEM rate by 6.5\%. Hence we conclude that our proposed scheme for estimating $\rho_{AB}$ and $\nu_{AB}$ with a neural network solver yields a reasonable accuracy in higher dimensions. 

\begin{table}[h]
    \centering
    \begin{tabular}{|c|c|c|c|}
    \hline
         & Simulations, optimal control & TPT, NN &TPT, FEM\\
    \hline
        $\rho_{AB}$ & NA &2.56e-4& 2.45e-4\\
    \hline
        $\mathbb{E}[\tau_{AB}]$ & 5.94e-2 $\pm$ 0.46e-2 & NA & NA\\
    \hline
        $\nu_{AB}$ & 4.31e-3, [4.0e-3,4.67e-3] & 5.23e-3 & 4.61e-3\\
    \hline
    \end{tabular}
    % FEM TPT: used to be 2.35e-4 and 4.90e-3
    \caption{Estimates for $\rho_{AB}$, $\mathbb{E}[\tau_{AB}]$ and $\nu_{AB}$ for the overdamped Langevin dynamics in the rugged Mueller potential \eqref{eq:Vrugged} at $\beta^{-1} = 10$ in $\mathbb{R}^{10}$ computed using three schemes. The value for $\nu_{AB} = 4.31e-3$ is obtained with $\rho_{AB}=2.56e-4$. The 95\% confidence interval is obtained using the error estimate for $\mathbb{E}[\tau_{AB}]$.
    }
    \label{table:ruggedMueller2}
\end{table}
}
%%%%%%%%%%%%%%%%%%%%%%%%%%%%%%%%%%%%%%%%%%%%%%%%%%%%%%%%%%%%%%%
        %%%%%%%%%%%%%%%%%%%%%%%%%%%%%%%%%%%%%%%%%%%%%%%%%%%%%%%%%%%%%%%
        %%%%%%%%%%%%%%%%%%%%%%%%%%%%%%%%%%%%%%%%%%%%%%%%%%%%%%%%%%%%%%%
        
 \subsection{Duffing Oscillator in 1D}
 \label{sec:Duffing}
 Now we test the proposed methodology on the bistable Duffing oscillator with mass $m=1$, friction coefficient $\gamma = 0.5$, and the potential energy function $V(x) = 0.25(x^2-1)^2$. The dynamics are governed by the Langevin SDE
\begin{equation}
\label{Duffing}
    d\begin{bmatrix} X_t \\ P_t \end{bmatrix} = \begin{bmatrix} P_t \\ -X_t(X_t^2 - 1) - \frac{1}{2} P_t \end{bmatrix}dt + \sqrt{\epsilon}\begin{bmatrix} 0 \\ dW_t \end{bmatrix}.
\end{equation}
The system has two stable equilibria $a = (-1,0)$ and $b=(1,0)$ and an unstable equilibrium at the origin. { The full energy of the system is $H(x,p) = \tfrac{1}{2}p^2 + V(x)$. The invariant probability density is $\mu(x,p) = Z_H^{-1}\exp(-H(x,p)/\epsilon)$.} The sets $A$ and $B$ are chosen to be ellipses with radii $r_x = 0.3$ and $r_y = 0.4$ centered at $a $ and $b$ respectively. Two values of the noise coefficient are used: $\epsilon = 0.1$ and $\epsilon = 0.05$.

%%%%%%%

\subsubsection{Computation of forward and backward committor functions}
\begin{figure}[hbt!]
\begin{subfigure}{.5\linewidth}
\includegraphics[width=\linewidth]{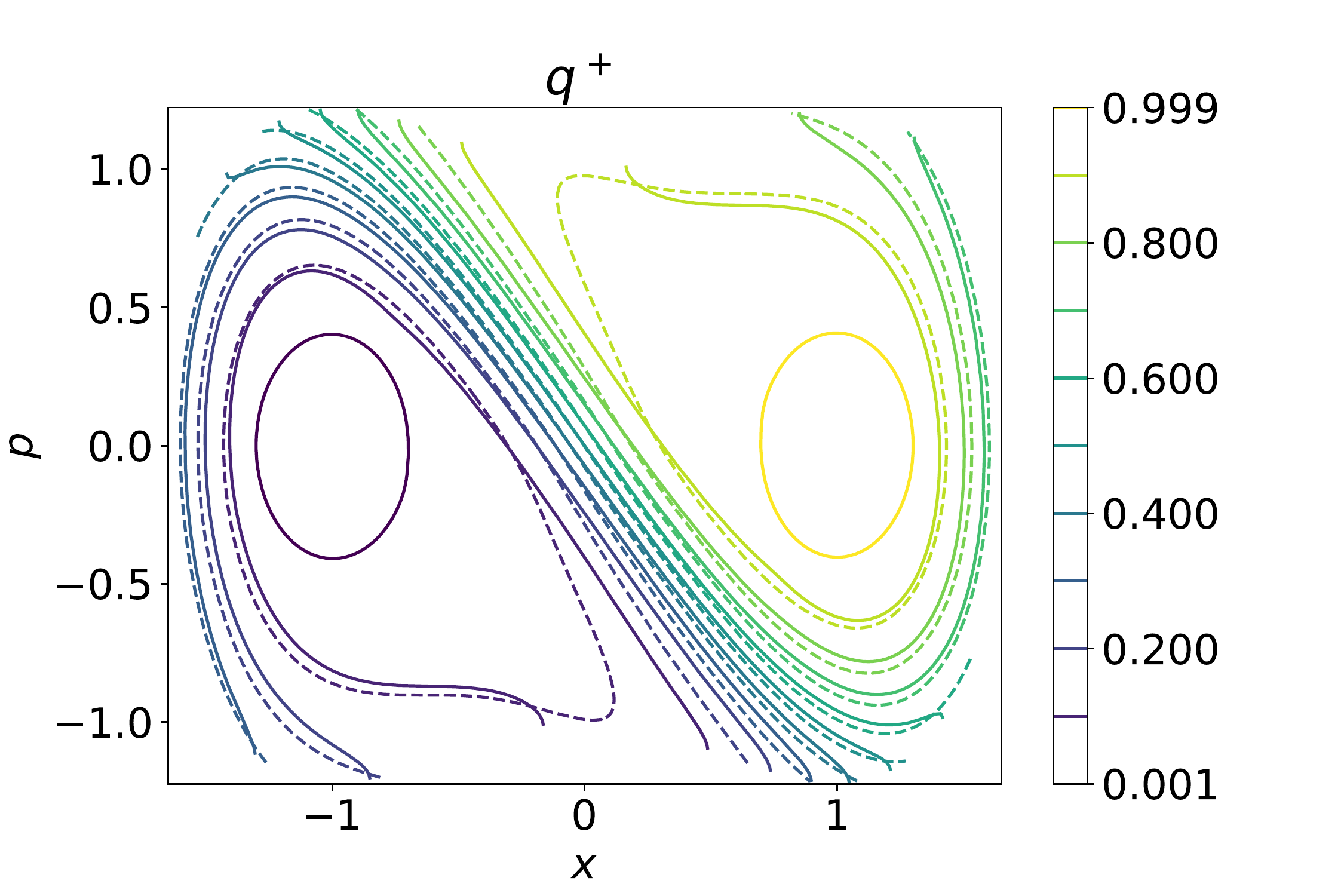}
  \caption{FEM (solid) \& PINN (dashed), $\epsilon = 0.1$}
\end{subfigure}\hfill % <-- "\hfill"
\begin{subfigure}{.5\linewidth}
\includegraphics[width=\linewidth]{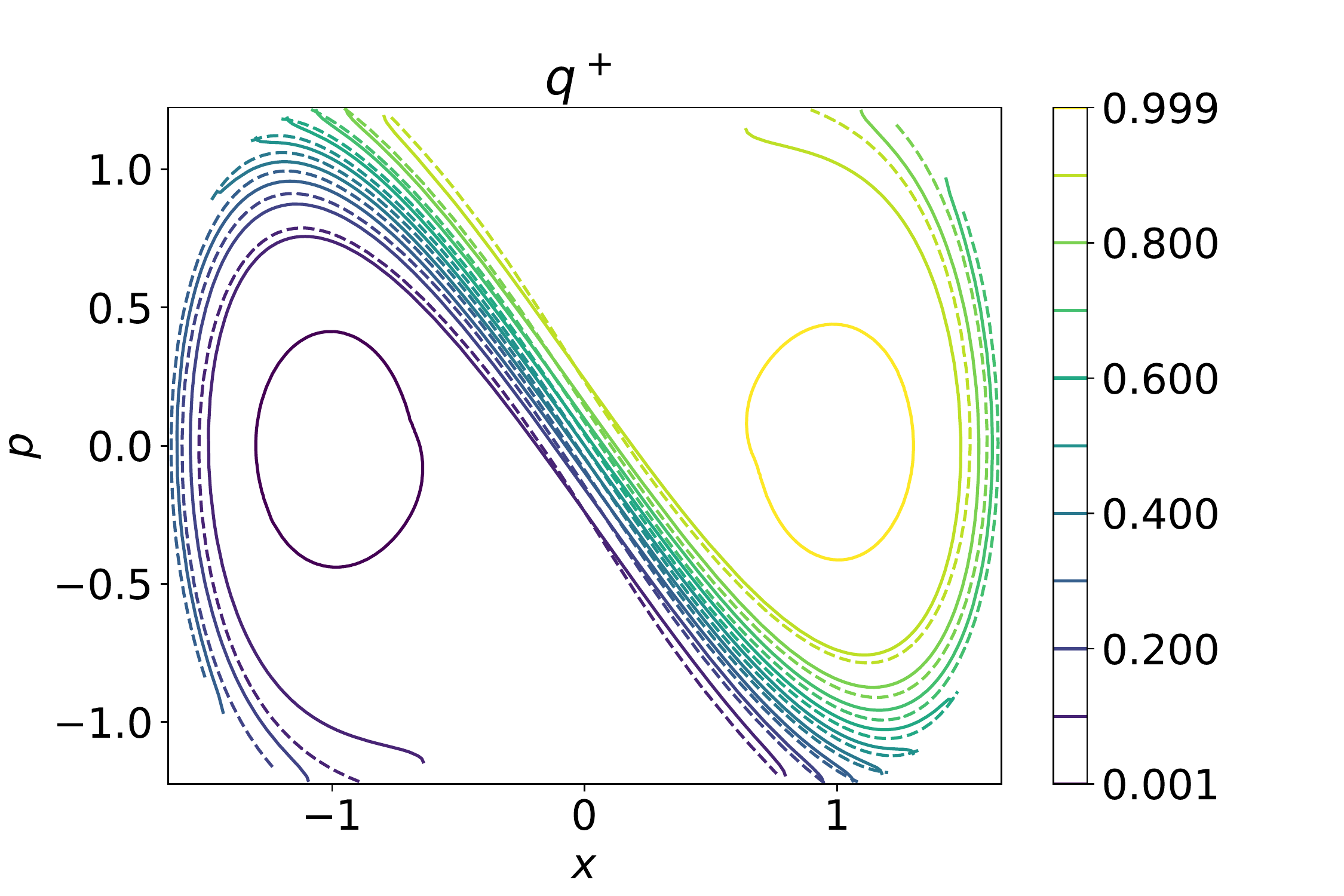}
  \caption{FEM (solid) \& PINN (dashed), $\epsilon = 0.05$}
\end{subfigure}
\caption{Comparison of contours of forward committor computed using 
FEM and PINNs, marked in solid and dashed lines respectively for Duffing Oscillator at 
$\epsilon = 0.1$ and $\epsilon = 0.05$.}
\label{fig:Duffing_q+}
\end{figure}

\begin{figure}[hbt!]
\begin{subfigure}{.5\linewidth}
\includegraphics[width=\linewidth]{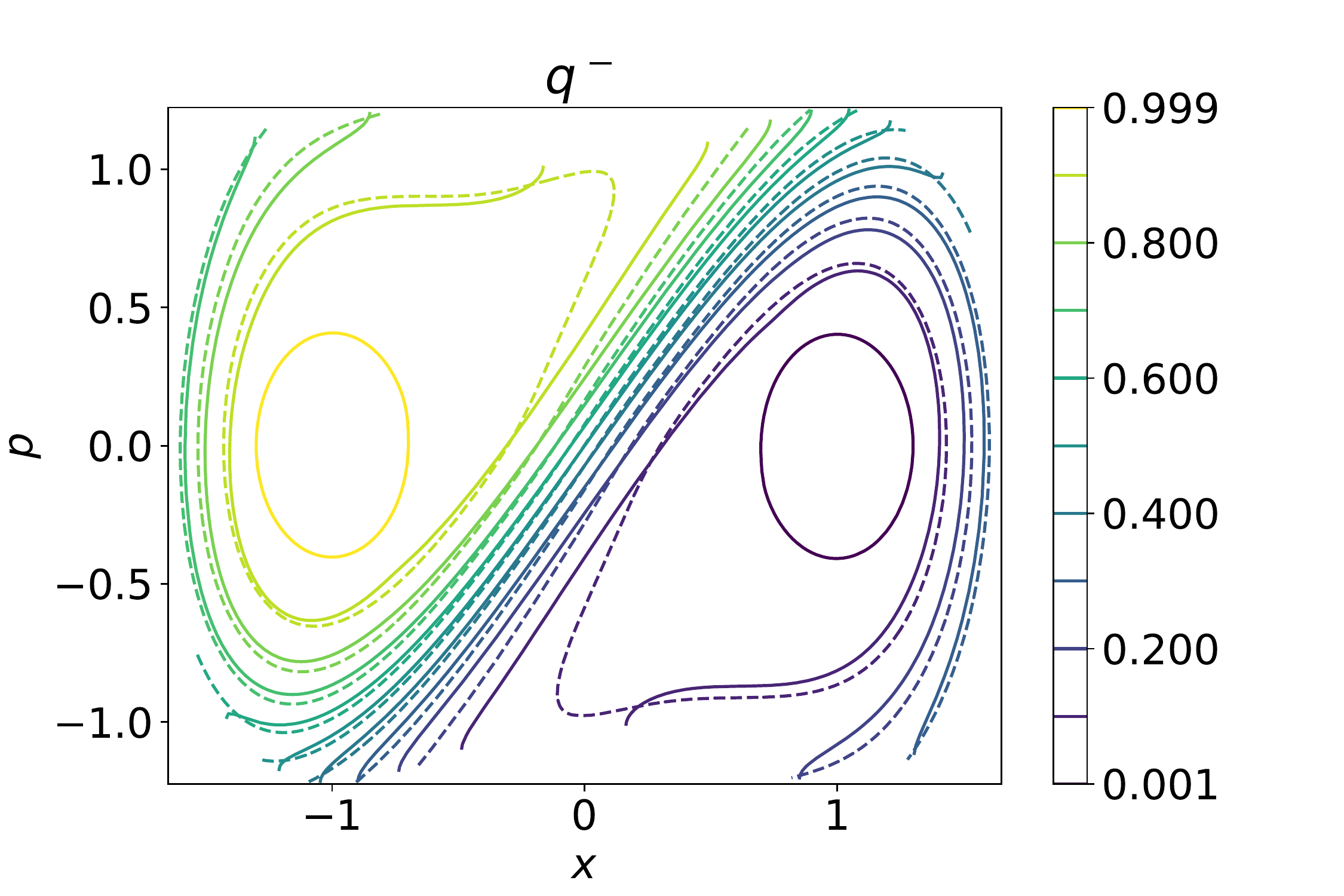}
  \caption{FEM (solid) \& PINN (dashed), $\epsilon = 0.1$}
\end{subfigure}\hfill % <-- "\hfill"
\begin{subfigure}{.5\linewidth}
\includegraphics[width=\linewidth]{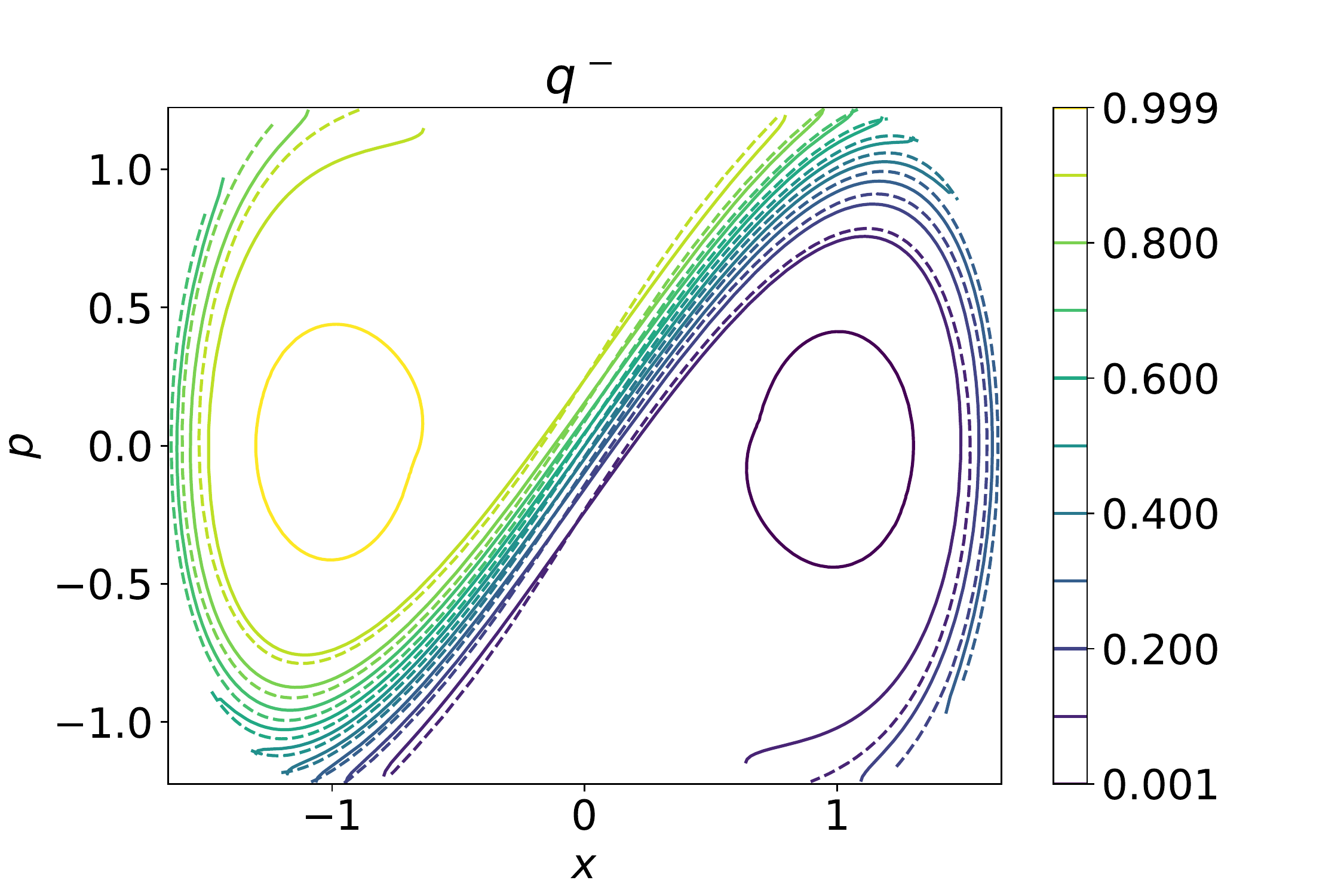}
  \caption{FEM (solid) \& PINN (dashed), $\epsilon = 0.05$}
\end{subfigure}
\caption{Comparison of contours of backward committor computed using 
FEM and PINNs, marked in solid and dashed lines respectively for Duffing Oscillator at 
$\epsilon = 0.1$ and $\epsilon = 0.05$.}
\label{fig:Duffing_q-}
\end{figure}

Since the dynamics \eqref{Duffing} are not time-reversible, we implement the Physics-Informed Neural Network (PINN) approach detailed in Section \ref{sec:PINN} to solve the committor problem.
% { AS: I think you should explain that the fact that it is not time-reversible means that it does not admit a variational formula and thus we have to use a PINN}. 
A uniform grid with a total of 16000 points in the rectangle $[-2.5, 2.5] \times [-2,2]$ is taken as training data. The architecture of the neural network is as in equation \eqref{eq:NN1} with a single hidden layer, $L = 1$, and $W=40$ neurons in it. The Adam optimizer is used with the learning rate $10^{-3}$ for 500 epochs. We also compute the committors using FEM as described in \ref{sec:appFEM2}. 
 The computed forward and backward committors, $q^{+}(x,p)$ and $q^{-}(x,p)$, for $\epsilon = 0.1$ and $\epsilon = 0.05$ are displayed in Figs. \ref{fig:Duffing_q+} and \ref{fig:Duffing_q-} respectively. The theoretical relationship between them is $q^{-}(x,p) = 1-q^{+}(x,-p)$. However, we still computed $q^{-}$ using FEM because the FEM mesh is not symmetric.

% Since the committor problem is hypoelliptic we do not assume that the FEM solution \jx{to be ground truth.}
% % is more accurate than the PINN one. 
% { AS: what does being hypoelliptic have to do with this assumption?} 
We call the discrepancies between the FEM and PINN solutions computed by  \eqref{eq:errormetric}
the weighted mean absolute difference (wMAD) and the weighted root mean square difference (wRMSD) with weights  at the training points $(x_i,p_i)$ given by
\begin{equation}
    w(x_i,p_i) = \frac{{q}^+(x_i,p_i){q}^-(x_i,p_i)e^{-H(x_i,p_i)/\epsilon}}
    { \sum_{j = 1}^{N_{\sf train}}{q}^+(x_j,p_j){q}^-(x_j,p_j)e^{- H(x_j,p_j)/\epsilon}},\quad i=1,\ldots,N_{\sf train}.
\end{equation}
Table \ref{table:Duffing1} shows the wMAD and wRMSD for forward and backward committors computed using PINN and FEM for $\epsilon = 0.1$ and $\epsilon = 0.05$.

\begin{table}[h]
    \centering
%    \begin{tabular}{|m{2cm}|m{2cm}|m{2cm}|m{2cm}|m{2cm}|m{2cm}|}
    \begin{tabular}{|l|l|l|l|l|l|}
    \hline
    $\epsilon$ & NN structure & $q^+$, wMAD & $q^+$, wRMSD & $q^-$, wMAD & $q^-$, wRMSD  \\
    \hline
    $\epsilon = 0.1$ & $L=1$, $W=40$ & 1.6e-2 & 2.0e-2 & 1.8e-2 & 2.2e-2 \\
    \hline
     $\epsilon = 0.05$ & $L=1$, $W=40$ & 1.3e-2 & 2.0e-2 & 1.3e-2 & 2.0e-2\\
    \hline
    \end{tabular}
    \caption{Comparison of the numerical solutions for the forward and backward committors computed using PINNs and FEM. Metrics wMAD and wRMSD are used.
    The parameter $L$ indicates the number of hidden layers (see \eqref{eq:NN1}), and $W$ indicates the number of neurons per hidden layer.
    }
    \label{table:Duffing1}
\end{table}
% { AS: It would be nice to have some interpretation of these results. Does this mean that the PINN does well?} \jx{Resolved}

%%%%%%%
% {
% \subsubsection{Computation of $\rho_{AB}$.}
%  The probability $\rho_{AB}$ is needed for estimating the transition rate. To compute $\rho_{AB}$ for the Duffing oscillators, we explored two distinct committor estimations using PINN and FEM. Given the PINN-based committor, we first rewrote 
% \begin{equation*}
%     \rho_{AB} = \int_{\Omega_{AB}} \mu q^+ q^- dx = \mathbb{E}_{ x \sim \mu,x \in \Omega_{AB}}[q^+(x) q^-(x)]
% \end{equation*}
% We obtained dataset comprising 16,000 uniformly sampled points within rectangular region $[-2.5, 2.5] \times [-2, 2]$ . These points also served as the training dataset for the PINN. Since $\mu(x) = \frac{1}{Z} e^{-\beta U(x)}$ where $Z = \int_{\Omega} e^{-\beta U(x)}dx$, we compute an empirical estimate for $Z$ using uniform points: $\hat{Z} = \sum_{i} e^{-\beta U(x)}dx$. Consequently, for each data point, we derived an approximate $\hat{\mu}(x)$ as $\hat{\mu}(x) = \frac{1}{\hat{Z}}e^{-\beta U(x)}$.

% }

\subsubsection{Estimation of the transition rate using the controlled process}
\label{sec:Duffing_rate}
\begin{figure}[hbt!]
\begin{subfigure}{.25\linewidth}
\includegraphics[width=1\linewidth]{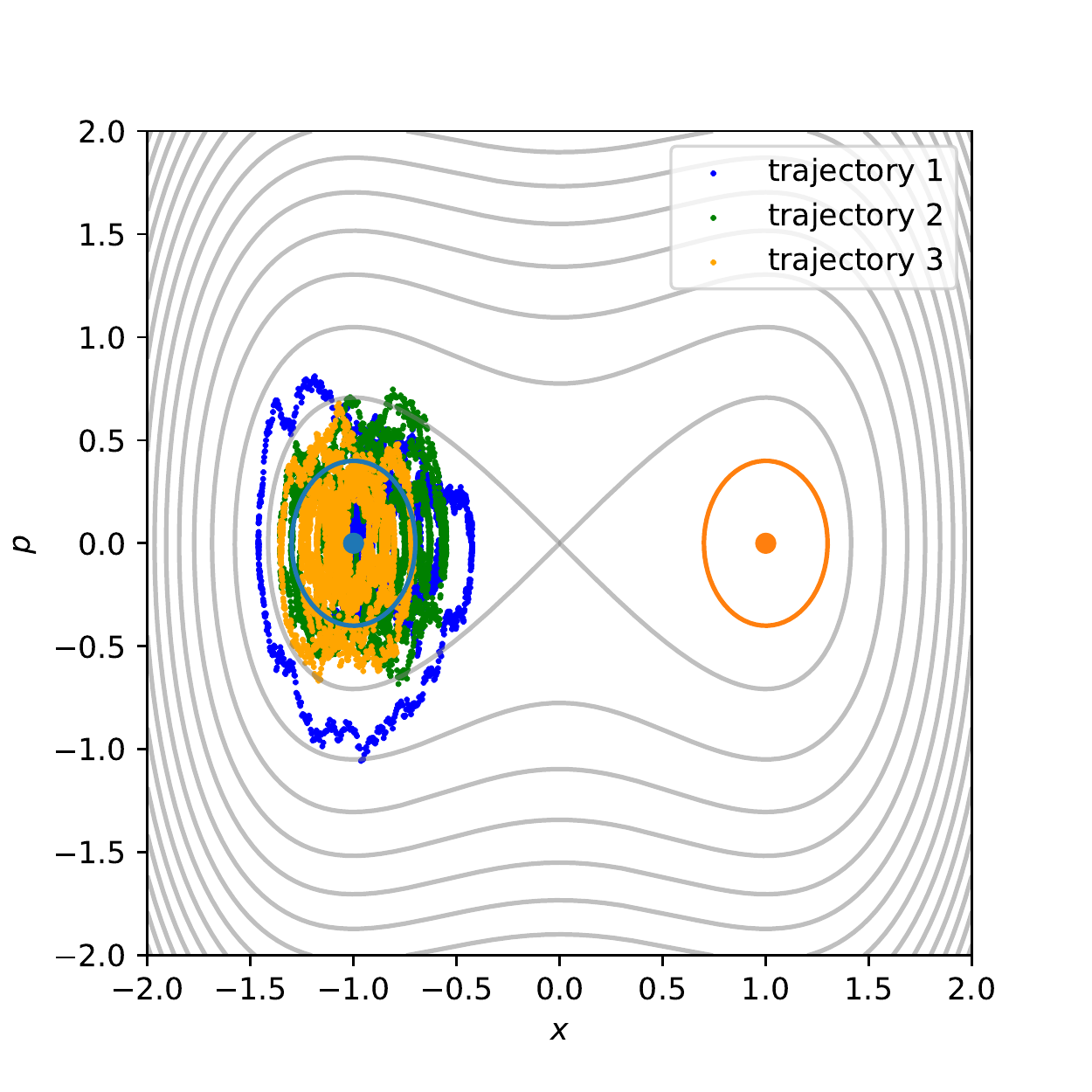}
  \caption{$\epsilon = 0.1$, w/o o/c}
\end{subfigure}\hfill % <-- "\hfill"
\begin{subfigure}{.25\linewidth}
\includegraphics[width=1\linewidth]{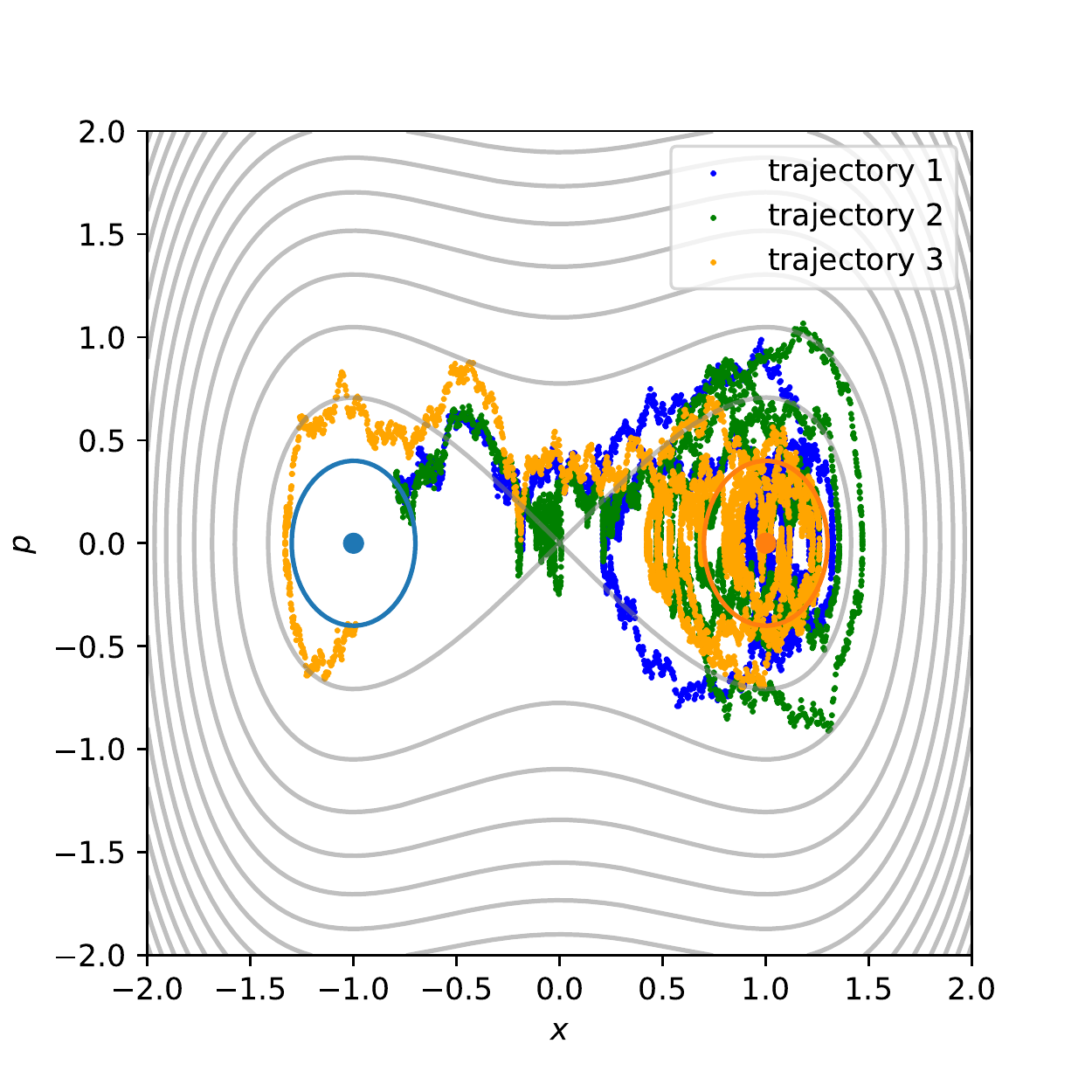}
  \caption{$\epsilon = 0.1$, o/c}
\end{subfigure}\hfill
\begin{subfigure}{.25\linewidth}
\includegraphics[width=1\linewidth]{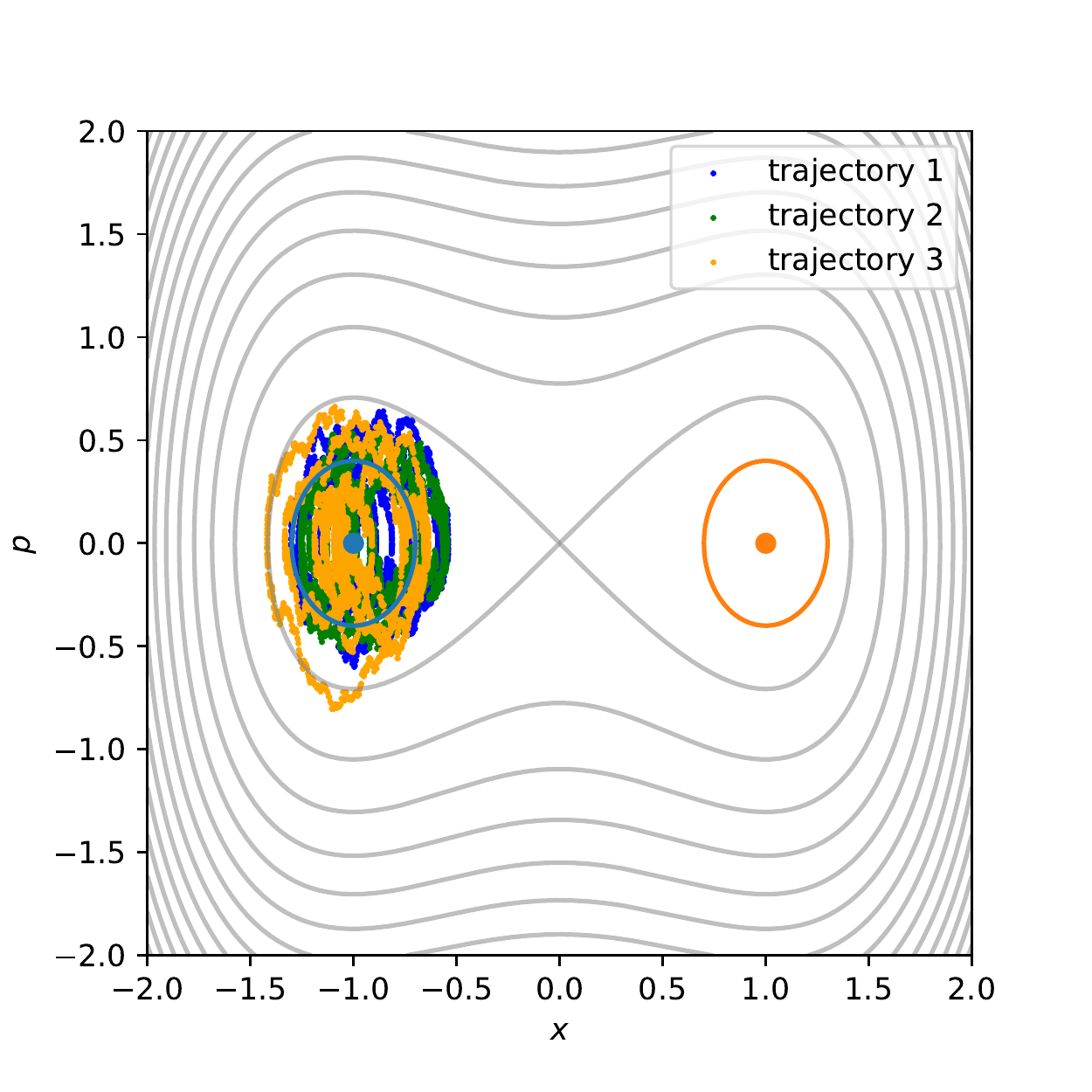}
  \caption{$\epsilon = 0.05$, w/o o/c}
\end{subfigure}\hfill % <-- "\hfill"
\begin{subfigure}{.25\linewidth}
\includegraphics[width=1\linewidth]{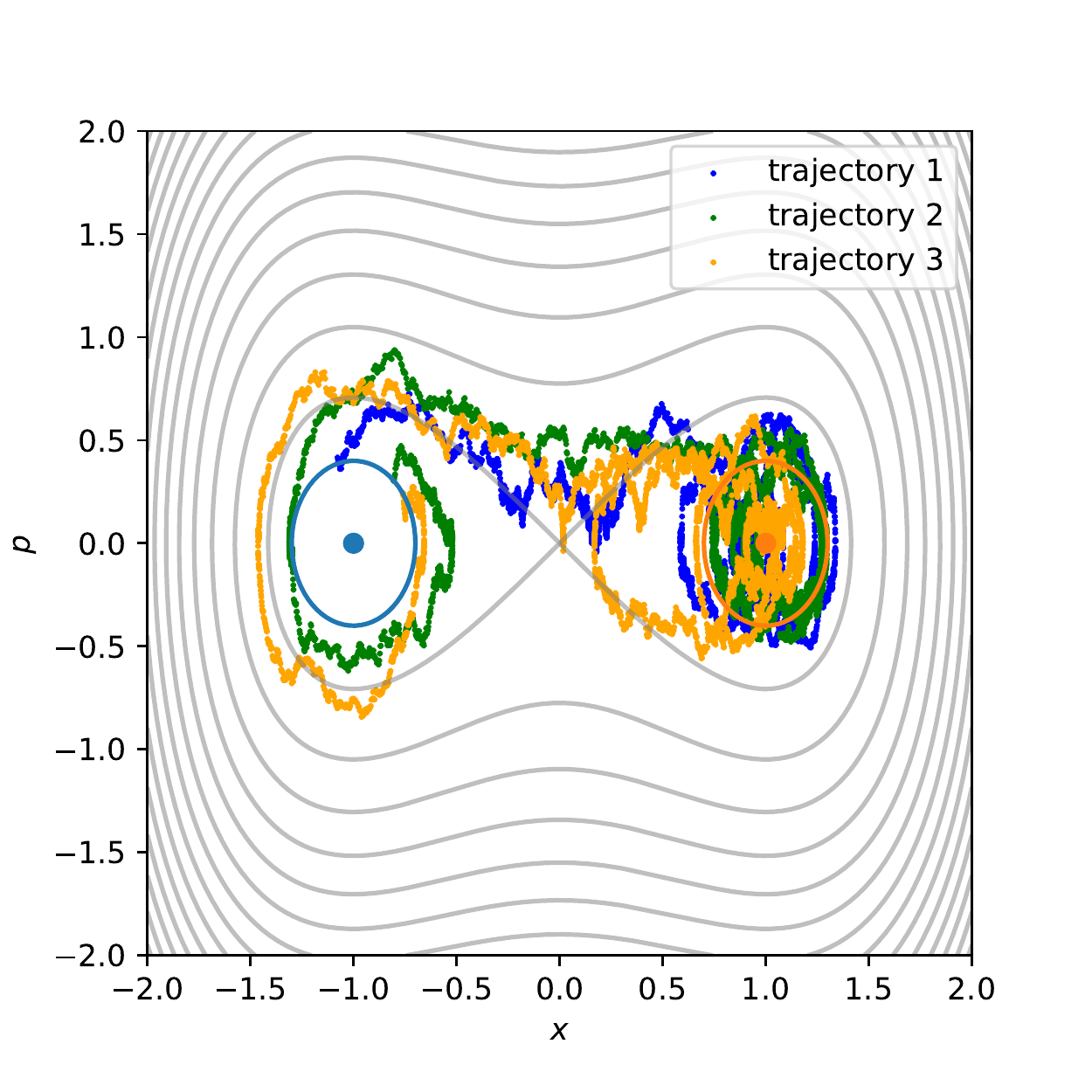}
  \caption{$\epsilon = 0.05$, o/c}
\end{subfigure}
\caption{Comparison of three uncontrolled and controlled processes with the same initial random position for Duffing oscillator at $\epsilon = 0.1$ and $\epsilon = 0.05$. (a) and (c) shows three trajectories under the original dynamics at $\epsilon = 0.1$ and $\epsilon = 0.05$ respectively. (b) and (d) shows three trajectories under the original dynamics for the two cases respectively. 5000 time steps are sampled with $\Delta t= 10^{-3}$. Region $A$ and $B$ are ellipses colored in blue and orange.}
\label{fig:Duffing_traj}
\end{figure}

The optimally controlled process for the Langevin dynamics is governed by
\begin{equation}
\label{eq:Duffing_controlled}
    \begin{cases}
        dX_t = P_t dt \\
        dP_t = - \left[X_t(X_t^2-1) +\tfrac{1}{2} P_t - \epsilon\frac{\partial}{\partial p}\log(q_{\sf nn}^+)\right] dt + \sqrt{ \epsilon } dW.
    \end{cases}
\end{equation}
Fig. \ref{fig:Duffing_traj} shows three sampled trajectories with and without the influence of the optimal control 
starting at the same initial position and the same realizations of the Brownian motion  for $\epsilon = 0.1$
(Fig. \ref{fig:Duffing_traj}(a,b)) and $\epsilon= 0.05$ (Fig. \ref{fig:Duffing_traj}(c,d)). The trajectories governed by the original Langevin dynamics stay near region $A$.  
In contrast, the trajectories governed by the optimally controlled dynamics \eqref{eq:Duffing_controlled} leave $A$ and reach region $B$. 
% { AS: You say some of the same stuff in this paragraph and the caption of the image. This could be consolidated}

% { AS: purely stylistic but I slightly prefer ``o.c.'' instead of ``o/c''. Feel free to ignore this}

\begin{figure}[hbt!]
\begin{center}
\includegraphics[width=0.5\linewidth]{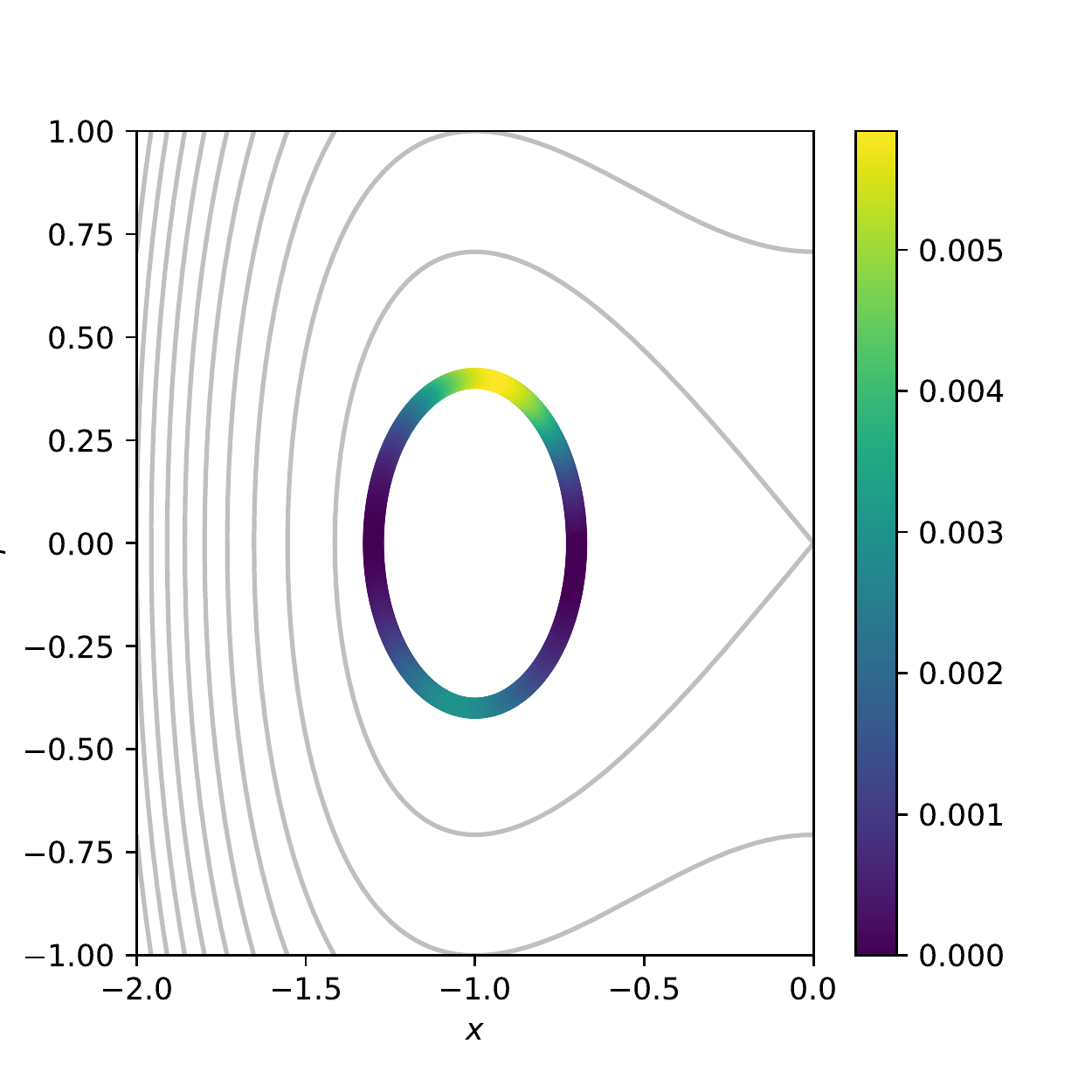}
\caption{The probability distribution of the starting points of transition trajectories on the boundary of $A$ for the Duffing oscillator is shown by color code.  
The gray curves are level sets of Hamiltonian.
}
\label{fig:Duffing_Abdry}
\end{center}
\end{figure}

Next, we find the transition rate $\nu_{AB}$ at $\epsilon = 0.1$ and $\epsilon = 0.05$ in four ways. The results are summarized in Tables \ref{table:Duffing2} and \ref{table:Duffing3}.
\begin{enumerate}
\item \emph{Simulations with optimal control.} The expected crossover time $\mathbb{E}[\tau_{AB}]$ is averaged over 250 trajectories governed by \eqref{eq:Duffing_controlled}. The distribution of the starting points of these trajectories defined by \eqref{distrA} is obtained as described in Section \ref{sec:Mueller_rate}. It is displayed in Fig. \ref{fig:Duffing_Abdry}.
{ The probability $\rho_{AB}$ is found by \eqref{rhoAB} similarly to how it is done in the test problem with Mueller's potential. The PINN committors and the PINN training points are used for Monte Carlo integration. The normalization constant for the invariant density is also found using Monte Carlo integration.}
Then formula \eqref{cross2} is used to find $\nu_{AB}$.

\item \emph{Simulations without optimal control.} Direct simulations of the uncontrolled Langevin dynamics \eqref{Duffing} are used to find $\nu_{AB}$. Ten simulations of $10^7$ time steps each with timestep $\Delta t = 5\times 10^{-3}$ were performed. 
\item \emph{TPT, NN.} The rate $\nu_{AB}$ was found by \eqref{nu5} using the gradient of the committor computed using PINNs.
\item \emph{TPT, FEM.} Likewise, except for the FEM committor was used. 
\end{enumerate}

\begin{table}[h]
    \centering
    \begin{tabular}{|c|c|c|c|c|}
    \hline
        \multicolumn{5}{|c|}{Duffing oscillator $\epsilon = 1/10$}\\
    \hline
         & Simul., o/c & Simul., w/o o/c & TPT, NN & TPT, FEM\\
    \hline
        $\rho_{AB}$ &NA & 4.31e-2 $\pm$ 0.12e-2 & 3.97e-2 & 4.04e-2 \\ % was 4.04e-2
    \hline
        $\mathbb{E}[\tau_{AB}]$ & \textbf{6.88 $\pm$ 0.34} & \textbf{7.32 $\pm$ 0.14}  & NA & NA\\
    \hline
        $\nu_{AB}$ & \textbf{[5.50e-3,6.07e-3]} & \textbf{[5.76e-3,6.01e-3]}  & 4.53e-3 & 5.74e-3 \\ % was 5.74
    \hline
    \end{tabular}
    \caption{Comparison of the estimates for $\rho_{AB}$, $\mathbb{E}[\tau_{AB}]$ and $\nu_{AB}$ for Duffing oscillator at $\epsilon = 0.1$. 
    }
    \label{table:Duffing2}
\end{table}

\begin{table}[h]
    \centering
    \begin{tabular}{|c|c|c|c|c|}
    \hline
        \multicolumn{5}{|c|}{Duffing oscillator $\epsilon = 1/20$}\\
    \hline
         & Simul., o/c & Simul., w/o o/c & TPT, NN & TPT, FEM\\
    \hline
        $\rho_{AB}$ &NA & 4.5e-3 $\pm$ 0.4e-3 & 4.23e-3 & 4.07e-3 \\ % was 4.07e-3
    \hline
        $\mathbb{E}[\tau_{AB}]$ & \textbf{7.34 $\pm$ 0.33} & \textbf{7.48 $\pm$ 0.49}  & NA & NA\\
    \hline
        $\nu_{AB}$ & \textbf{[5.53e-4,6.06e-4]} & \textbf{[5.49e-4,6.51e-4]} & 4.72e-4 & 5.49e-4 \\ % was 5.49e-4 
    \hline
    \end{tabular}
    \caption{Comparison of the estimates for $\rho_{AB}$, $\mathbb{E}[\tau_{AB}]$ and $\nu_{AB}$ for Duffing oscillator at $\epsilon = 0.05$.}
    \label{table:Duffing3}
\end{table}

The results in tables \ref{table:Duffing2} and \ref{table:Duffing3} show that the 95\% confidence intervals for the transition rate $\nu_{AB}$ at $\epsilon = 0.1$ and $\epsilon = 0.05$ estimated by means of simulations with and without optimal control largely overlap. The 95\% confidence intervals expected crossover times $\mathbb{E}[\tau_{AB}]$ largely overlap for $\epsilon = 0.05$ and slightly overlap for $\epsilon=0.1$. The estimates for probability $\rho_{AB}$ that a trajectory at a random time is reactive obtained by direct simulations of uncontrolled dynamics and TPT\&NN and TPT\&FEM are all consistent for $\epsilon = 0.05$ and both TPT-based estimates for $\rho_{AB}$ are smaller than those by direct simulations for $\epsilon=0.1$. 
At both values of $\epsilon$, the TPT estimates for $\rho_{AB}$  obtained using the FEM and PINN committors are consistent, 
while there is a notable discrepancy between the estimates for $\nu_{AB}$ by TPT\&NN and TPT\&FEM. This discrepancy must be caused by the fact that $\nu_{AB}$ uses the gradient of the committors while $\rho_{AB}$ involves the committors themselves as shown in \ref{app:LLexample}. 
% { AS: why must it be caused by this?} 
At both values of $\epsilon$,  the TPT\&FEM estimate for $\nu_{AB}$ falls into the 95\% confidence intervals obtained using simulations, controlled or uncontrolled, while the TPT\&NN seems to underestimate the transition rate.

       %%%%%%%%%%%%%%%%%%%%% 

\subsection{Lennard-Jones-7 in 2D}
\label{sec:LJ7}
Finally, we apply the proposed methodology to estimate the transition rate between the trapezoidal and the hexagonal configurations of the Lennard-Jones-7 cluster (LJ7) in a plane. This is a popular test problem in chemical physics \cite{Dellago1998LJ7,Passerone2001LJ7,Wales2002LJ7,Coifman2008LJ7}. 
In this example, we will compute the committor for the reduced model 2D and use it to construct an approximation to the optimal controller in the original 14D model. The expected crossover time for the original 14D model and the estimate for $\rho_{AB}$ for the 2D model will be used to determine the transition rate. The result will be compared with those obtained via brute force simulations of the original uncontrolled dynamics in 14D. 

\begin{figure}[h!]
\includegraphics[width = 0.9\textwidth]{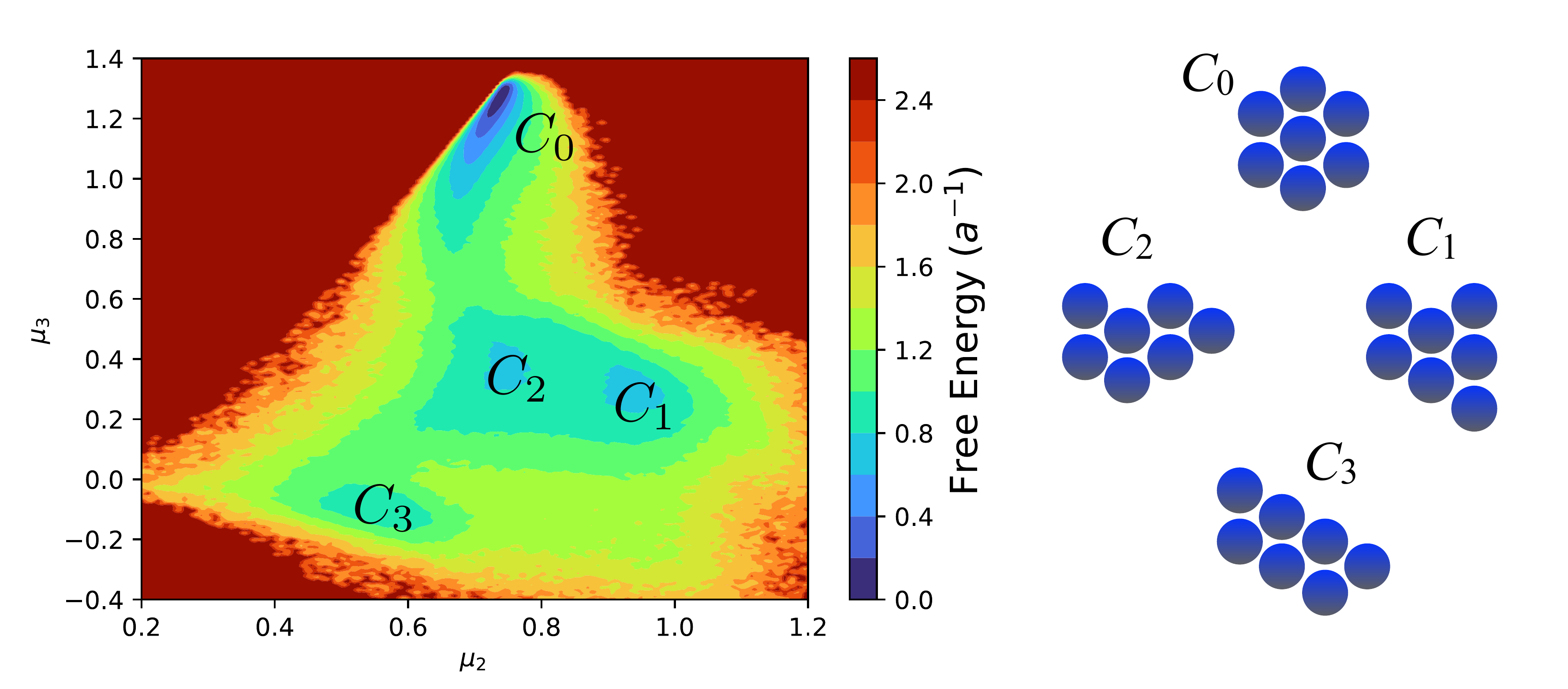}
\caption{The free energy of LJ7 in 2D in CVs second and third central moment of coordination number $\mu_2, \mu_3$ and the four geometrically distinct local minima.}
\label{fig:LJ7main}
\end{figure}

We consider seven two-dimensional particles interacting according to the Lennard-Jones pair potential 
\begin{equation}
   V_{\sf pair}(r) = 4a\left [\left(\frac{\sigma}{r} \right)^{12} - \left(\frac{\sigma}{r} \right)^6 \right]
\end{equation}
where $\sigma >0$ and $a > 0$ are parameters controlling the range and strength of interparticle interaction respectively. We set $\sigma = 1$ and $a = 1$. 
The potential energy of the system 
\begin{equation}
\label{LJ7pot}
    V_{\sf LJ}(x) = \sum_{\substack{i,j = 1\\i < j}}^{7} V_{\sf pair}(\|x_i - x_j\|)
\end{equation}
has four geometrically distinct local minima denoted by $C_0, C_1, C_2, C_3$ shown in Fig. \ref{fig:LJ7main}. 

We assume that the original system is evolving according to the overdamped Langevin dynamics \eqref{sde2}. This choice is dictated by our wish to construct a controller using the committor for the reduced model in collective variables. Since the collective variables are functions only of $x$, the committor for the reduced model lifted to the original phase space will depend only on $x$ and not on the momenta $p$. Therefore, it cannot give an approximation to the optimal control for the Langevin dynamics -- see equation \eqref{sde5c}.
We set $\beta = 5$ as in \cite{evans2021computing}.

%%%%%%%%%%%%%%%%%%%

\subsubsection{The reduced model}
Following  \cite{Tribell02010LJ7,Tsai2019LJ7,evans2021computing}, we pick the 2nd and 3rd central moments of the coordination numbers as the collective variables (CVs) for LJ7.
The coordinate number of particle $i$ is a smooth function approximating the number of nearest neighbors of $i$:
\begin{equation}
    c_i(x) = \sum_{i \neq j} \frac{1 - (\frac{r_{ij}}{1.5\sigma})^8}{1 - (\frac{r_{ij}}{1.5\sigma})^{16}}, \quad r_{ij} = \|x_i - x_j\|.
\end{equation}
The $k$-th central moment of $c_i(x)$ is defined by
\begin{equation}
    \mu_{k}(x) = \frac{1}{7} \sum_{i=1}^7(c_i(x) - \Bar{c}(x))^k, \quad{\rm where}\quad \Bar{c}(x) = \frac{1}{7}\sum_{j=1}^7 c_j(x).
\end{equation}

 The reduced model is governed by the overdamped Langevin dynamics in collective variables \eqref{sde4}
{
\begin{align}
    d\left[\begin{array}{c}\mu_2\\\mu_3\end{array}\right] &= \left[-\nabla M(\mu_2,\mu_3)F(\mu_2,\mu_3) +\beta^{-1}\cdot\nabla M(\mu_2,\mu_3)\right]dt \notag \\
    &+ \sqrt{2\beta^{-1}}M^{1/2}(\mu_2,\mu_3)\left[\begin{array}{c}dw_1\\dw_2\end{array}\right]. \label{eq:LJ7sde}
\end{align}
The corresponding generator is given by \eqref{gen4}.
Fig. \ref{fig:LJ7main} displays the free energy\footnotemark[1] $F(\mu_2,\mu_3)$.    
The diffusion matrix $M(\mu_2,\mu_3)$ varies significantly throughout the accessible free energy region (see Fig. 8 in \cite{evans2021computing}).
The computation of $F(\mu_2,\mu_3)$ and $M(\mu_2,\mu_3)$ is detailed in Appendix A in \cite{evans2021computing}.}
\footnotetext[1]{We thank Luke Evans for sharing with us the free energy and the diffusion matrix in CVs $\mu_2,$ and $\mu_3$.}

Regions $A_{\sf CV}$ and $B_{\sf CV}$ are chosen around minima $C_3$ (the trapezoid) and $C_0$ (the hexagon) respectively.  
We use the subscript to {\sf CV} to indicate that these regions are defined in the set of collective variables. 
Region $A_{\sf CV}$ is a circle centered at $(0.5526,-0.0935)$ of radius $r = 0.1$ and while region $B_{\sf CV}$ is a tilted ellipse defined by the equation 
\begin{equation}
    \frac{\left((x-c_x)\cos \theta + (y - c_y)\sin \theta \right)^2}{r_x^2} + \frac{\left((x-c_x)\sin \theta + (y - c_y)\cos \theta \right)^2}{r_y^2} = 1
\end{equation}
where $(c_x, c_y) = (0.7184,1.1607), r_x = 0.15, r_y = 0.03$ and $\theta = 5\pi/12$. 

%%%%%%%%%%%%%%%%%%%

\subsubsection{Computation of the committor for the reduced model}
\begin{figure}[hbt!]
\begin{subfigure}{.5\linewidth}
\includegraphics[width=1\linewidth]{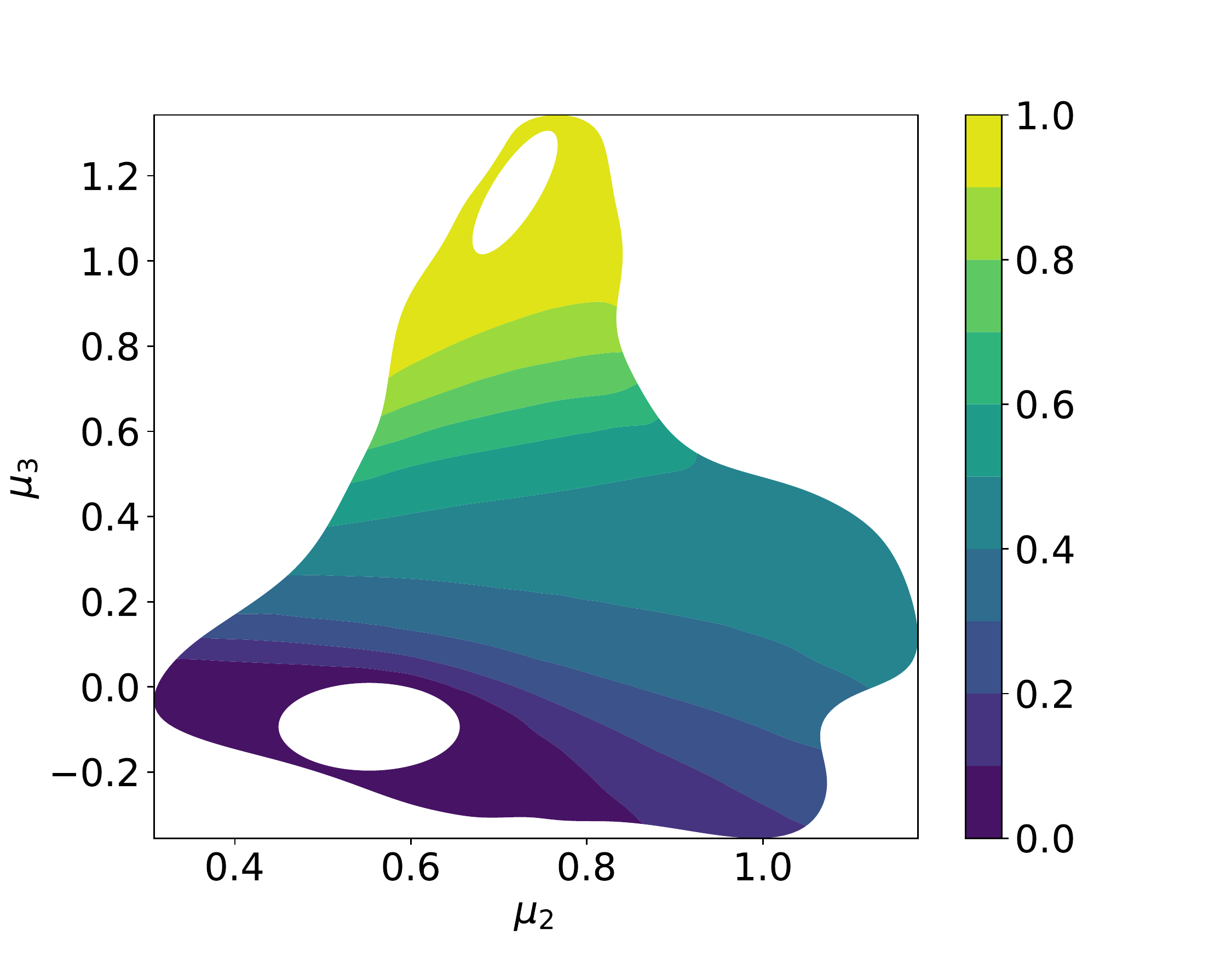}
  \caption{ FEM}
\end{subfigure}\hfill % <-- "\hfill"
\begin{subfigure}{.5\linewidth}
\includegraphics[width=1\linewidth]{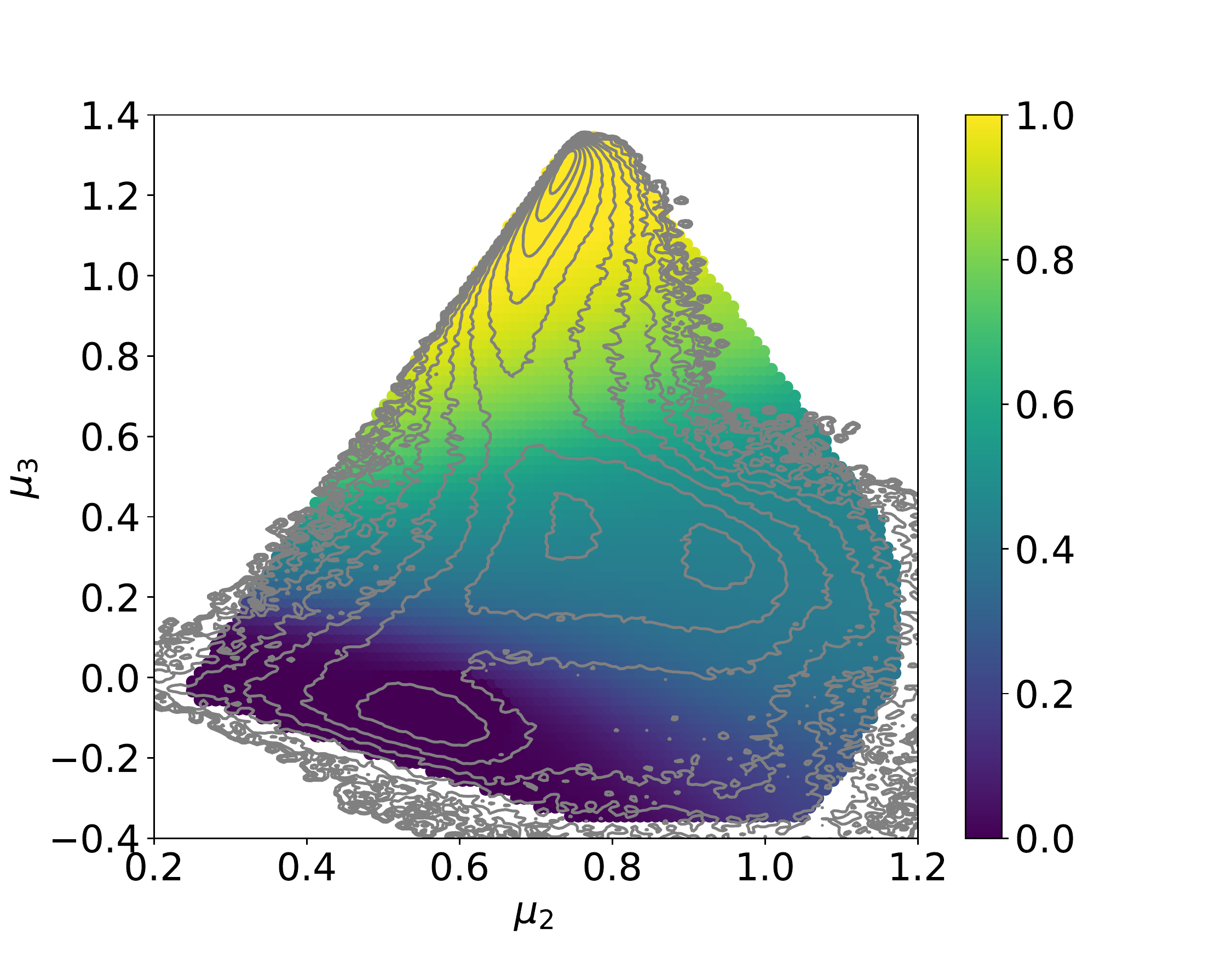}
  \caption{Variational NN}
\end{subfigure}
\caption{The forward committor for LJ7 in 2D computed in the space of collective variables at $\beta = 5$ using
(a) FEM and (b) the variational NN-based solver. The gray curves in (b) are level sets of the free energy.}
\label{fig:LJ7committor}
\end{figure}

The committor for the reduced model is computed in two ways: using FEM and the variational NN-based solver described in Section \ref{sec:NN1}.
For the variational NN, the neural network \eqref{eq:NN1} with $L=2$ hidden layers and $W=10$ neurons per layer has been used { to minimize the loss \eqref{loss1}}. 
The training points were $10^4$ trajectory data projected to the space $(\mu_2,\mu_3)$ 
and assumed to be distributed according to the invariant density { e.g., $\rho \sim e^{\beta F}$ where $F$ is the free energy in $(\mu_2,\mu_3)$} -- see Section 4.2.2 in Ref. \cite{evans2021computing} for more details. { The resulting loss to be minimized hence becomes
\begin{equation}
    {\sf Loss}(\theta) = \frac{1}{K}\sum_{i = k}^K \left[\nabla q(x_k;\theta)^\top M(x_k) \nabla q(x_k;\theta)\right].
\end{equation}}
The results are displayed in Fig. \ref{fig:LJ7committor}. The wMAE and wRMSE are given in Table \ref{table:LJ7_comm}.
% are 1.53e-2 and 2.43e-2 respectively after 2000 epochs of training.

\begin{table}[h]
    \centering
    \begin{tabular}{|c|c|c|c|}
    \hline
    Temperature & NN structure & wMAE & wRMSE  \\
    \hline
    $\beta = 5$ & $L=2$, $W=10$ & 1.5e-2 & 2.4e-2\\
    \hline
    \end{tabular}
    \caption{The errors wMAE and wRMSE of the committor function for LJ7 computed using the variational NN-based solver for the two-dimensional reduced model. The committor computed using FEM is taken as the ground truth.}
    \label{table:LJ7_comm}
\end{table}

%%%%%%%%%%%%

\subsubsection{Estimation of the transition rate using the reduced model and the controlled process}
We set up the controlled process in the original 14-dimensional coordinate space as
\begin{equation}
\label{LJ7controlled}
   dX_t = \left(- \nabla V_{\sf LJ}(X_t) + 2\beta^{-1} \nabla_x \ln q_{\sf nn}(\mu_2(X_t), \mu_3(X_t)) \right)dt + \sqrt{2\beta^{-1}}dW_t,
\end{equation} 
where $V_{\sf LJ}$ is defined by \eqref{LJ7pot} and $q_{\sf nn}$ is the committor computed for the reduced model using the variational NN-based solver. 
The Metropolis-Adjusted Langevin Algorithm (MALA) \cite{Gareth1996_MALA}  with the time step $\Delta t= 5\times 10^{-5}$ has been used for time integration to prevent very large moves 
of the system that can occur due to extremely strong repulsive forces. 

The sets $A$ and $B$ in the original coordinate space $\mathbb{R}^{14}$ are defined by lifting the sets $A_{\sf CV}$ and $B_{\sf CV}$:
\begin{align*}
%\label{LJ7_AB}
A&:=\{x\in\mathbb{R}^{14}~|~(\mu_2(x),\mu_3(x))\in A_{\sf CV}\},\\
 B&:=\{x\in\mathbb{R}^{14}~|~(\mu_2(x),\mu_3(x))\in B_{\sf CV}\}.
\end{align*}

Three trajectories of the controlled process \eqref{LJ7controlled} and three trajectories of the uncontrolled overdamped Langevin dynamics \eqref{sde2} in the 14D with the same three realizations of the Brownian motion projected to the space $(\mu_2,\mu_3)$ are displayed in Fig. \ref{fig:LJ7_traj} (a) and (b) respectively.

\begin{figure}[hbt!]
\begin{center}
\includegraphics[width=0.5\linewidth]{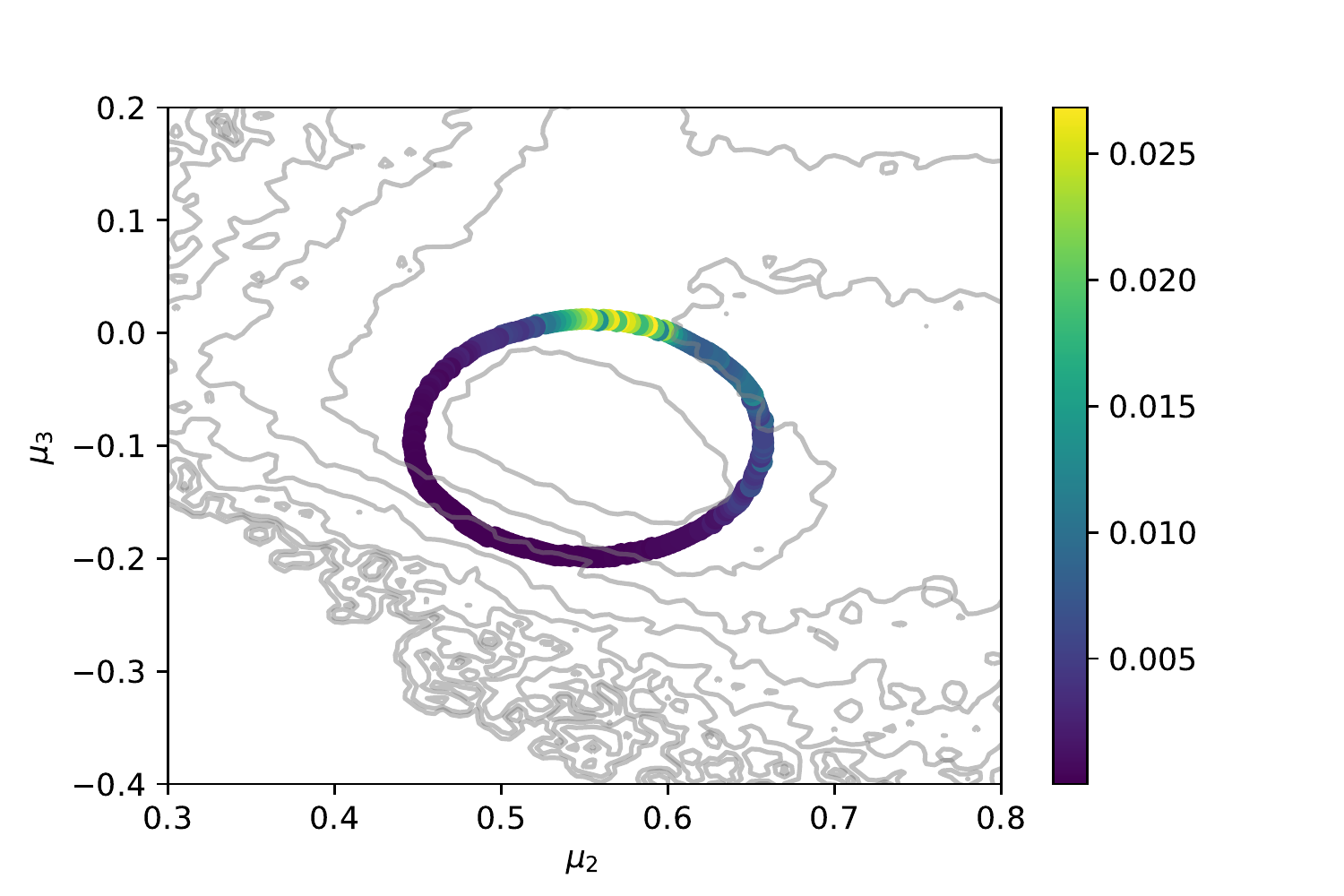}
\caption{The probability distribution of the starting points of transition trajectories on the boundary of $A$ for LJ7  is shown by color code.  
The gray curves are level sets of the free energy.}
\label{fig:LJ7_Abdry}
\end{center}
\end{figure}

\begin{figure}[hbt!]
\begin{subfigure}{.5\linewidth}
\includegraphics[width=1.0\linewidth]{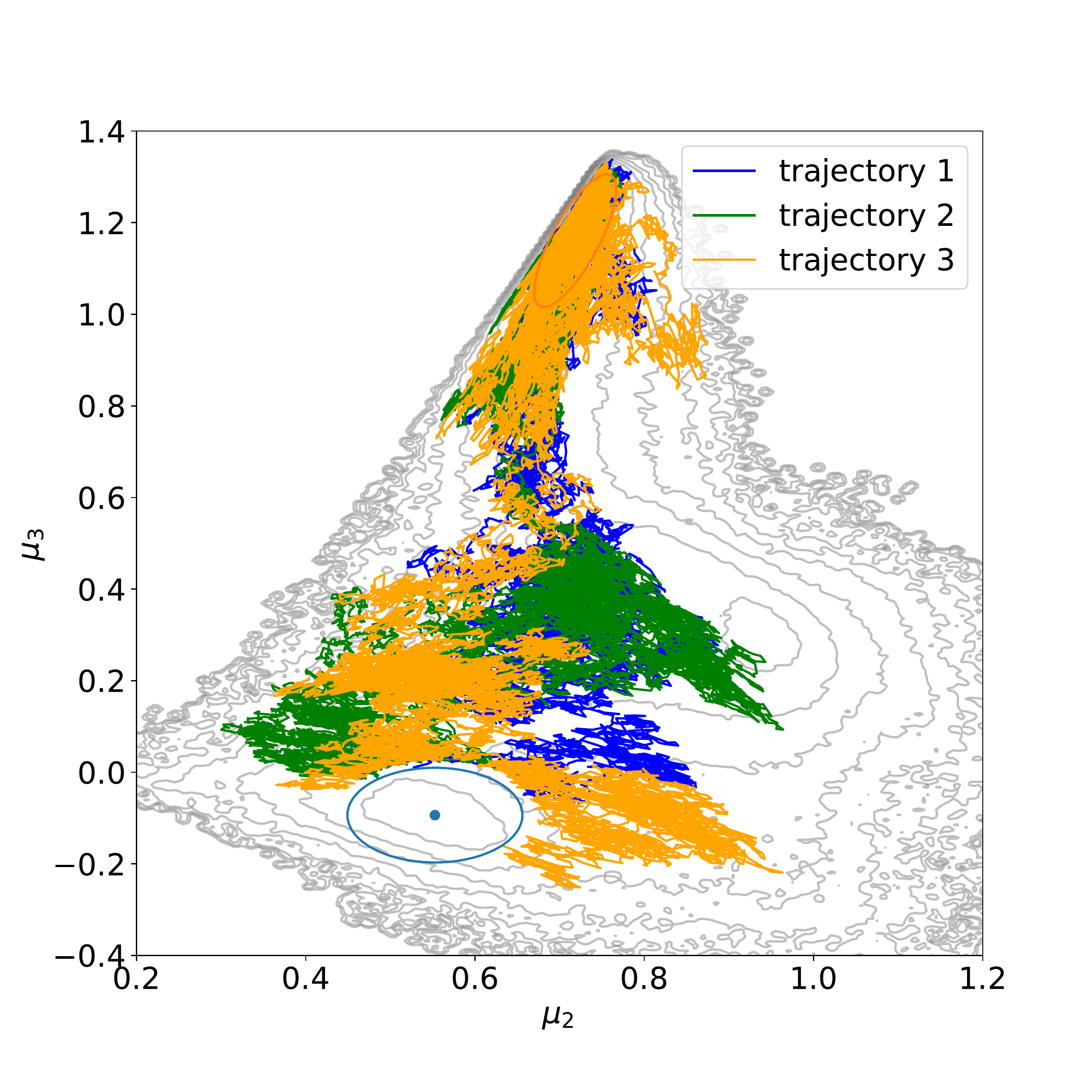}
  \caption{$\beta = 5$, o/c}
\end{subfigure}\hfill % <-- "\hfill"
\begin{subfigure}{.5\linewidth}
\includegraphics[width=1.0\linewidth]{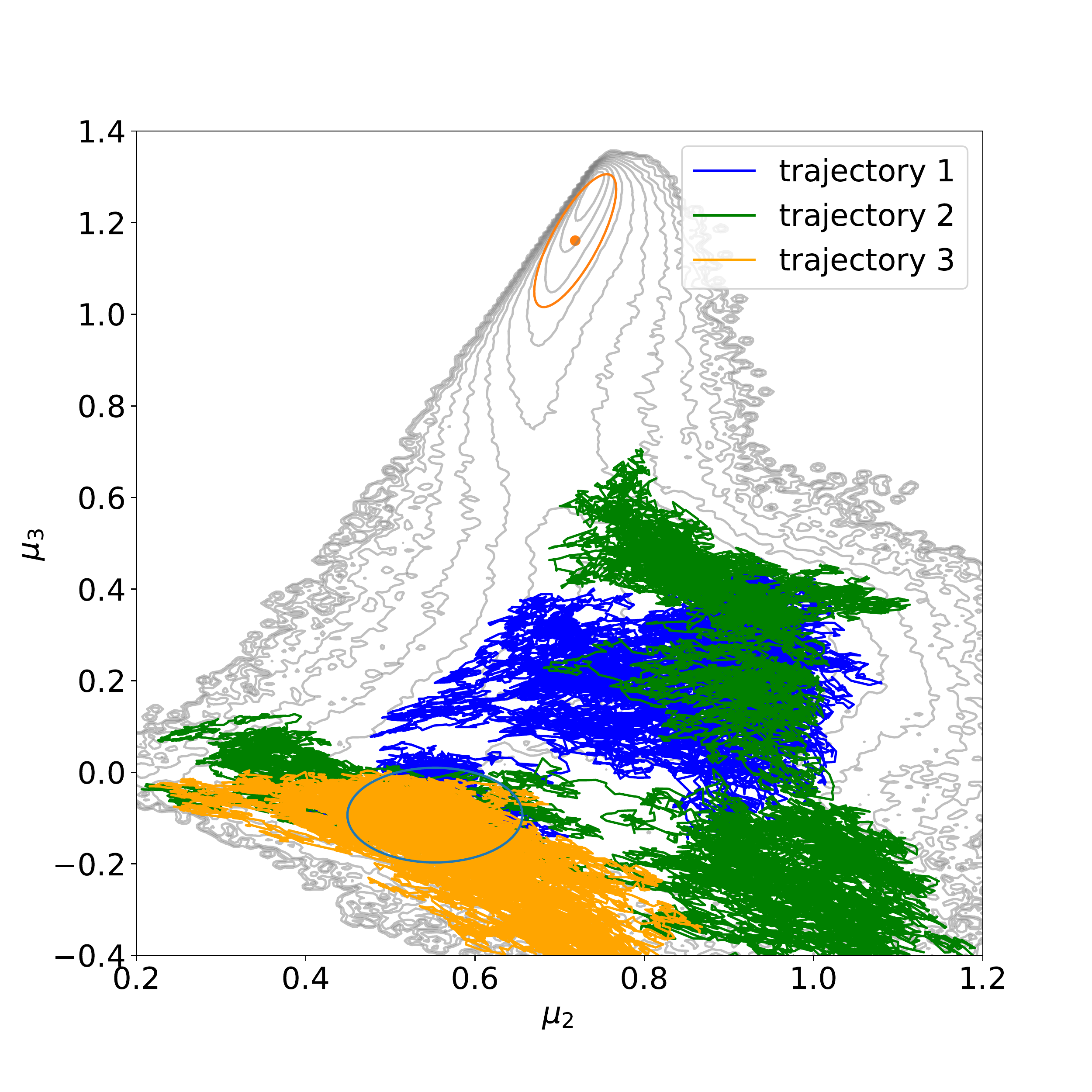}
  \caption{$\beta = 5$, w/o o/c}
\end{subfigure}
\caption{Three trajectories of (a) the controlled process \eqref{LJ7controlled} and (b) the uncontrolled process projected into the space of collective variables, i.e., the overdamped Langevin dynamics \eqref{sde2} in the 14D projected into the space of collective variables with the same starting points and the same realizations of Brownian motion.}
\label{fig:LJ7_traj}
\end{figure}

The expected crossover time $\mathbb{E}[\tau_{AB}]$ is averaged over 254 trajectories of SDE \eqref{LJ7controlled}. 
The starting points for these trajectories near the boundary of $A$ are sampled as follows. 
First, the probability weights of the points on the boundary of $A_{\sf CV}$ are computed as described in Section \ref{sec:Mueller_rate} (see Fig. \ref{fig:LJ7_Abdry}). 
Then, points on $\partial A_{\sf CV}$ are sampled according to these weights and lifted to the original coordinate space by running biased simulations as described in Appendix A of \cite{evans2021computing} (see equations (A.3) and (A.4) there).

We compute the transition rate $\nu_{AB}$ in four ways analogous to those for the Duffing oscillator (Section \ref{sec:Duffing_rate}).
\begin{enumerate}
\item \emph{Simulations, optimal control, 14D.} 
The rate is calculated using \eqref{cross2}. 
The expected crossover time is estimated as described above using the controlled process  \eqref{LJ7controlled} in 14D.
The probability $\rho_{AB}$ is obtained for the reduced 2D model using \eqref{rhoABfromNN}--\eqref{ZVnn} with the committor computed by the variational NN-based solver 
and a uniform grid of points rather than the actual training points for the neural network.
\item \emph{Simulations without optimal control, 14D.} Ten runs of direct simulations of the uncontrolled overdamped Langevin dynamics \eqref{sde2} in 14D 
of $10^8$ timesteps with $\Delta t = 5\cdot10^{-5}$ were executed. The numbers of transitions from $A$ to $B$ that occurred in these runs were  28, 138, 160, 146,  93,  63, 158, 165, 171, and 160. 
% { AS: don't think we need to write out all of these. Perhaps we could just give the average?} \jx{Resolved}
\item \emph{TPT, NN.} The rate $\nu_{AB}$ was found by \eqref{nuAB} for the reduced 2D model using the gradient of the committor computed using the variational NN-based solver.
\item \emph{TPT, FEM.} Likewise, except for the FEM committor was used. 
\end{enumerate}
The results are displayed in Table \ref{table:LJ7rates}. 
The following observations can be made.
\begin{enumerate}
\item
For the transition rates obtained using controlled and uncontrolled simulations in the 14D, the relative error is about 14\%. There is a large overlap of the 95\%  confidence intervals. 
\item
There is a discrepancy in the expected crossover time $\mathbb{E}[\tau_{AB}]$ computed using controlled and uncontrolled simulations in the 14D: the expected crossover time for the controlled process exceeds that for the uncontrolled process by approximately 35\%. 
This discrepancy is primarily caused by the fact that the controller obtained by lifting the committor computed for the reduced model is not optimal. It still drives the trajectories away from the boundary of $A$ but, unlike the optimal controller, somewhat affects the statistics for the transition trajectories.
\item 
The rates obtained using the reduced 2D model are highly exaggerated (by the factor of approximately four) as one can expect given the Zhang-Hartmann-Schuette rate formula  \eqref{ZHS2016}. 
Indeed, the collective variables $\mu_2$ and $\mu_3$ are chosen due to their ability to separate the four geometrically distinct local minima of LJ7 while there is no indication that they are supposed to represent the dynamics accurately. 
\item 
The values of $\rho_{AB}$ and $\nu_{AB}$ computed for the reduced 2D model using the FEM and the variational NNs are in good agreement with each other.
\end{enumerate}

\begin{table}[h]
    \centering
    \begin{tabular}{|c|c|c|c|c|}
%    \hline
%        \multicolumn{5}{|c|}{LJ7 $\beta = 5$}\\
    \hline
         & Simul., o/c & Simul., w/o o/c & TPT, NN & TPT, FEM\\
         \hline
         Dimension & 14D & 14D & 2D & 2D\\
    \hline
%    $\rho_{A}$ &NA & 0.279 $\pm$ 0.078 & 0.328 & 0.348 \\
%
%    \hline
        $\rho_{AB}$ &NA & 0.080 $\pm$ 0.024 & 0.106 & 0.108 \\ % 0.108
    \hline
        $\mathbb{E}[\tau_{AB}]$ & 4.88 $\pm$ 0.48 & 3.16 $\pm$ 0.21 & NA & NA\\
    \hline
        $\nu_{AB}$ & {\bf 0.022},~{\bf [0.020,\,0.024] }&{\bf 0.025},~ {\bf [0.019,\,0.033] }& 0.097 & 0.086 \\
        % was 0.092 for FEM TPT
    \hline
    \end{tabular}
    \caption{Comparison of the estimated $\rho_{AB}$, $\mathbb{E}(\tau_{AB})$ and $\nu_{AB}$ for LJ7 at $\beta = 5$. 
    }
    \label{table:LJ7rates}
\end{table}

%%%%%%%%%%%%%%%%%%%%%%%%%%%%%%%%%%%%%%%%%%%%
%%%%%%%%%%%%%%%%%%%%%%%%%%%%%%%%%%%%%%%%%%%%
%%%%
%%%%. C O N C L U S I O N

\section{Conclusion}
\label{sec:conclusion}
In this work, we have proposed a methodology for sampling transition trajectories and estimating transition rates in systems governed by SDEs using optimal control and perhaps model reduction.

Our main theoretical contribution is the proof of Theorem \ref{thm:T1} establishing the optimality of the control obtained from the committor via the Doob $h$-transform  for a broad class of processes including the Langevin dynamics and the overdamped Langevin dynamics in collective variables. 

We have elaborated on a number of practical aspects related to the use of neural network-based solvers, finite element methods, and sampling reactive trajectories.

We have conducted in-depth case studies of three benchmark systems.
In particular, we have demonstrated that the optimal control and the estimate for the probability of a trajectory to be reactive 
at a random moment of time obtained for the reduced model result is a reasonably good estimate of the transition rate even 
if the collective variables do not represent the dynamics accurately. 

Our codes are published on GitHub \cite{margotyjx,mar1akc}.

Further improvement of the proposed methodology can be done in the following two directions. First, the design of collective variables is important for an 
accurate representation of the dynamics. Autoencoders (see e.g. \cite{Belkacemi2021} and references therein) 
with an appropriate choice of the loss function seem to be a promising tool. 
Second, the neural network-based techniques for solving the committor problem are promotable to higher dimensions. 
% { AS: Are there any computational difficulties of promoting these to higher dimensions?} \jx{Resolved}
In this work, we intentionally calculated all required quantities for the use of the transition path theory without meshing the space. 
We did not attempt, though, to use these techniques in higher dimensions.
We are leaving these research topics for future work.

  \section{Acknowledgements}
  We thank Dr. Luke Evans for providing us with the free energy and diffusion matrix data for the Lennard-Jones-7 test problem.
 We also thank the UMD REU students Luke Triplett, Dmitry Pinchuk, Prisca Calkins, and William Clark for their investigation into neural network-based committor solvers. This work was partially supported by AFOSR MURI grant FA9550-20-1-0397 and by 
  NSF REU grant DMS-2149913.

 \appendix
\setcounter{equation}{0}
\renewcommand{\theequation}{\Alph{section}-\arabic{equation}}
%\section*{Appendix A}

\section{ Proof of equation \eqref{cross1}:  $\nu_{AB} = {\rho_{AB}}/{\mathbb{E}[\tau_{AB}]}$}
\label{sec:app0}
\begin{proof}
Let $X_t$, $0 \leq t \leq T$, be a long trajectory. We decompose the interval $[0,T]$ into two subsets $[0,T] = \mathcal{I}_A\cup\mathcal{I}_B$ where
\begin{equation}
    \begin{cases}
    \mathcal{I}_A([0,T]):=\{t\in[0,T] ~|~\tau_A^-(X_t) <\tau_B^-(X_t) \},\\
    \mathcal{I}_B([0,T]):= \{t\in[0,T] ~|~\tau_B^-(X_t) <\tau_A^-(X_t)\}.
    \end{cases}
\end{equation}
In words, $\mathcal{I}_A([0,T])$ is the set of all moments of time $t$ in the interval $[0,T]$ such that the trajectory at time $t$, $X_t$, last visited $\bar{A}$ rather than $\bar{B}$. 
The set $\mathcal{I}_B([0,T])$ is described likewise.
Let $T_A$ and $T_B$ be the total lengths of $\mathcal{I}_A([0,T])$ and $\mathcal{I}_B([0,T])$ respectively. 
The probabilities $\rho_A$ and $\rho_B$ that the trajectory at a randomly picked moment of time last visited $A$ or $B$ are, respectively,
\begin{equation}
\rho_A =\lim_{T\rightarrow\infty} \frac{T_A}{T},\quad \rho_B = \lim_{T\rightarrow\infty}\frac{T_B}{T}.
\end{equation}
%Therefore, for large enough $T$, we can approximate $T_A$ and $T_B$ as 
%\begin{equation}\label{eq: TA}
%    T_A = \lim_{T\rightarrow\infty} \rho_A T, \quad T_B = \lim_{T\rightarrow\infty} \rho_B T = \lim_{T\rightarrow\infty} (1- \rho_A) T.
%\end{equation}
The set $\mathcal{I}_A([0,T])$ is further decomposed into two subsets $\mathcal{I}_A([0,T])=\mathcal{I}_{AA}([0,T])\cup\mathcal{I}_{AB}([0,T])$ of total lengths $T_{AA}$ and $T_{AB}$ respectively where
\begin{equation}
    \begin{cases}
    \mathcal{I}_{AA}([0,T]):=\{t\in[0,T] ~|~ \tau_A^-(X_t) <\tau_B^-(X_t)~ \&~ \tau_A^+(X_t) <\tau_B^+(X_t) \},\\
    \mathcal{I}_{AB}([0,T]):= \{t\in[0,T]~|~ \tau_A^-(X_t) <\tau_B^-(X_t) ~\&~ \tau_B^+(X_t) <\tau_A^+(X_t) \}.
    \end{cases}
\end{equation}
I.e., $ \mathcal{I}_{AA}([0,T])$ is the set of moments of times $t\in[0,T]$ such that the trajectory at time $t$, $X_t$, last visited $\bar{A}$ rather than $\bar{B}$ and going to hit next $\bar{A}$ rather than $\bar{B}$, while $ \mathcal{I}_{AB}([0,T])$ is the subset of moments of time $t\in[0,T]$ such that the trajectory is reactive. Respectively, the probability $\rho_{AA}$ that a trajectory at a randomly picked time $t$ last hit $\bar{A}$ rather than $\bar{B}$ and is not reactive and the probability $\rho_{AB}$ that a trajectory at a randomly picked time $t$ is reactive are given by
\begin{equation}
    \rho_{AA} = \lim_{T\rightarrow\infty} \frac{T_{AA}}{T}, \quad
    \rho_{AB} =\lim_{T\rightarrow\infty} \frac{T_{AB}}{T}.
\end{equation}
Now we recall the definitions of the transition rate $\nu_{AB}$ and the expected crossover time $\mathbb{E}[\tau_{AB}]$:
\begin{equation}
\nu_{AB} =  \lim_{T\rightarrow\infty} \frac{N_{AB}}{T},\quad
\mathbb{E}[\tau_{AB}] =  \lim_{T\rightarrow\infty} \frac{T_{AB}}{N_{AB}}.
\end{equation}
Hence the expected crossover time $\mathbb{E}[\tau_{AB}]$ can be written as
\begin{equation}
\nu_{AB} =  \lim_{T\rightarrow\infty} \frac{N_{AB}}{T} =   \lim_{T\rightarrow\infty} \frac{N_{AB}}{T_{AB}}\frac{T_{AB}}{T}  =\frac{\rho_{AB}}{\mathbb{E}[\tau_{AB}]}.
\end{equation}
\end{proof}

%%%%%%%%%%%%%%%%%%%%%%%%%

 \section{Proof that $-\nabla\cdot J_{AB} = 0$ in $\Omega_{AB}$}
 \label{sec:appA}
 Let us show that the divergence of the reactive current, or, equivalently, the stationary current of the transition path process, vanishes in $\Omega_{AB}$. 
 We will need a formula for the divergence of a matrix-vector product. It can be checked directly that for any $A\in\mathbb{R}^{d\times d}$ and any $y\in\mathbb{R^d}$,
 \begin{equation}
 \label{aux1}
 \nabla\cdot(Ay) = {\sf div}A\cdot y + {\sf tr}\left(A\nabla y\right),\quad{\rm where}\quad 
 \nabla y = \left[\begin{array}{c}\nabla y_1^\top\\\vdots\\\nabla y_d^\top\end{array}\right].
 \end{equation}
 Using \eqref{aux1} we calculate:
 \begin{align*}
\nabla\cdot J_{AB} & = -\nabla\cdot\left[ \left(b\mu -\frac{1}{2}{\sf div}\left(\sigma\sigma^\top \mu \right)\right)q^+q^- 
+\frac{1}{2}\mu\sigma\sigma^\top\left(q^-\nabla q^+  - q^+\nabla q^-\right)\right] \\
& = q^+q^-\nabla \cdot  \left(b\mu -\frac{1}{2}{\sf div}\left(\sigma\sigma^\top \mu \right)\right) \\
&+ \left(b\mu -\frac{1}{2}{\sf div}\left(\sigma\sigma^\top \mu \right)\right) \cdot  \left(q^+\nabla q^- + q^-\nabla q^+\right) \\
&+
\frac{1}{2}{\sf div}\left(\sigma\sigma^\top \mu\right)\cdot \left(q^+\nabla q^- - q^-\nabla q^+\right)
+\frac{1}{2}{\sf tr}\left(\sigma\sigma^\top \nabla\nabla q^+\right)\mu q^- + \\
 &+\frac{1}{2}{\sf tr}\left(\sigma\sigma^\top \mu\nabla q^-[ \nabla q^+]^\top\right)
 -\frac{1}{2}{\sf tr}\left(\sigma\sigma^\top \mu\nabla q^+[ \nabla q^-]^\top\right)
 -\frac{1}{2}{\sf tr}\left(\sigma\sigma^\top \nabla\nabla q^-\right)\mu q^+\\
 &= q^+q^- \nabla\cdot J + \mu q^-\mathcal{L}q^+ - \mu q^+\mathcal{L}^\dagger q^- = 0.
 \end{align*}
In the last expression, $J$ is the stationary current for the invariant density $\mu$, and hence 
$\mathcal{L}^*\mu = -\nabla\cdot J = 0$ in $\Omega$. The last two terms are zero as by \eqref{eg:commproblem}. 

%%%%%%%%%%%%%%%%%%%%%%%%%%%%
%%%%%%%%%%%%%%%%%%%%%%%%%%%%
\section{Proof of Theorem \ref{thm:T1}}
\label{sec:AppB}
\begin{proof}
This proof combines ideas from Gao et al. (\cite{Gao2021_OptimalControl}, the proof of Theorem 3.3) and L.~C. Evans's notes on the control theory \cite{LCEvans_notes} 

{\bf Step 1. Regularization.}
We first consider a regularized optimal control problem in which the exit cost \eqref{eq:exitcost} is replaced with a finite exit cost
\begin{equation}
\label{eq:exitcost1}
        g_N(x) = 
        \begin{cases}
           N, & x \in \overline{A} \\
            0, & x \in \overline{B}
        \end{cases},\quad {\rm where}~N~\text{is a large number}.
\end{equation}
Let $c_N^{\ast}(x)$ be the infimum of the cost functional $C_x[N;v(\cdot)]$ with the regularized exit cost \eqref{eq:exitcost1}  among all admissible controls. 
Note that the admissible set $\mathcal{V}$ is not empty because the time $\tau_{AB}<\infty$ 
almost surely since the system is ergodic and the domain $\Omega$ is compact.  
Furthermore, $c_N^{\ast}(x) <\infty$ as $v\equiv 0$ is an admissible control and the
corresponding cost functional is $C_x[N;v(\cdot)] = N(1-q^+(x)) < \infty$.  
Indeed, if $v\equiv 0$, then
the process hits $\partial B$ first with probability $q^+(x)$ and scores zero and hits $\partial A$ first with probability $1-q^+(x)$ and scores $N$.

We also define a regularized forward committor as the solution to the
following boundary-value problem
    \begin{equation}
    \label{commN}
    \begin{cases}
    \mathcal{L}q_N^+ = 0,& x\in\Omega_{AB},\\
    q_N^+ = e^{-N},&x\in\partial A\\
    q_N^+ = 1,&x\in\partial B\\
    \frac{\partial q_N^+}{\partial\hat{n}} = 0,&x\in\partial \Omega
    \end{cases},
    \end{equation}
    where $\mathcal{L}$ is the generator for \eqref{sdeB1}. It is easy to check that
    \begin{equation}
    \label{commNcomm}
    q_N^+ = q^+ + (1-q^+)e^{-N}.
    \end{equation}
    
%%%%
{\bf Step 2. Show that $c_N^* \ge -\log q_N^+$.}
The regularized forward committor can be written as 
\begin{equation}
\label{commN1}
q_N^+(x) = \mathbb{E}_{P}\left[e^{-g_N(X_{\tau_{AB}})}~|~X_0 = x\right] \equiv  \mathbb{E}_{P,x}\left[e^{-g_N(X_{\tau_{AB}})}\right].
\end{equation}
Indeed, the process $X_t$ governed by \eqref{sdeB1} with $X_0=x$  reaches $\partial A$ at time  the stopping time $\tau_{AB}$ with probability $(1-q^+)$ and scores $e^{-N}$, and reaches $\partial B$ at $\tau_{AB}$ and scores $1$. This results in the expectation given by the right-hand side of \eqref{commNcomm} which is equal to $q_N^+(x)$.

Let $Y_t$ be the controlled process governed by \eqref{eq:sde1c} with a control $\sigma v$, $v\in C^1(\Omega_{AB}$,  satisfying Novikov's condition
\begin{equation}
\label{Novikov}
\mathbb{E}_P\left[\exp\left(\int_0^{\tau{AB}}\frac{1}{2}\|\sigma^T(Y_s) v(Y_s)\|^2ds\right)\right] < \infty
\end{equation}
and $P_v$ be the probability measure on the path space of this process.
%The admissibility condition \eqref{admiss} is Novikov's condition in the Girsanov theorem (Theorem 8.6.5, p.158 in \cite{Oksendal}).
According to the Girsanov theorem  (Theorem 8.6.5, p.158 in \cite{Oksendal}),
\begin{equation}
\label{Gir1}
\mathbb{E}_{P,x}\left[e^{-g_N(X_{\tau_{AB}})}\right] = \mathbb{E}_{P_v,x}\left[e^{-g_N(Y_{\tau_{AB}})}\right]  = \mathbb{E}_{P,x}\left[e^{-g_N(Y_{\tau_{AB}})}\frac{dP_v}{dP}\right],
\end{equation}
where the Radon-Nikodym derivative
\begin{equation}
\label{RN0}
\frac{dP_v}{dP} = \exp\left\{-\int_0^{\tau_{AB}}\sigma^\top v(Y_s)\cdot dW_s - \frac{1}{2}\int_0^{\tau_{AB}}\|\sigma^\top v(Y_s)\|^2ds\right\} > 0~~{\rm P-a.s.}
\end{equation}
Therefore,
\begin{equation}
\label{RN1}
\mathbb{E}_{P,x}\left[e^{-g_N(X_{\tau_{AB}})}\right]  = \mathbb{E}_{P,x}\left[e^{-g_N(Y_{\tau_{AB}}) 
-\int_0^{\tau_{AB}}\sigma^\top v(Y_s)\cdot dW_s -
 \frac{1}{2}\int_0^{\tau_{AB}}
\|\sigma^\top v(Y_s)\|^2ds  }\right].
\end{equation}
By Jensen's inequality, for any smooth convex function $\phi$ and a random variable $Z$ we have $\phi(\mathbb{E}[Z])\le \mathbb{E}[\phi(Z)]$. 
Applying it to the right-hand side of \eqref{RN1} we get
\begin{equation}
\label{RN2}
e^{-\mathbb{E}_{P,x}\left[  g_N\left(Y_{\tau_{AB}}\right) +
\int_0^{\tau_{AB}}\sigma^\top v(Y_s)dW_s +
 \frac{1}{2}\int_0^{\tau_{AB}}
\|\sigma^\top v(Y_s)\|^2ds\right]} \le \mathbb{E}_{P,x}\left[e^{-g_N\left(X_{\tau_{AB}}\right)}\right].
\end{equation}
Since the expectation of the Ito stochastic integral is zero, i.e.,
\begin{equation*}
\mathbb{E}_{P,x}\left[ \int_0^{\tau_{AB}}\sigma^\top v(Y_s)\cdot dW_s \right] = 0,
\end{equation*}
and the cost functional $C_x[N;v] $ defined in \eqref{eq:costfunction} is exactly 
\begin{equation*}
C_x[N;v] \equiv \mathbb{E}_{P,x}\left[  g_N\left(Y_{\tau_{AB}}\right) +
 \frac{1}{2}\int_0^{\tau_{AB}}
\|\sigma^\top v(Y_s)\|^2ds\right],
\end{equation*}
and recalling \eqref{commN1} we get
\begin{equation}
\label{RN3}
e^{-C_x[N;v]} \le \mathbb{E}_{P,x}\left[e^{-g_N\left(X_{\tau_{AB}}\right)}\right] \equiv q_N^+(x).
\end{equation}
Taking logarithms of the left- and right-hand side of \eqref{RN3} and multiplying the result by $-1$ we obtain the following lower bound for the cost functional:
for any control $\sigma^\top v$ satisfying Novikov's condition \eqref{Novikov},
\begin{equation}
\label{RN4}
C_x[N;v] \ge -\log q_N^+(x).
\end{equation}
Since the admissible set $\mathcal{V}$ is closed, the bound \eqref{RN4} holds for any admissible $v$. This means that for any $v\in\mathcal{V}$,
\begin{equation}
\label{RN5}
c^*_N(x)  \ge -\log q_N^+(x).
\end{equation}

%%%
{\bf Step 3. Derive the Hamilton-Jacobi-Bellman equation for the minimal cost $c^*_N$.} 
This upper bound will be derived via the Hamilton-Jacobi-Bellman equation. 
Let $v(x)\in C^1(\Omega_{AB})$ be such that Novikov's condition \eqref{Novikov} holds and let $h$ be a small positive number. 
Then for the process $Y_t$ governed by the controlled SDE \eqref{eq:sde1c} with the control $\sigma^\top v$ we have the following upper bound:
\begin{equation}
\label{ineq0}
c_N^{\ast}(x) \le  \mathbb{E}_{P,x}\left[\frac{1}{2} \int_0^{h\land\tau_{AB}}\|\sigma^T(Y_s) v(Y_s)\|^2 ds + c_N^{\ast}(Y_{h\land\tau_{AB}})\right].
\end{equation}
The equality is reached if $v$ is an optimal controller. We observe that if $v\equiv 0$ than $C_x[N;0] = N(1-q^+(x)) < N$ for $x\in\Omega_{AB}$.
Therefore, $c_N^{\ast}(x)  \le C_x[N;0] <N<\infty$. Therefore, we subtract $c_N^{\ast}(x) $ from both sides of the inequality \eqref{ineq0} and get 
\begin{align}
\label{ineq1}
0 & \le  \mathbb{E}_{P,x}\left[\frac{1}{2} \int_0^{h\land\tau_{AB}}\|\sigma^T(Y_s) v(Y_s)\|^2 ds + c_N^{\ast}(Y_{h\land\tau_{AB}})- c_N^{\ast}(x) \right] \\
& =  \mathbb{E}_{P,x}\left[\frac{1}{2} \int_0^{h\land\tau_{AB}}\|\sigma^T(Y_s) v(Y_s)\|^2 ds\right] +\mathbb{E}_{P,x}\left[ c_N^{\ast}(Y_{h\land\tau_{AB}}) \right]- c_N^{\ast}(x).
\end{align}
Dividing by $h$ and letting $h\rightarrow 0$ we obtain:
\begin{equation} 
\label{ineq2}
0 \le \frac{1}{2}\|\sigma^\top v(x)\|^2 + \lim_{h\rightarrow 0} \frac{ \mathbb{E}_{P,x}\left[c_N^{\ast}(Y_h) \right] - c_N^{\ast}(x)}{h}.
\end{equation}
Here we took into account that the optimal control is continuously differentiable in $\Omega_{AB}$ and $b$ and 
$\sigma$ are smooth. Therefore the drift and the diffusion in \eqref{eq:sde1c} are finite and hence the probability that $\tau_{AB}<h$ tends to zero as $h\rightarrow 0$.
Furthermore, we note that 
\begin{equation}
\label{ineq22}
  \lim_{h\rightarrow 0} \frac{ \mathbb{E}_{P,x}\left[c_N^{\ast}(Y_h) \right] - c_N^{\ast}(x)}{h} =
 \lim_{h\rightarrow 0} \frac{ \mathbb{E}_{P_v,x}\left[c_N^{\ast}(Y_h) \right] - c_N^{\ast}(x)}{h} = \mathcal{L}_vc^*_N(x),
\end{equation}
where $\mathcal{L}_v$ is the generator of the controlled process \eqref{eq:sde1c}.
The equality \eqref{ineq22} follows from the fact that 
$$
 \mathbb{E}_{P_v,x}\left[c_N^{\ast}(Y_h) \right] =  \mathbb{E}_{P,x}\left[\frac{dP_v}{dP}c_N^{\ast}(Y_h) \right] ,
 $$
 where 
 $$
 \frac{dP_v}{dP} = \exp\left\{-\int_0^{h}\sigma^\top v(Y_s)\cdot dW_s - \frac{1}{2}\int_0^{h}\|\sigma^\top v(Y_s)\|^2ds\right\}\rightarrow 1~~{\rm a.s.}~~{\rm as}~~h\rightarrow 0.
 $$

Therefore, \eqref{ineq2} is equivalent to
\begin{equation}
\label{ineq3}
 \frac{1}{2}\|\sigma^\top v\|^2 + \left[b+\sigma\sigma^\top v\right]\cdot\nabla c^*_N +\frac{1}{2}{\sf tr}\left(\sigma\sigma^\top \nabla\nabla c_N^*\right) \ge 0.
 \end{equation}
 Furthermore, the equality is reached if and only if the control $v$ is optimal, i.e.,
 \begin{equation}
\label{HJB0}
\inf_{v\in\mathcal{V}} \left[\frac{1}{2}\|\sigma^\top v\|^2 + \sigma\sigma^\top v \cdot\nabla c^*_N\right] + 
b \cdot\nabla c^*_N+\frac{1}{2}{\sf tr}\left(\sigma\sigma^\top \nabla\nabla c_N^*\right) = 0.
 \end{equation}
  The function in the square brackets in \eqref{HJB0} is convex quadratic in $v$. To minimize it, we take its gradient and set it to zero:
   \begin{equation}
   \label{minv}
         \nabla_v\left[\sigma\sigma^{\top} v \cdot \nabla c_N^*+ \frac{1}{2} v^{\top} \sigma\sigma^{\top} v\right] = \sigma\sigma^{\top} \nabla c_N^*+ \sigma\sigma^{\top} v = 0.
         \end{equation}
         Therefore, a minimizer $v_N^*$ must satisfy $  \sigma \sigma^{\top} (\nabla c^*_N+ v_N^*)  = 0$. 
         Since columns of $\sigma$ are linearly independent, this condition is equivalent to 
         \begin{equation}
         \label{vstarN}
          \sigma^{\top} (\nabla c_N^* + v_N^*)= 0\quad{\rm or}\quad \sigma^\top v_N^* = - \sigma^{\top}\nabla c_N^*.
          \end{equation}
Plugging this into \eqref{HJB0} we obtain the following equation for the minimal cost $c^*_N$:
    \begin{align}
    0&=    \frac{1}{2} {\sf tr}\left(\sigma \sigma^{\top} \nabla \nabla c_N^*\right) + b \cdot \nabla c_N^*- \sigma\sigma^{\top} \nabla c_N^*\cdot \nabla c_N^*+  \frac{1}{2} (\nabla c_N^*)^{\top} \sigma\sigma^{\top} \nabla c_N^*\notag \\
    & =   \frac{1}{2} {\sf tr}\left(\sigma \sigma^{\top} \nabla \nabla c_N^*\right) + b \cdot \nabla c_N^*-  \frac{1}{2} (\nabla c_N^*)^{\top} \sigma\sigma^{\top} \nabla c_N^*. \label{HJB1}
    \end{align}

{\bf Step 4. Show that $c^*_N = -\log q_N^+(x)$ is the solution to the HJB equation.}
Plugging $c_N^* =-\log q_N^+$ into  \eqref{HJB1} we get
      \begin{align*}
     - &\frac{1}{2q_N^+}  {\sf tr}\left(\sigma\sigma^{\top}\nabla\nabla q_N^+\right)+ \frac{1}{2(q_N^+)^2} \left\|\sigma^{\top}  \nabla q_N^+\right\|^2 
     - b \cdot \nabla \log q_N^+-  \frac{1}{2(q_N^+)^2}\left\|\sigma^{\top}  \nabla q_N^+\right\|^2 =\\
        - &\frac{1}{2q_N^+}{\sf tr}\left( \sigma\sigma^{\top}\nabla\nabla q_N^+\right) - b \cdot \nabla \log q_N^+ =\\
      -& \frac{1}{q_N^+} \mathcal{L}q_N^+ = 0. 
    \end{align*}
   The last equality follows from the fact that $ \mathcal{L}q_N^+ = 0$ in $\Omega_{AB}$. 
The boundary conditions for $c_N^* =-\log q_N^+$  are readily checked: $-\log q_N^+ = N$ on $\partial A$, $-\log q_N^+ = 0$  on $ \partial B$, and 
$$
 \frac{\partial }{\partial\hat{n}} \left(-\log q^+_N\right) = -\frac{1}{q^+_N} \frac{\partial q_N^+}{\partial\hat{n}}  = 0,\quad x\in\partial\Omega.
 $$
The optimal control associated with $c^*_N = -\log q_N^+(x)$ given by \eqref{vstarN} is
\begin{equation}
\label{vstarN1}
\sigma^\top v_N^* = - \sigma^{\top}\nabla c_N^* = \sigma^{\top}\nabla \log q^+_N.
\end{equation}

{\bf Step 5. Show that the control $\sigma^\top v^*_N=\sigma^\top\log q^+_N$ is admissible.}
Equation \eqref{commNcomm} implies that
\begin{equation}
\label{optvN}
\sigma^\top v^*_N = 
-\sigma^\top \nabla c^*_N = \sigma^\top\nabla\log q_N^+ =
 \sigma^\top\left[ \frac{\nabla q^+\left(1-e^{-N}\right)}{ q^+\left(1-e^{-N}\right) + e^{-N}} \right].
\end{equation}
Hence 
\begin{equation}
\label{bvN}
\left\|\sigma^\top v^*_N\right\| \le e^N\max_{x\in\Omega_{AB}}\left\|\nabla q^+_N(x)\right\| < \infty.
\end{equation}
The stopping time $\tau_{AB}<\infty$ a.s. as the system is ergodic and the domain is compact. 
Therefore
$$
\frac{1}{2}\int_0^{\tau_{AB}}\|\sigma^\top v^*_N(Y_s)\|^2ds <\infty~~{\rm a.s.}
$$
and hence
\begin{equation}
\label{adm}
\mathbb{E}_{P,x}\left[e^{\frac{1}{2}\int_0^{\tau_{AB}}\|\sigma^\top v^*_N(Y_s)\|^2ds}\right] < \infty,
\end{equation}
i.e. $v_N^*$ is admissible.

%%%
{\bf Step 6. Take the limit $N\rightarrow\infty$.}
Letting $N\rightarrow\infty$ in \eqref{RN5}  and taking into account the explicit expression  \eqref{commNcomm} for $q^+_N(x)$ we conclude that
\begin{equation}
\label{lowerbound}
c^*(x) \ge -\log q^+(x).
\end{equation}
On the other hand, as we have shown in Step 3,  \eqref{RN5} is actually an equality, and the corresponding optimal control satisfies
\begin{equation}
\label{optvN}
\sigma^\top v^*_N = 
-\sigma^\top \nabla c^*_N = \sigma^\top\nabla\log q_N^+ =
 \sigma^\top\left[ \frac{\nabla q^+\left(1-e^{-N}\right)}{ q^+\left(1-e^{-N}\right) + e^{-N}} \right].
\end{equation}
Taking limit $N\rightarrow\infty$ we obtain
\begin{equation}
\label{vstar1}
\sigma^\top v^* = \sigma^\top \frac{\nabla q^+}{q^+} = \sigma^\top\nabla\log q^+.
\end{equation}
Since the admissible set $\mathcal{V}$ is closed, $ v^* =\nabla\log q^+ \in\mathcal{V}$.
One can readily check that the corresponding solution the Hamilton-Jacobi-Bellman equation \eqref{HJB0} with the boundary conditions $c^* = +\infty$ on $\partial A$, $c^* = 0$ on $\partial B$, and $\tfrac{\partial c^*}{\partial\hat{n}} = 0$ on $\partial \Omega$ is 
\begin{equation}
\label{cstar}
c^*(x) = -\log q^+(x).
\end{equation}
This completes the proof of Theorem \eqref{thm:T1}.
\end{proof}

%%%%%%%%%%%%%%%%%%%%%%%%%%%%%%%%
\section{Errors due to model reduction: an example}
\label{app:LLexample}
We will illustrate the error due to model reduction in the transition rate 
$\nu_{AB}$ as well as in the probabilities $\rho_A$ and $\rho_{AB}$ on the example used in \cite{LL2010}.  
A system is evolving according to the overdamped Langevin dynamics \eqref{sde2} with the potential given by
\begin{equation}
\label{LLpot}
V(x,y) = (x^2-1)^2 + \epsilon^{-1}(y+x^2-1)^2,
\end{equation}
where $\epsilon$ is a small parameter.
The second term in \eqref{LLpot} effectively restricts the dynamics to a small neighborhood of the parabola $y = 1-x^2$. 
It is shown in \cite{LL2010} that $x$ is a suboptimal choice of a collective variable because the gradient of 
$x$ with respect to $(x,y)$ is not orthogonal to the normal vector to the manifold near which the dynamics live.  

\begin{figure}[htbp]
\begin{center}
\includegraphics[width = 0.9\textwidth]{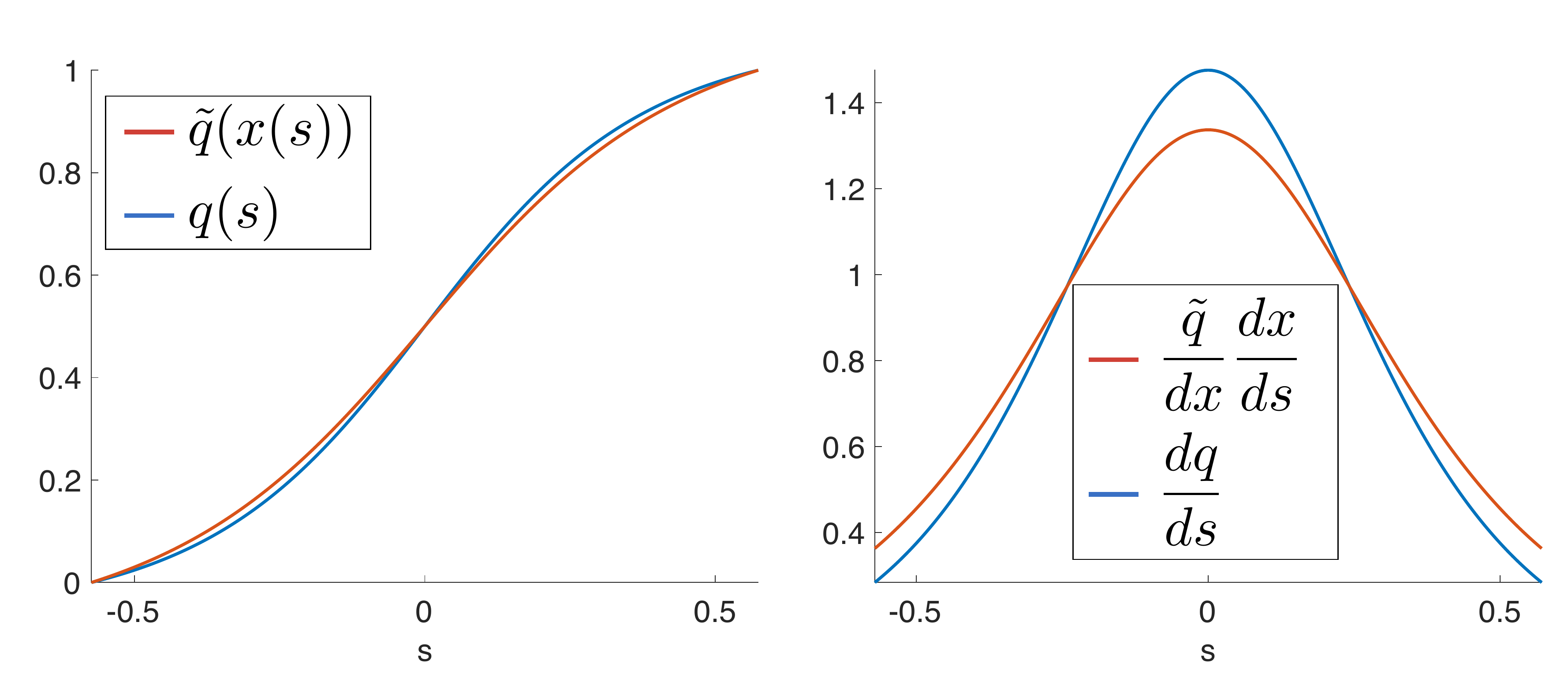}
\caption{An illustration to \ref{app:LLexample}.}
\label{fig:LL}
\end{center}
\end{figure}

Let us consider the signed arclength parameter along the parabola $y = 1-x^2$
\begin{equation}
\label{sx}
s(x) = \int_0^x \sqrt{1 + 4z^2}dz = \frac{1}{2}x\sqrt{1+4x^2} + \frac{1}{4}\log(2x + \sqrt{1+4x^2})
\end{equation}
as a collective variable. The function $s(x)$ is monotone and hence invertible.
In the limit $\epsilon\rightarrow 0$, the dynamics are one-dimensional and governed by
\begin{equation}
\label{dyn1}
ds = -\frac{d V_0(x(s))}{ds}dt + \sqrt{2\beta^{-1}}dw,
\end{equation}
where $V_0:=(x^2-1)^2$. We set $\beta = 3$ choose the sets $A$ and $B$ as in \cite{LL2010}:
\begin{equation}
\label{LL_AB}
A = \{x < a\},\quad B = \{x >b\},\quad a= -0.5,\quad b = 0.5.
\end{equation}
We calculate the committor $q(s)$ and $\tilde{q}(x)$ using the exact formula for the one-dimensional case:
\begin{equation}
\label{qLL}
q(s) = \frac{\int_{s(a)}^{s}e^{\beta V_0(x(s'))}ds'}{\int_{s(a)}^{s(b)}e^{\beta V_0(x(s'))}ds'},\quad
\tilde{q}(x) = \frac{\int_{a}^{x}e^{\beta V_0(x')}dx'}{\int_{a}^{b}e^{\beta V_0(x')}dx'}.
\end{equation}
Using the inverse of $s(x)$, we obtain $\tilde{q}(x(s))$. The plots of $q(s)$ and $\tilde{q}(x(s))$ and their derivatives in $s$ are displayed in Fig. \ref{fig:LL}(left). It is evident that the difference between their derivatives is notably larger. 
% { AS: the figure doesn't compile. Also, you might want to fix it in place}
Next, we use $q(s)$ and $\tilde{q}(x)$ to calculate the transition rate { $\nu_{AB}$ from $A$ to $B$ via \eqref{nu1} and the probability $\rho_{AB}$ via \eqref{rhoAB}}. 
The notation with tilde will indicate the results obtained using $x$ as a collective variable. We get:
\begin{align*}
\nu_{AB} &= 9.46\cdot10^{-3},& \quad &\rho_{AB} = 9.47\cdot10^{-3},\\
\tilde{\nu}_{AB}& = 2.19\cdot10^{-2},&\quad  &\tilde{\rho}_{AB} = 7.98\cdot10^{-3}.
\end{align*}
The transition rate estimated using $x$ as a collective variable exceeds to true rate by the factor of approximately 2.3, while the error in the estimate of the probability to be reactive is about 16\%.
%%%%%%%%%%%%%%%%%%%%%%%%%%%%%%%%
\section{Robustness of the crossover time: an example.}
\label{app:E}
{
In this appendix, we will examine how the quality of the committor estimate affects the crossover time in the controlled dynamics on the example of the bistable Duffing oscillator \eqref{Duffing} with $\epsilon = 0.05$. The controlled dynamics of the Duffing oscillator are governed by SDE \eqref{eq:Duffing_controlled} where $q^+_{\sf nn}$ is the estimate to the forward committor computed using PINN.  
%through example of the bistable Duffing oscillators at $\beta = 20$. 

The PINN committor solver computes the committor via an optimization process called \emph{training the neural network}. The number of training steps is measured in epochs\footnotemark[2].
\footnotetext[2]{One epoch comprises the number of iterations necessary to use all training data $N_{\sf train}$ one time. For example, if the optimization algorithm is deterministic and all training data are used for computing the gradient of the objective function in each iteration, then one epoch is equal to one iteration. If an optimizer is stochastic and a subset, a \emph{batch}, of $N_{\sf batch}$ training data is used at each iteration to evaluate the direction of the step, then one epoch consists of ${\sf round}(N_{\sf train}/N_{\sf batch})$ iterations.
}
Stopping the training process too early results in a rough approximation to the forward committor. The results reported in Section \ref{sec:Duffing} are obtained as a result of training the neural network for 500 epochs. Hence, we form a sequence of approximations to the forward committor by evaluating the solution model after 100, 125, 150, 300, and 500 epochs of training visualized in the figures in Table \ref{table:Duffing_tau} using dashed contour plots. The solution becomes progressively closer to the final solution $q_{\sf nn}^+$ evaluated at 500 epochs, which, in turn, in close to the FEM solution depicted using solid contour plots. The discrepancies MAD and RMSD between the FEM forward committor and the approximations to it progressively shrink.

For each of these ``undertrained" solutions, we evaluate the expected crossover time $\mathbb{E}[\tau_{AB}]$ by averaging the crossover times of 250 transition trajectories governed by the controlled SDE \eqref{eq:Duffing_controlled} with the corresponding ``undertrained" forward committor. The results are shown in the last column of Table \ref{table:Duffing_tau}.
}

% Therefore, we obtain a sequence of a
% In particular, we took five snapshots of the PINN model during the training of the committor. Case 0 is the final PINN model after 500 epochs of training, also included in Table \ref{table:Duffing3}; Case 1 to Case 4 are PINN models saved after 300, 150, 125 and 100 epochs respectively.  
% The qualities of the committor estimates are assessed via visual comparison with FEM results and numerical comparison using wMAD and wRMSD. 
% The influence on the controlled dynamics is examined through the expected crossover time $\mathbb{E}[\tau_{AB}]$ and the results are reported in Table \ref{table:Duffing_tau}. 
\begin{table}[htbp]
    \centering
    \begin{tabular}{|M{2.1cm}|M{1.8cm}|M{1.8cm}|M{5cm}|M{2.5cm}|}
    \hline
          & wMAD & wRMSD & Visual comparison of $q^+$ of PINNs and FEM &$\mathbb{E}[\tau_{AB}]$ \\
    \hline
        Case 0 (epoch 500) & 1.3e-2 & 2.0e-2 & \includegraphics[width=0.3\textwidth]{forward_beta20_ellip.pdf} & 7.34 $\pm$ 0.33\\
    \hline
        Case 1 (epoch 300)& 4.9e-2 & 6.1e-2 &\includegraphics[width=0.3\textwidth]{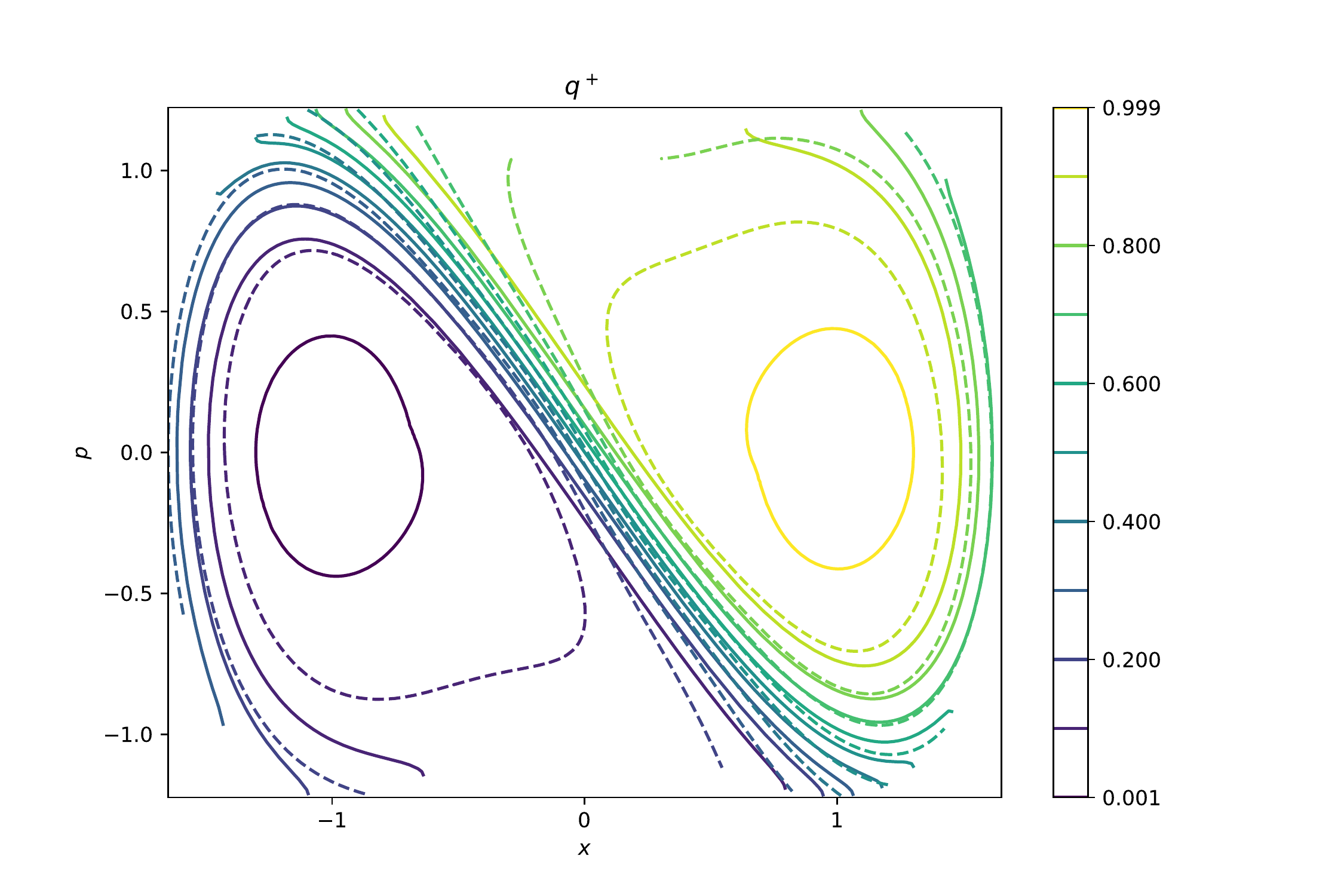}  & 7.88 $\pm$ 0.36\\
    \hline
        Case 2 (epoch 150) & 6.0e-2& 7.7e-2 & \includegraphics[width=0.3\textwidth]{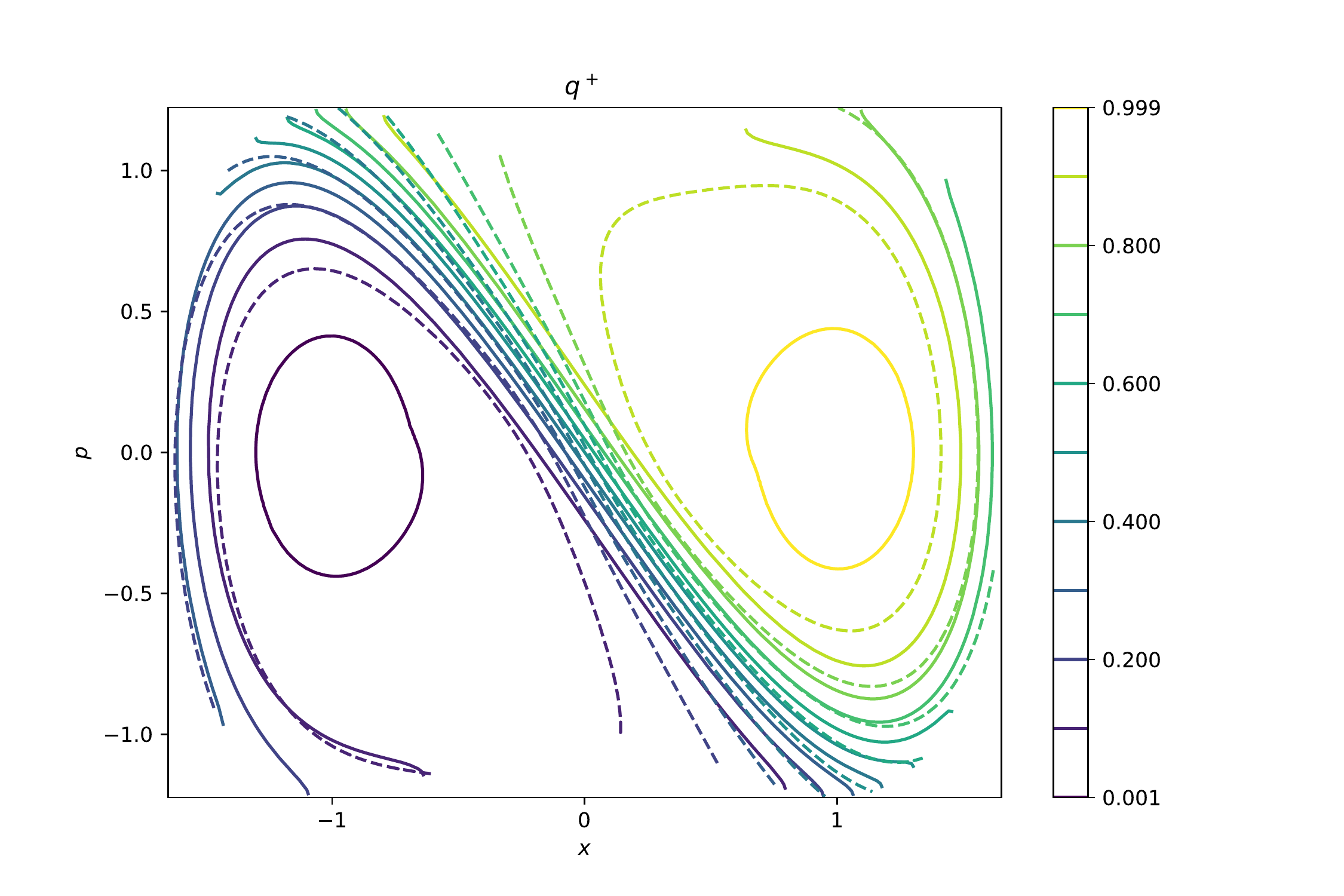} &8.48 $\pm$ 0.84\\
    \hline
    Case 3 (epoch 125)& 7.9e-2& 10.2e-2 & \includegraphics[width=0.3\textwidth]{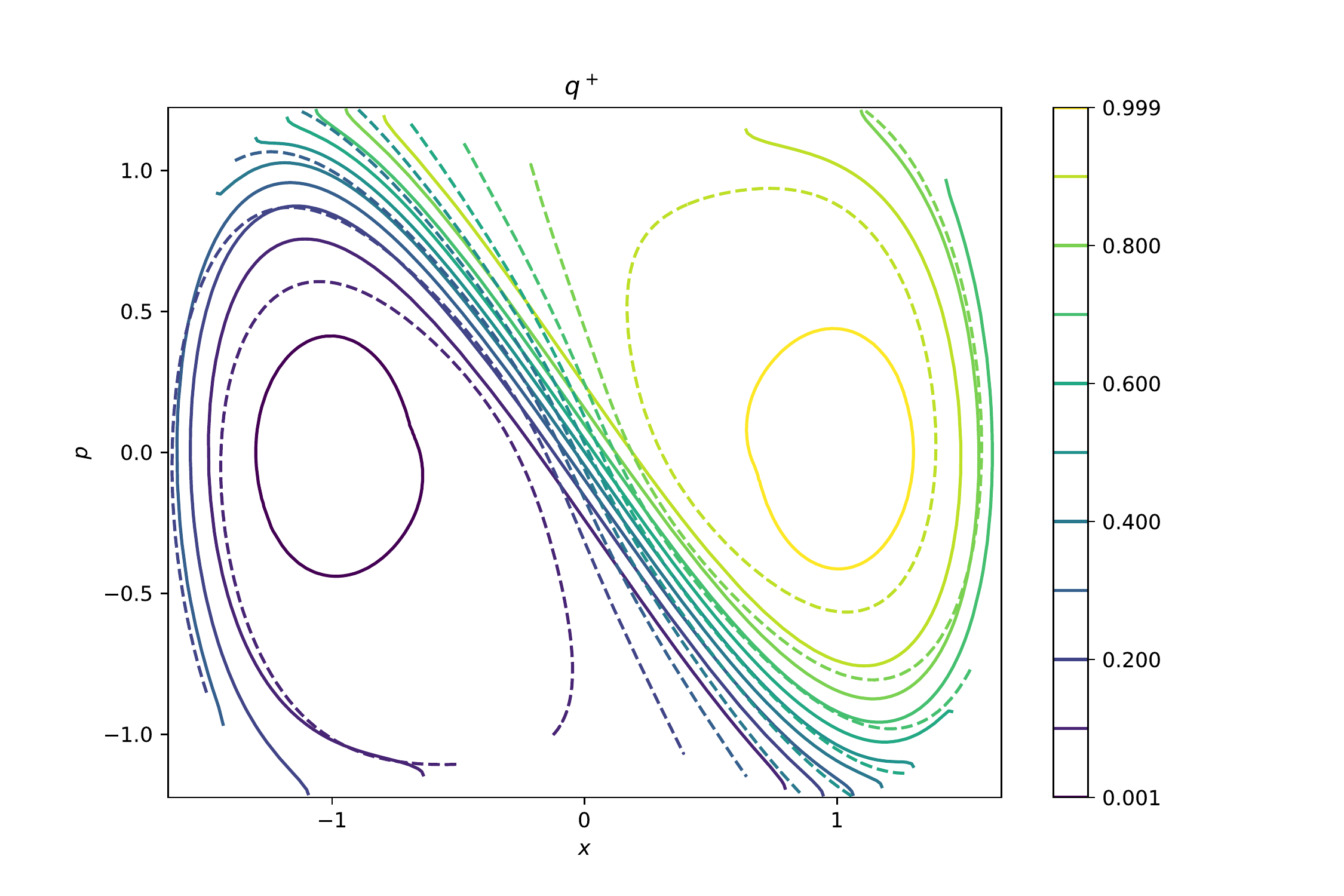} &12.18 $\pm$ 2.49\\
    \hline
    Case 4 (epoch 100) & 13.6e-2& 16.7e-2 & \includegraphics[width=0.3\textwidth]{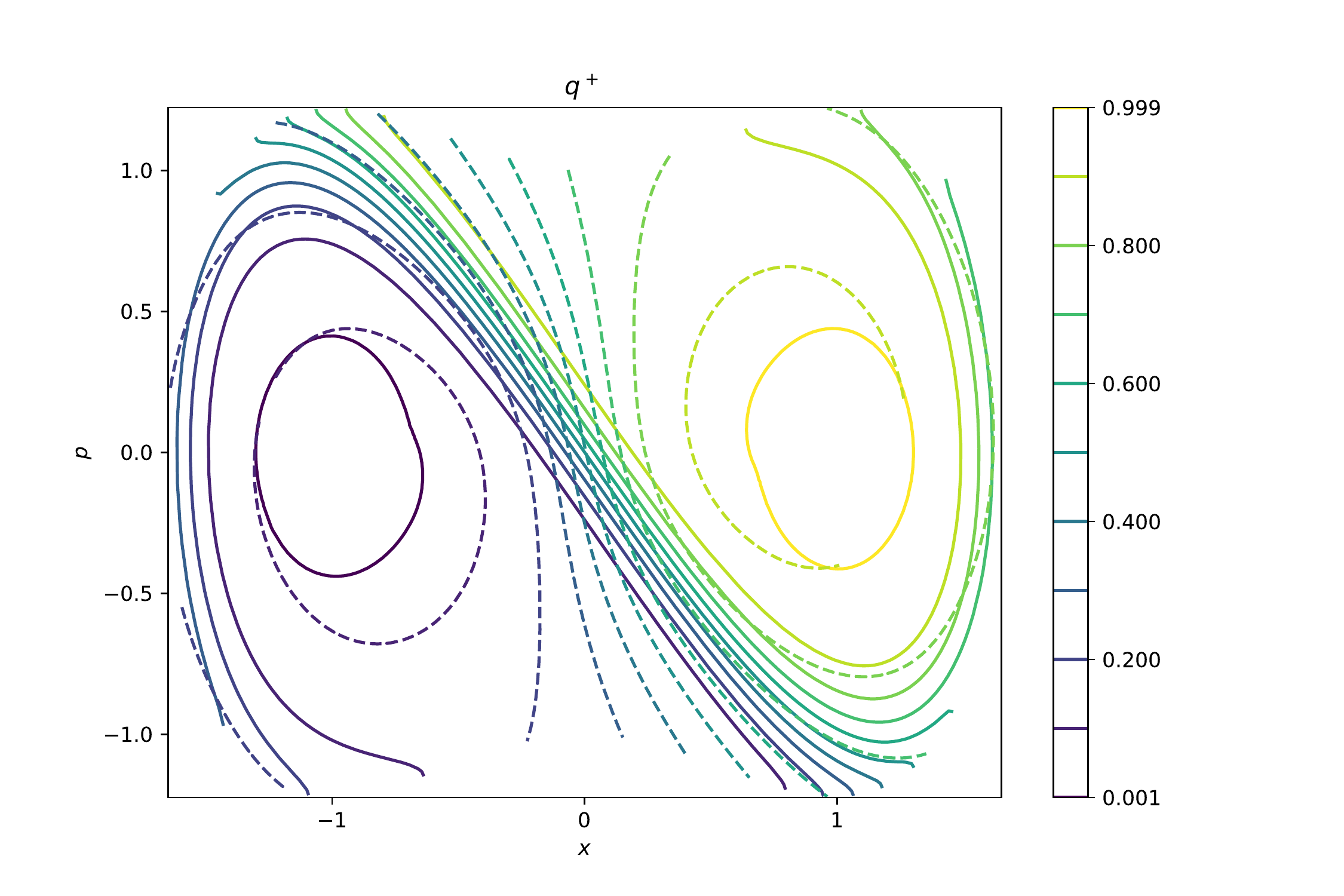} FEM(solid), PINN(dashed) &48.86 $\pm$ 5.71\\
    \hline
    \end{tabular}
    % used to be 0.35\textwidth
    \caption{Estimates of wMAD, wRMSD and $\mathbb{E}[\tau_{AB}]$ computed for five different approximations to the forward committor (dashed contours) of the bistable Duffing oscillator \eqref{Duffing} with $\epsilon = 0.05$ obtained at various stages of training the neural network-based solution model with the PINN loss function. The FEM forward committor (solid contours) and the corresponding $\mathbb{E}[\tau_{AB}]=\mathbf{7.48 \pm 0.49}$ are treated as the ground truth. 
    % $\mathbb{E}[\tau_{AB}]=\mathbf{7.48 \pm 0.49}$ from brute force simulation is $\mathbf{7.48 \pm 0.49}$. Visual comparison shows contours of the forward committor computed using FEM and PINNs, marked in solid and dashed lines respectively. 
    }
    \label{table:Duffing_tau}
\end{table}

\begin{figure}[htbp]
\begin{center}
\includegraphics[width = 0.9\textwidth]{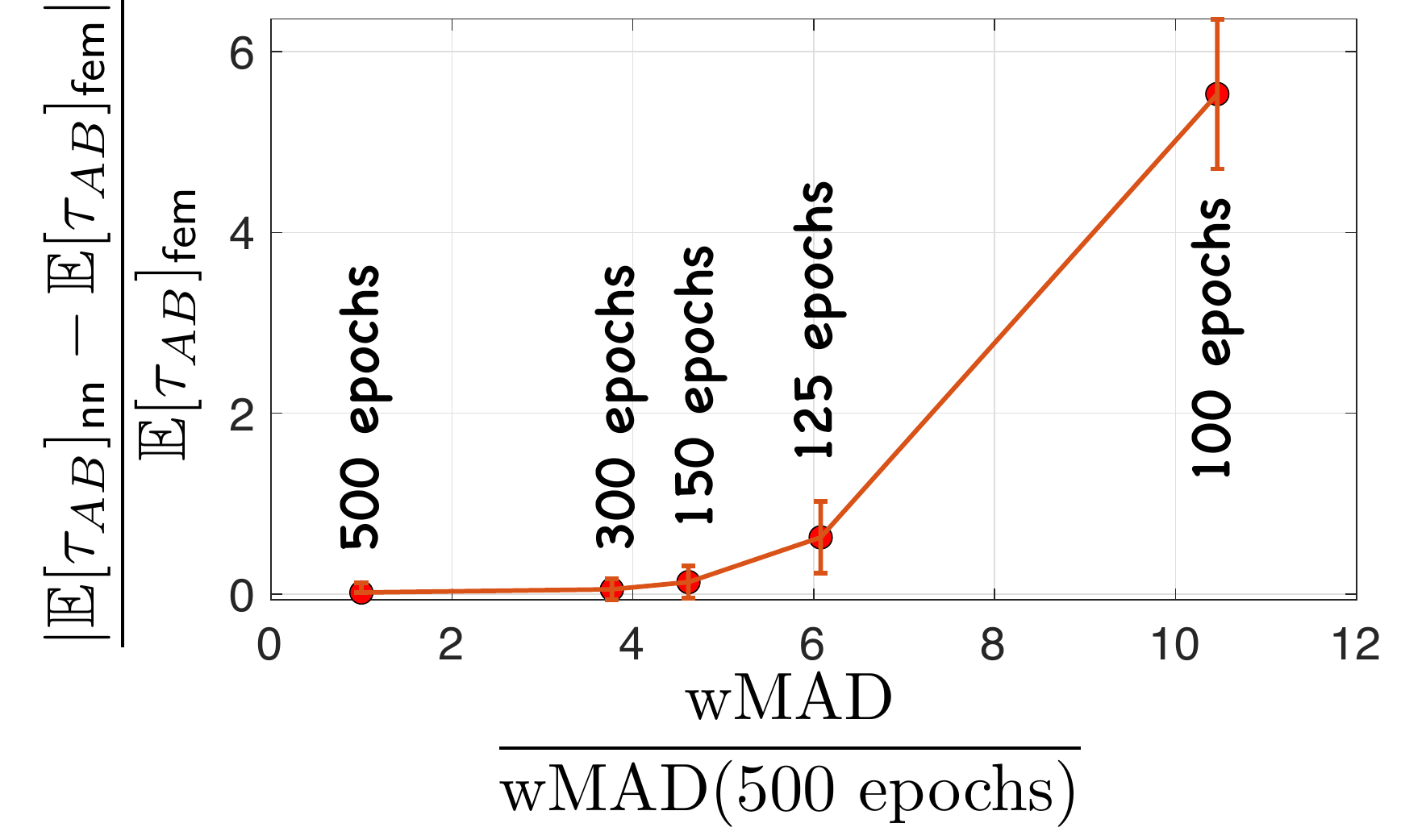}
\caption{Test case: the bistable Duffing oscillator with $\epsilon = 0.05$. The growth of the relative error in the expected crossover time $\mathbb{E}[\tau_{AB}]$ versus the relative growth of the weighted mean absolute difference (wMAD) of the ``undertrained" approximations to the forward committors and the FEM committor. The approximate committors $q_{\sf nn}^+$ are obtained by evaluating the solution model after $100$, $125$, $150$, $300$, and $500$ epochs of training.}
\label{fig:EtauAB}
\end{center}
\end{figure}

{
Table \ref{table:Duffing_tau} shows that as the accuracy of the committor decreases, the expected crossover time $\mathbb{E}[\tau_{AB}]$ increases. Nonetheless, the relative increment in $\mathbb{E}[\tau_{AB}]$ is notably smaller compared to the discrepancy in the committors as evident from Figure \ref{fig:EtauAB}.
}
{
In summary, this investigation suggests that the estimate of the expected crossover time obtained by means of sampling controlled trajectories remains reasonably accurate even if the approximation to the forward committor is rough.
}
% Nonetheless, the increase in $\mathbb{E}[\tau_{AB}]$ is comparably smaller than the discrepancy in the committors.
% not noticeable so long as the approximation of the committor stays at a reasonable level. 
% Compared to Case 0, PINN-based committor in Case 1 has wRMSD approximately 3 times of the wRMSD in Case 0, $\mathbb{E}[\tau_{AB}]$ on the other hand increases by 7.4\%. 
% In Case 2, committor estimates of PINN model has wRMSD 3.9 times of the wRMSD in Case 0, the disparity of $\mathbb{E}[\tau_{AB}]$ is about 15.5\%. 
% Case 3 has larger discrepancy in the estimated mean crossover time. 
% The magnitude of the increase however, is still significantly smaller than that of the committor function. 
% Meanwhile, shown in Case 4, as the error of estimates for committor function surpasses a certain level, the sampled crossover time under the controlled process becomes almost useless.
% From this investigation, we also observe an increase in variance of the simulated $\mathbb{E}[\tau_{AB}]$ as the accuracy of the committor decreases. 

%%%%%%%%%%%%%%%%%%%%%%%%%%%%%%%%
\section{Finite element method for the committor problem}
\label{sec:appFEM}
\subsection{Time-reversible dynamics}
\label{sec:appFEM1}
If the governing SDE is time-reversible as it is in the case of the overdamped Langevin dynamics \eqref{sde2} or the overdamped Langevin dynamics in collective variables \eqref{sde4}, the committor problem  \eqref{eg:commproblem} with the generator \eqref{gen4} is self-adjoint. In this case, we proceed in the standard way detailed in \cite{50linesFEM}. First, we decompose the committor $q$ into $q = q_1 + q_0$ where is a prescribed function such that $q_1 = 1$ on $\partial B$ and $q_1 = 0$ outside a small neighborhood of $\partial B$ and $q_0$ needs to be found. 
The boundary value problem for $q_0$ is 
\begin{align}
\mathcal{L}q_0 &= -\mathcal{L}q_1,\quad x\in\Omega_{AB},\label{bvpq0}\\
q_0 &= 0,\quad x\in\partial A\cup\partial B,\label{bvpq0_d}\\
\frac{\partial q_0}{\partial \hat{n}} & = 0,\quad x\in\partial \Omega.\label{bvpq0_n}
\end{align}
Second, we multiply multiply \eqref{bvpq0} by a test function $w\in H^1_0(\bar{\Omega})$ where the subscript $0$ means that $w=0$ on $\partial A\cup\partial B$, integrate over $\Omega_{AB}$ and apply the generalized divergence theorem to both parts.
The result is the following integral equation for $q_0$ that must hold for all $w\in C^1_0(\bar{\Omega})$:
\begin{equation}
\label{FEM1}
\int_{\Omega_{AB}} e^{-\beta F(x)}w(x)^\top M(x) \nabla q_0 =  -\int_{\Omega_{AB}} e^{-\beta F(x)}w(x)^\top M(x) \nabla q_1.
\end{equation}
Next, we triangulate $\Omega_{AB}$ and denote the associated finite element space by $S_h$ with the standard piecewise-linear basis 
$\{v_j(x)\}_{i\in\mathcal{I}}$ where $\mathcal{I}$ is the set of vertices of the triangles \cite{50linesFEM}. The subset of vertices that do not belong to $\partial A\cup \partial B$ is denoted by $\mathcal{I}_{\rm free}$. 
We choose $q_1\in S_h$ so that $q_1 = 1$ only at the nodes lying on $\partial B$ and $q_1 = 0$ at all other nodes.  
We seek the finite element solution for $q_0$ of the form 
\begin{equation}
\label{FEM2}
q_0 = \sum_{k\in\mathcal{I}_{\rm free}}v_k(x)(q_0)_k
\end{equation}
where the vector $\{(q_0)_k\}_{k\in\mathcal{I}_{free}}$ is the solution to the linear system
\begin{equation}
\label{FEM3}
\sum_{k\in\mathcal{I}_{\rm free}}{\sf A}_{jk}(q_0)_k =  -\sum_{k\in\partial B}{\sf A}_{jk}
\end{equation}
with the matrix elements ${\sf A}_{jk}$ given by
\begin{equation}
\label{FEM4}
{\sf A}_{jk} = \int_{\Omega_{AB}}e^{-\beta F(x)}\nabla v_j(x)^\top M(x)\nabla v_k(x)dx.
\end{equation}
The integral in \eqref{FEM4} is the sum of the integrals over all triangles.
In each triangle, the gradients of the basis functions are constant, and $F(x)$ and $M(x)$ are approximated by their values at the center of mass.
Finally, the finite element solution $q_{\sf fem}$ is found at the sum $q_{\sf fem} = q_0 + q_1$.

%%%%%%%%%%%%%
\subsection{The Langevin dynamics}
\label{sec:appFEM2}
For the Langevin dynamics \eqref{sde5}, the committor problem \eqref{eg:commproblem} 
with the generator \eqref{gen5} is hypoelliptic, and the application of the finite element method (FEM) requires care.
We design a FEM solver for this case motivated by the article by Morton on FEM for non-self-adjoint problems \cite{FEMmorton} that suggests to make the problem as close to self-adjoint as possible. Since in the case of Langevin dynamics FEM is practical only if the space $(x,p)$ is two-dimensional, $x$ and $p$ will be one-dimensional in the presentation below.

As in \ref{sec:appFEM1}, we start by decompositing $q^{+}$ into $q^{+} = q^{+}_0 + q^{+}_1$. The boundary value problem for $q^{+}_0$ is of the form \eqref{bvpq0}--\eqref{bvpq0_n} except for $x$ is replaced with $(x,p)$. 
Then we multiply the PDE for $q^{+}_0$ by $\epsilon^{-1}\exp\left(-\tfrac{p^2}{2\epsilon}\right)$ and get:
\begin{equation}
\label{FEM5}
\frac{e^{-\frac{p^2}{2\epsilon}}}{\epsilon} \mathcal{L}q^{+}_0 = 
\frac{e^{-\frac{p^2}{2\epsilon}}}{\epsilon}\left[\frac{p}{m}\frac{d q_0}{dx} - V'(x)\frac{dq^{+}_0}{dp}\right] + 
\gamma\frac{d}{dp}\left( e^{-\frac{p^2}{2\epsilon}}\frac{dq^{+}_0}{dp}\right)  = -\frac{e^{-\frac{p^2}{2\epsilon}}}{\epsilon}\mathcal{L}q^{+}_1. 
\end{equation}
We denote $(x,p)$ by $z$ and the gradient with respect to $z$ by $\nabla$ and rewrite \eqref{FEM5} in a matrix form:
\begin{equation}
\label{FEM6}
\gamma\nabla \cdot\left(e^{-\frac{p^2}{2\epsilon}}\left[\begin{array}{cc}0&0\\0&m\end{array}\right]\nabla q^{+}_0\right) +
\frac{e^{-\frac{p^2}{2\epsilon}}}{\epsilon}\left[\begin{array}{c}p/m\\V'(x)\end{array}\right]\cdot\nabla q^{+}_0
= - \frac{e^{-\frac{p^2}{2\epsilon}}}{\epsilon} \mathcal{L}q^{+}_1.
\end{equation}
Then we follow the steps in \ref{sec:appFEM1}. We multiply \eqref{FEM6} by a test function 
$w(z)\in H_0^1(\bar{\Omega}_{AB})$, integrate over $\Omega_{AB}$ and apply the generalized divergence theorem.
This results in the following integral equation for $q^{+}_0$ 
that must hold for all $w(z)\in H_0^1(\bar{\Omega}_{AB})$:
\begin{align}
&\gamma\int_{\Omega_{AB}}e^{-\frac{p^2}{2\epsilon}}\nabla w^\top \left[\begin{array}{cc}0&0\\0&m\end{array}\right]\nabla q^{+}_0dx 
- \frac{e^{-\frac{p^2}{2\epsilon}}}{\epsilon}\int_{\Omega_{AB}} w \left[\begin{array}{c}p/m\\V'(x)\end{array}\right]\cdot\nabla q^{+}_0  \notag \\ 
=& 
-\gamma\int_{\Omega_{AB}}e^{-\frac{p^2}{2\epsilon}}\nabla w^\top \left[\begin{array}{cc}0&0\\0&m\end{array}\right]\nabla q^{+}_1dx 
+ \frac{e^{-\frac{p^2}{2\epsilon}}}{\epsilon}\int_{\Omega_{AB}} w \left[\begin{array}{c}p/m\\V'(x)\end{array}\right]\cdot\nabla q^{+}_1dx.\label{FEM7}
\end{align}
Then we triangulate $\Omega_{AB}$, introduce the standard FEM basis $\{v_i\}_{i\in\mathcal{I}}$ in $S_h$, and represent $q^{+}_0(z)$ as a linear combination of the basis functions associated with the nodes not in $\partial A\cup\partial B$, and obtain the following linear system for the coefficients $(q^{+}_0)_k$:
\begin{equation}
\label{FEM7}
\sum_{k\in\mathcal{I}_{\rm free}}{\sf A}_{jk}(q^{+}_0)_k - \sum_{k\in\mathcal{I}_{\rm free}}{\sf B}_{jk}(q^{+}_0)_k 
= -\sum_{k\in\partial B}{\sf A}_{jk} +  \sum_{k\in\partial B}{\sf B}_{jk}.
\end{equation}
The matrix elements in \eqref{FEM7} are given by
\begin{equation}
\label{FEM8}
{\sf A}_{jk} = \gamma \int_{\Omega_{AB}}e^{-\frac{p^2}{2\epsilon}}\nabla v_j^\top \left[\begin{array}{cc}0&0\\0&m\end{array}\right]\nabla v_k dx,\quad
{\sf B}_{jk} =  \frac{e^{-\frac{p^2}{2\epsilon}}}{\epsilon}\int_{\Omega_{AB}} v_j \left[\begin{array}{c}p/m\\V'(x)\end{array}\right]\cdot\nabla v_k dx.
\end{equation}
Computing the integrals in \eqref{FEM8} over each triangle, all nonlinear functions are approximated by 
their values at the centers of mass of the triangle.
Finally, $q^+_{\sf fem} = q^+_0 + q^+_1$. The backward committor is readily found by 
$q^-_{\sf fem}(x,p) = 1 - q^+_{\sf fem}(x,-p)$.

\section{Confidence interval}
\label{app:confidence_interval}
To compute confidence intervals for simulated transition time from $A$ to $B$, we first compute mean $\Bar{\tau}_{AB}$ and standard error of the mean $\text{se}_{\tau}$ using \texttt{scipy.stats.sem}, which is the sample standard deviation divided by square root of the sample size: $\text{se}_{\tau} = \frac{\sigma}{\sqrt{n}}$. The confidence interval then is obtained using t distribution:
\begin{equation}
    \text{confidence interval} = [\Bar{\tau}_{AB} - \text{se}_{\tau} t(\alpha), \Bar{\tau}_{AB} + \text{se}_{\tau} t(\alpha)]
\end{equation}
where $t(\alpha)$ satisfies $\int_{t(\alpha)}^\infty f(x)dx = \alpha$, and can be found using \texttt{scipy.stats.t.ppf}. For all examples, $95\%$ confidence intervals are used.

 \bibliographystyle{elsarticle-num} 
% \bibliography{reference}


%\begin{thebibliography}{00}


%

%\end{thebibliography}
\end{document}